\documentclass[11pt]{article}
\usepackage[OT1]{fontenc}

\usepackage{amsmath,amssymb,amsfonts,amsthm}
\usepackage{mathtools,stmaryrd}
\usepackage{array,booktabs,longtable,tabularx}
\usepackage[shortlabels]{enumitem}
\usepackage{xcolor}
\usepackage[margin=1in]{geometry}
\usepackage[expansion=false]{microtype}
\usepackage{hyperref}
\hypersetup{colorlinks=true,linkcolor=black,citecolor=red,urlcolor=red}
\numberwithin{equation}{section}

\newcommand{\defterm}[1]{\textbf{\textup{#1}}}
\newcommand{\R}{\mathbb R}
\newcommand{\C}{\mathbb C}
\newcommand{\HH}{\mathcal H}
\newcommand{\dd}{\mathrm d}
\newcommand{\vol}{\operatorname{vol}}
\newcommand{\tr}{\operatorname{tr}}

\newcommand{\reg}{\operatorname{reg}}
\newcommand{\sing}{\operatorname{sing}}
\newcommand{\supp}{\operatorname{supp}}
\newcommand{\dist}{\operatorname{dist}}
\renewcommand{\llbracket}{[\![}
\renewcommand{\rrbracket}{]\!]}

\newtheorem{theorem}{Theorem}[section]
\newtheorem*{introbranchtheorem}{Theorem A [\ref*{thm:intro-branch}]}
\newtheorem*{flatnessrestatement}{Theorem~\ref*{thm:intro-flatness}}
\newtheorem*{continuationrestatement}{Theorem~\ref*{thm:intro-continuation}}
\newtheoremstyle{introlettered}
  {\topsep}{\topsep}{\itshape}{}{\bfseries}{.}{.5em}
  {\thmname{#1}\thmnumber{ [#2]}\thmnote{ (#3)}}
\theoremstyle{introlettered}
\newtheorem{introflatnesstheorem}[theorem]{Theorem B}
\newtheorem{introcontinuationtheorem}[theorem]{Theorem C}
\theoremstyle{plain}
\newtheorem{proposition}[theorem]{Proposition}
\newtheorem{lemma}[theorem]{Lemma}
\newtheorem{corollary}[theorem]{Corollary}
\theoremstyle{definition}
\newtheorem{definition}[theorem]{Definition}
\newtheorem{example}[theorem]{Example}

\newtheorem{remark}[theorem]{Remark}

\title{Branch Obstructions for Area-Minimizing Currents\\
and Flat Critical Points of {\color{black}Weakly-Conformal} Maps}
\author{Anton Iliashenko\\
  {\small Beijing Institute of Mathematical Sciences and Applications}\\
  {\small Beijing, China}\\
  {\small \texttt{antoniliashenko@bimsa.cn}}}
\date{November 2026}

\begin{document}
\maketitle

\begin{abstract}
We study two problems motivated by the local geometry of conformal $\phi$-curves.
First, let $T$ be a $k$-dimensional locally area-minimizing integral current
in Euclidean space, with $k\geq3$, codimension at least two, no boundary,
and a locally conformally flat induced metric on $\reg T$.
We prove the following branch obstruction: at any rapid-decay branch point
whose planar frequency differs from $2$, every nonzero cyclic component
of every cylindrical normalized tangent function has repetition multiplicity
at least two.
Consequently, the nonintegral-frequency points of densities $2$ and $3$
in the rapid-decay branch set have Hausdorff dimension at most $k-3$.
On the full branch set, they have zero $\mathcal H^{k-2}$ measure.
These conclusions apply to Euclidean image currents of conformal $\phi$-curves
with locally conformally flat domains, provided the image currents satisfy
the stated hypotheses.

Second, for $C^j$ weakly-conformal maps between smooth Riemannian manifolds,
with domain dimension $k\geq3$ and $j\geq2$, we prove that all derivatives
of orders $1$ through $j$ vanish at every critical point.
As a consequence, for smooth weakly-conformal maps with connected domain, we establish
a curvature condition under which the map is constant or has no critical points.
The flatness and continuation results apply in particular to conformal
$\phi$-curves with the corresponding regularity and curvature assumptions.
Finally, we examine the Bor\r{u}vka cone, which illustrates the role of regularity and the distinction
between critical points and branch points.
\end{abstract}\tableofcontents

\section{Introduction}\label{sec:intro}
Our starting point is the class of conformal $\phi$-curves:
weakly-conformal maps between Riemannian manifolds that, away from their
critical sets, locally parametrize submanifolds calibrated by a fixed
calibration $\phi$. Our motivation comes from calibrated analogues of
Kobayashi--Royden hyperbolicity, where these maps play the role of
holomorphic disks \cite{BIM}. The equicontinuity and compactness results
of Iliashenko--Karigiannis--Madnick \cite{IKM} make their regularity and
behavior at critical points relevant to that theory. The results here
concern the local geometry underlying these questions.

More precisely, let $(L^k,g)$ and $(M^n,h)$ be Riemannian manifolds,
with $L$ oriented and $M$ carrying a \emph{calibration} $\phi$
of degree $k$, that is, a smooth $k$-form satisfying $\dd\phi=0$
and $\phi(e_1,\ldots,e_k)\leq1$ for every $h$-orthonormal $k$-tuple. A $C^1$ map $u:L\to M$ is a
\emph{conformal $\phi$-curve} if
\[
 u^*h=\lambda^2g,\qquad u^*\phi=\lambda^k\vol_g
\]
for some function $\lambda:L\to[0,\infty)$.
The first equality defines \emph{weak conformality} and implies that
$\lambda=|du|/\sqrt{k}$.
Conversely, with $\lambda=|du|/\sqrt{k}$, the second identity implies the first.
The function $\lambda$ is the \emph{dilation}. Weak conformality implies that, at each $p\in L$,
$du_p$ has rank $k$ or is zero. The \emph{critical set} is
\[
 Z_u:=\{p\in L:du_p=0\}=\{p\in L:\lambda(p)=0\}.
\]
With $\lambda=|du|/\sqrt{k}$, the second equality implies that, away from $Z_u$,
the oriented images of tangent planes are calibrated by $\phi$.

These maps motivate two related questions: which branch points can
occur in their images, and how can the maps degenerate at critical
points? Weak conformality constrains their differentials, while
calibration gives area minimization for suitable image currents.   For smooth maps with  locally conformally flat  domain, the
induced metric on each such submanifold is locally conformally
flat as well. We study the resulting restrictions on branch points
and critical points, when the domain dimension is at least $3$.

We treat the current and map problems separately.   For a locally area-minimizing integral current $T$ in Euclidean space,
we prove that local conformal flatness of the induced metric on $\reg T$
restricts the
multiplicities of cylindrical normalized tangent functions at
rapid-decay branch points. For weakly-conformal maps between arbitrary
smooth Riemannian manifolds, we prove finite-order flatness at critical
points and derive a continuation result under an additional curvature
bound. The map notations and definitions are collected in
Section~\ref{sec:scope}. The current-theoretic definitions are given in
Section~\ref{sec:setup}.

\subsection{Branch obstructions}

Theorem~\ref{thm:intro-branch} is the main result of this part of
the paper. It concerns a locally area-minimizing
integral current $T$ in an open set $\mathcal U\subset\R^{k+m}$,
where $k\geq3$ and $m\geq2$, with $\partial T=0$ in $\mathcal U$.
We first recall the objects in its statement.

An \emph{integral current} records integration over an oriented
rectifiable set with integer multiplicity and locally finite mass, whose
boundary is also integer rectifiable with locally finite mass. Write $\|T\|$ for its mass measure, $\mathbf M$ for mass, and $\partial$ for the boundary operator. The current is \emph{locally area minimizing} if, for every
relatively compact open $W\subset\mathcal U$ and every integral
$(k+1)$-current $S$ with compact support in $W$, we have
\[
 \mathbf M(T\llcorner W)\leq
 \mathbf M((T+\partial S)\llcorner W),
\]
where $\llcorner W$ means restriction to $W$.
The set $\reg T$ consists of points near which
$T=q\llbracket\Sigma\rrbracket$ for a smooth embedded minimal
$k$-submanifold $\Sigma$ through the point and an integer $q\geq1$,
called the \emph{multiplicity}.
The notation $\llbracket\Sigma\rrbracket$ means integration over
$\Sigma$ with multiplicity one and a chosen orientation.
We assume that the induced metric on $\reg T$ is
\emph{locally conformally flat}, meaning that near each point it has
the form $e^{2f}\sum_{i=1}^k(dx^i)^2$ in suitable coordinates,
with $f$ smooth.

A \emph{branch point} is a point of
$\supp T\setminus\reg T$ with a tangent cone $q\llbracket P\rrbracket$,
where $P$ is a $k$-plane and the integer $q\geq2$ is the multiplicity of the tangent cone.
A tangent cone is a limit of magnified currents
$(\eta_{Z,\rho_i})_\#T$, where
$\eta_{Z,\rho}(W):=(W-Z)/\rho$ and $\rho_i\downarrow0$, and $\#$ denotes pushforward. The multiplicity $q$ of this tangent cone equals the density of $T$ at $Z$.
The \emph{rapid-decay branch set} $\mathcal R_q(T)$ consists
of branch points of density $q$ with unique tangent cone
$q\llbracket P_Z\rrbracket$ for which
\[
 \int_{B_r(Z)}\dist^2(W,Z+P_Z)\,\dd\|T\|(W)
 \leq C r^{k+2+2\mu}\qquad(0<r<r_0)
\]
for some $\mu\in(0,1)$, $C<\infty$, and $r_0>0$.

A \emph{normalized tangent function}
$\widehat\Phi:P_Z\to\mathcal A_q(P_Z^\perp)$ is a limit of
$q$-valued approximations to the current after rescaling in the
tangent and normal directions, as specified in
Definition~\ref{def:tangent-function}.
Here $\mathcal A_q(V)$ is the space of unordered $q$-tuples in $V$,
written $\sum_{\ell=1}^q\llbracket v_\ell\rrbracket$, and
$\bigl|\sum_{\ell=1}^q\llbracket v_\ell\rrbracket\bigr|^2
:=\sum_{\ell=1}^q|v_\ell|^2$.
At $Z\in\mathcal R_q(T)$ these functions satisfy
\[
 \widehat\Phi(tX)=t^\nu\widehat\Phi(X)\quad(t>0),\qquad
 \frac1{\omega_k}\int_{B_1^{P_Z}(0)}|\widehat\Phi|^2=1,
\]
where $\omega_k$ is the volume of the unit ball in $\R^k$ and scalar multiplication acts on each entry of the tuple.
The degree $\nu=N_{T,\mathrm{Pl}}(Z)>1$ is the
\emph{planar frequency}, defined by the limit in
Section~\ref{sec:planar-frequency}.
We call $\widehat\Phi$ \emph{cylindrical} if
\[
 \dim\{v\in P_Z:\widehat\Phi(X+v)=\widehat\Phi(X)
                   \text{ for all }X\in P_Z\}=k-2,
\]
meaning that, in suitable orthonormal coordinates $X=(z,t)\in\C\times\R^{k-2}$,
it is independent of $t$.

For a cylindrical tangent function, $\nu=p_0/q_0$ is rational in
lowest terms and
\[
 \widehat\Phi=m_0\llbracket0\rrbracket+
       \sum_{j=1}^Jm_j\widehat\Phi_j,\qquad
 q=m_0+q_0\sum_{j=1}^Jm_j.
\]
The distinct nonzero \emph{cyclic components} have the form
\[
 \widehat\Phi_j(z,t)=
 \sum_{\ell=0}^{q_0-1}
 \llbracket
 \operatorname{Re}\!\left(c_jr^\nu
       e^{i\nu(\theta+2\pi\ell)}\right)
 \rrbracket,\qquad z=re^{i\theta},
\]
where $c_j\in P_Z^\perp\otimes\C\setminus\{0\}$.
The formula is independent of the choice of $\theta$ as an
unordered tuple, and $\widehat\Phi_j(0,t):=q_0\llbracket0\rrbracket$.
The integer $m_0\geq0$ is the multiplicity of the zero component,
and $m_j\geq1$ is the
\emph{repetition multiplicity} of $\widehat\Phi_j$.

Our pointwise conclusion is $m_j\geq2$ for every nonzero
component when $\nu\ne2$, without assuming uniqueness of
$\widehat\Phi$. Combined with stratification, it gives the
Hausdorff-dimension bound below. We write $\mathcal H^s$ for
$s$-dimensional Hausdorff measure and $\dim_{\HH}$ for Hausdorff
dimension. Section~\ref{sec:setup} gives the full definitions
and the results of Krummel--Wickramasekera used here.

{
\begin{introbranchtheorem}[Branch obstruction and dimension bound]
Let $T$ be a $k$-dimensional locally area-minimizing integral current in an
open set $\mathcal U\subset\R^{k+m}$, with $k\geq3$, $m\geq2$, and no
boundary in $\mathcal U$. Suppose that the induced metric on $\reg T$
is locally conformally flat. For every integer $q\geq2$, the following
statements hold.
\begin{enumerate}[label=(\roman*),leftmargin=2em]
\item Let $Z\in\mathcal R_q(T)$ and
$\nu:=N_{T,\mathrm{Pl}}(Z)\ne2$. Every cylindrical normalized tangent
function
\[
 \widehat\Phi=m_0\llbracket0\rrbracket+
                     \sum_{j=1}^Jm_j\widehat\Phi_j
\]
has $m_j\geq2$ for every nonzero component. If its rational degree is
$\nu=p_0/q_0$ in lowest terms, then
\[
 2q_0J\leq q-m_0.
\]
For $k\geq4$ this pointwise statement also holds at $\nu=2$ whenever
the normalized tangent function is cylindrical.
\item Define
\[
 \mathcal E_q:=
 \{Z\in\mathcal R_q(T):\nu:=N_{T,\mathrm{Pl}}(Z)\notin\mathbb Q
 \ \text{or}\;
 \nu=p_0/q_0\notin\mathbb Z,\ \gcd(p_0,q_0)=1,\ 2q_0>q\}.
\]
Then $\dim_{\HH}\mathcal E_q\leq k-3$. In particular,
\[
 \dim_{\HH}\{Z\in\mathcal R_2(T)\cup\mathcal R_3(T):
                    N_{T,\mathrm{Pl}}(Z)\notin\mathbb Z\}\leq k-3.
\]
\item At $\HH^{k-2}$-almost every branch point of density $q$ and planar
frequency $\nu\ne2$, the unique normalized tangent function $\widehat\Phi$ has repetition
multiplicity $m_j \geq 2$ for every nonzero component and satisfies
$2q_0J\leq q-m_0$. In particular, the nonintegral-frequency
points of densities $2$ and $3$ have zero $\HH^{k-2}$ measure.
\end{enumerate}
\end{introbranchtheorem}
}

The proof connects a local curvature calculation with the geometry of
the current. The Gauss equation gives a nonzero leading Cotton tensor
in dimension three, except at the quadratic frequency, and a nonzero
leading Weyl tensor in higher dimensions
(Theorem~\ref{thm:local}). To apply this obstruction, we first match
the multiplicities of the multivalued approximation with those of
currents with pairwise disjoint supports over profile-separated
annular sheet charts. A component of repetition multiplicity one then
gives an exact sheet-height function whose normalized derivatives
converge, allowing us to apply the curvature calculation.

The existence and estimates for tangent functions, and the
stratification, are due to Krummel--Wickramasekera
\cite{KWIII}.
Our contribution is the restriction that local conformal flatness
places on the component multiplicities and the resulting bound on
the dimension of the branch set.
For image currents of conformal $\phi$-curves, the application requires
the hypotheses discussed in Section~\ref{sec:main-hypotheses}, including
local finiteness and absence of boundary.

\subsection{Taylor flatness and conditional continuation}
Our second main result concerns critical points of weakly-conformal
maps. It applies in arbitrary codimension and for arbitrary smooth
domain and target metrics. In local coordinates, $D^bu$ denotes
the collection of partial derivatives of order $b$, with $D^1u=du$.
A smooth map is \emph{Taylor-flat at $p$} if
\[
 D^bu(p)=0\qquad\text{for every integer }b\geq1.
\]
The notation $C^{j,\alpha}$ means that the derivatives of order $j$
are locally H\"older continuous with exponent $\alpha$.
Section~\ref{sec:derivative-notation} fixes the norms and
convergence conventions.

\begin{introflatnesstheorem}[Finite-order flatness]\label{thm:intro-flatness}
Let $u:(L^k,g)\to(M^n,h)$ be $C^j$, where $j\geq2$, $n\geq k\geq3$,
and the metrics are smooth.  If $u$ is {\color{black}weakly-conformal}, then
for every $p\in L$ with $du(p)=0$, we have
\[
 D^bu(p)=0\quad(2\leq b\leq j)
\]
in any choice of local coordinates on the domain and target.
In local domain coordinates centered at $p$, as $x\to0$,
\[
 \begin{aligned}
 u(x)-u(p)&=o(|x|^j),\\
 D^bu(x)&=o(|x|^{j-b})\qquad(1\leq b\leq j).
 \end{aligned}
\]
If, in addition, $u\in C^{j,\alpha}$, $0<\alpha\leq1$, then
\[
 \begin{aligned}
 u(x)-u(p)&=O(|x|^{j+\alpha}),\\
 D^bu(x)&=O(|x|^{j+\alpha-b})\qquad(1\leq b\leq j).
 \end{aligned}
\]
In particular, every critical point of a smooth {\color{black}weakly-conformal} map is Taylor-flat.
\end{introflatnesstheorem}

For a smooth map it gives Taylor flatness at every
critical point. It does not by itself exclude critical points.

The proof starts with the rigidity of homogeneous weakly-conformal
polynomials in dimensions $n\geq k\geq3$, also proved by
P\'erez Fern\'andez--P\'erez \cite{PFP}.
We give a proof using spherical harmonics in
Theorem~\ref{thm:homogeneous-polynomial}, then apply this rigidity to
the first nonzero homogeneous Taylor term.

Baird--Wood use this Taylor-term argument in equal dimensions
\cite{BW}, and meanwhile Theorem~\ref{thm:intro-flatness} applies
in arbitrary codimension.

An additional curvature bound lets us pass from Taylor flatness to
the conclusion that a nonconstant smooth map has no critical points.

For $Z_u\ne\varnothing$, define
\[
 r_u(x):=\dist_g(x,Z_u)
       =\inf_{z\in Z_u}\dist_g(x,z),\qquad
 \gamma:=u^*h=\lambda^2g\quad\text{on }L\setminus Z_u.
\]
Here $\dist_g$ is Riemannian distance, and
$\operatorname{Scal}_\gamma$ denotes the scalar curvature of
$\gamma$. Each point of $L\setminus Z_u$ has a neighborhood $V$
on which $u$ is an embedding. On $u(V)$, the \emph{second
fundamental form} is
\[
 \mathrm{II}(A,B):=(\nabla^M_A B)^\perp
\]
for tangent vector fields $A,B$, where $\nabla^M$ is the
Levi--Civita connection of $h$ and $\perp$ denotes projection
onto the normal space. We pull its norm back to $V$ and write
$|\mathrm{II}|_\gamma$, since $u|_V$ identifies $\gamma$ with the
induced metric. We give the full definitions in
Sections~\ref{sec:energy} and~\ref{sec:energy-minimality}.

\begin{samepage}
\begin{introcontinuationtheorem}[Continuation under an inverse-distance curvature bound]\label{thm:intro-continuation}
Let $L$ be connected and let $u:(L^k,g)\to(M^n,h)$ be smooth and weakly-conformal, $k\geq3$. Assume $Z_u\ne\varnothing$ and that for every
compact $K\subset L$ there is a constant $C_K\geq0$ such that
\[
 \lambda(x)^2\operatorname{Scal}_\gamma(x)
       \geq-\frac{C_K}{r_u(x)^2}
       \qquad(x\in K\setminus Z_u).
\]
Then $u$ is constant.
  It suffices that
$r_u\lambda|\mathrm{II}|_\gamma$  is bounded on $K\setminus Z_u$
for every compact $K\subset L$.
\end{introcontinuationtheorem}
\end{samepage}

The theorem applies in particular to smooth conformal $\phi$-curves
satisfying the stated curvature bound, since these maps are
weakly-conformal.

We use the Yamabe equation to obtain a lower bound for a positive
power of the dilation $\lambda$ along a geodesic approaching a critical
point. Corollary~\ref{cor:uniform} implies that this quantity vanishes
faster than every positive power of the distance to that point,
giving a contradiction.

The Yamabe equation is the one used by Baird--Wood
\cite[equation~(11.4.3), p.~331]{BW}.
In equal dimensions, Baird--Wood prove that nonconstant smooth
weakly-conformal maps have no critical points
\cite[Theorem~11.4.6]{BW}.
It remains open whether nonconstant smooth weakly-conformal maps in
positive codimension can have critical points without the curvature
assumption.

\subsection{Organization}
Section~\ref{sec:scope} collects the map notation and definitions and
explains the relation between the map and current hypotheses.
Section~\ref{sec:setup} introduces currents, tangent functions, and
their normalizations, and states the branch theorem.
Section~\ref{sec:local} proves the local curvature obstruction.
Section~\ref{sec:extensions} treats graphs in a smooth ambient metric
under an explicit expansion hypothesis and gives the
Bor\r{u}vka-cone example explaining the quadratic exception in dimension
three.
Section~\ref{sec:current} proves the branch theorem.
Section~\ref{sec:flatness} proves Taylor flatness and the continuation
criterion. Within that section, Section~\ref{sec:boruvka} uses the
Bor\r{u}vka cone to illustrate the regularity and curvature hypotheses
and the distinction between critical points and branch points.
Section~\ref{sec:conclusion} records the remaining questions.

\subsection{Acknowledgements}
The first author acknowledges the support of the International Scientists Project BJNSF \#25ISA023.
The author thanks Spiro Karigiannis and Kotaro Kawai for useful comments.

\subsection{AI disclosure}
Artificial Intelligence was used to help identify the connection between the theory of
area-minimizing currents and the constraints imposed by local conformal
flatness. The author takes full responsibility for the contents of this paper.

\section{Preliminaries}\label{sec:scope}
We print terms in \defterm{bold} when defining them, including both standard notions and terminology introduced for this paper.

\subsection{Notation and regularity}\label{sec:derivative-notation}
Throughout the paper, $(L^k,g)$ and $(M^n,h)$ are smooth Riemannian
manifolds  without boundary,  of dimensions $k$ and $n$, respectively, with   $2\leq k\leq n$,
and $u:(L^k,g)\to(M^n,h)$ denotes a map between them.
The term \textbf{domain} refers to the manifold $L$.
We write $\vol_g$ for the \defterm{Riemannian volume measure} and, when an
orientation is fixed, for the corresponding \defterm{volume form}.
Our \defterm{Laplacian} is $\Delta_g:=\operatorname{div}_g\nabla_g$.

Throughout the paper, for a $C^1$ map $u:L \rightarrow M$, $du$ denotes the \defterm{first differential}
$du_x:T_xL\to T_{u(x)}M$, also written as a matrix after choosing bases.
In domain coordinates $x^1,\ldots,x^k$ and target coordinates
$u^1,\ldots,u^n$, its matrix entries are $\partial_i u^\mu$, where
$\partial_i:=\partial/\partial x^i$. For an integer $b\geq2$, write
\[
 D^bu(x):=
 \left(\partial_{i_1}\cdots\partial_{i_b}u^\mu(x)\right)_{
 1\leq\mu\leq n,\ 1\leq i_1,\ldots,i_b\leq k},
\]
  where the entries are all coordinate derivatives of order $b$.
These derivatives  depend on the chosen coordinates. The same notation applies to other maps and scalar
functions. Subscripts specify individual derivatives: for example,
$f_{ij}:=\partial_i\partial_j f$ is an entry of $D^2f$.

  In fixed coordinates, $|\cdot|$ denotes  the Euclidean
norm. We write $h\in C^b(A,\R^N)$ if $A\subset\R^k$ is open and $h$ has
continuous derivatives through order $b$. Set the \defterm{$C^0$ norm}
\[
 \|h\|_{C^0(A)}:=\sup_{X\in A}|h(X)|.
\]
For $b\geq1$, we use the \defterm{$C^b$ norm}
\[
 \|h\|_{C^b(A)}:=\|h\|_{C^0(A)}
 +\sup_{X\in A}|dh(X)|
 +\sum_{j=2}^b\sup_{X\in A}|D^jh(X)|.
\]
 Convergence in
$C^b_{\mathrm{loc}}(A)$ means convergence of all derivatives of
orders at most $b$, uniformly on each compact subset of $A$.
For a compact set $K$ and a function of class $C^b$ on a neighborhood of $K$,
$\|\cdot\|_{C^b(K)}$ is defined by the same formula, with suprema over $K$.
For a tensor field in fixed coordinates, we apply these norms and convergence definitions to the vector of its component functions.

For functions or tensors on a fixed domain, $h_\rho=O(\rho^\mu)$ in a
specified norm $\|\cdot\|$ means $\|h_\rho\|\leq C\rho^\mu$ for all sufficiently
small $\rho>0$, with $C$ independent of $\rho$. $h_\rho=o(\rho^\mu)$ in that
norm means $\|h_\rho\|/\rho^\mu\to0$ as $\rho\downarrow0$.
This includes $o(1)$ by taking $\mu=0$.
For scalar quantities, we use absolute value.

A map between manifolds is $C^b$ if it is $C^b$ in local coordinates.
  For $A\subset\R^k$ and $v:A\to\R^N$, the
\defterm{H\"older seminorm of exponent $\alpha\in(0,1]$} is
\[
 [v]_{\alpha;A}:=
 \sup_{\substack{X,Y\in A\\X\ne Y}}
 \frac{|v(X)-v(Y)|}{|X-Y|^\alpha}.
\]
With $D^0h:=h$ and $D^1h:=dh$, set
$\|h\|_{C^{b,\alpha}(A)}:=\|h\|_{C^b(A)}+[D^bh]_{\alpha;A}$.
We write $h\in C^{b,\alpha}_{\mathrm{loc}}$ if this norm is finite on
every relatively compact coordinate ball. For maps between manifolds,
 $C^{b,\alpha}$   means this condition in coordinates.
We also write $C^\alpha=C^{0,\alpha}$.
 A smooth map is \defterm{Taylor-flat at $p$} if every positive-order
coordinate derivative vanishes at $p$. The value $u(p)$ need not be zero.

\subsection{Weak conformality and local conformal flatness}\label{sec:energy}
A smooth Riemannian metric is \defterm{locally conformally flat} if
it can be written as $e^{2f}\delta$ in local coordinates around each
point, where $f$ is smooth and $\delta$ is the Euclidean metric.
This property is preserved by pullback under a local diffeomorphism
and by multiplication by a smooth positive function.

For a linear map $A$ between inner-product spaces, its
\defterm{Hilbert--Schmidt norm} is
\[
 |A|_{\rm HS}^2:=\sum_{i=1}^k|Ae_i|^2,
\]
where $(e_i)$ is any orthonormal basis of the domain.
We write $|du|$ for $|du|_{\rm HS}$ and represent $du$ in
orthonormal frames in the matrix formulas below.

For inner-product spaces $(V,g)$ and $(W,h)$ with $\dim V=k$, define
\[
 \mathrm{CO}(V,W):=\{{\color{black}A:V\to W\text{ linear}:}A^*h=\lambda^2g
                         \text{ for some }\lambda\geq0\}.
\]
A $C^1$ map is \defterm{{\color{black}weakly-conformal}} if
$du_p\in\mathrm{CO}(T_pL,T_{u(p)}M)$ at every point.
Equivalently, $u^*h=\lambda^2g$ for a function
$\lambda:L\to[0,\infty)$, called \defterm{dilation}. At each point, either $du$ has rank $k$
or $du=0$.
In orthonormal frames,
\[
 \mathrm{CO}(k,n):=\{A\in\R^{n\times k}:A^TA=\lambda^2I_k
                                \text{ for some }\lambda\geq0\}.
\]
Thus $du\in\mathrm{CO}(k,n)$ means precisely $(du)^Tdu=\lambda^2I_k$.
This implies that $\lambda=|du|/\sqrt{k}$, which gives that it is continuous everywhere and smooth wherever $u$ is smooth and $\lambda>0$.
For a {\color{black}weakly-conformal} map, a \defterm{critical point} is a point $p$
where $du_p=0$, and the \defterm{critical set} is $Z_u:=\{p:du_p=0\}$.
The   \defterm{pullback metric} is
$\gamma:=u^*h=\lambda^2g$ on $L\setminus Z_u$.
  When $u$ is smooth, this is a smooth metric, whose  scalar curvature
is denoted by $\operatorname{Scal}_\gamma$.
 All curvature expressions for $\gamma$ below are taken on $L\setminus Z_u$.
 For a smooth weakly-conformal map $u$, each point of $L\setminus Z_u$
has an open neighborhood $V\subset L\setminus Z_u$ such that
$u|_V$ is a smooth embedding. The induced metric on $u(V)$
pulls back under $u|_V$ to $\gamma|_V=\lambda^2g|_V$.
Thus, if $g$ is locally conformally flat, the induced metric on
each such submanifold $u(V)$ is locally conformally flat.
If $Z_u\ne\varnothing$, put
\[
 r_u(x):=\dist_g(x,Z_u):=\inf_{z\in Z_u}\dist_g(x,z),
\]
where $\dist_g(x,z)$ is the infimum of the $g$-lengths of piecewise
$C^1$ curves joining $x$ to $z$.

\subsection{Energy and minimality}\label{sec:energy-minimality}
For an open set $V\Subset L$, meaning that its closure is compact in $L$,
the \defterm{$k$-energy} is
\[
 E_k(u;V):=\frac{1}{(\sqrt{k})^k}\int_V|du|^k\,\vol_g.
\]
A $C^1$ map is \defterm{$k$-harmonic} if
  $\left.\frac{d}{dt}\right|_{t=0}E_k(u_t;V)=0$ for every relatively
compact open $V$ and every $C^1$ variation $u_t$ with $u_0=u$
that equals $u$ outside a fixed compact subset of $V$.
Here $t\mapsto u_t$ is continuously differentiable in the $C^1$
topology in local coordinates.
Equivalently, its \defterm{weak $k$-harmonic equation} is
\[
 \int_V |du|^{k-2}
       \sum_{i=1}^k h\bigl(du(e_i),\nabla^u_{e_i}\zeta\bigr)
       \,\vol_g=0
 \qquad
 \bigl(\zeta\in C_c^1(V,u^*TM)\bigr).
\]
Here $C_c^1$ denotes continuously differentiable sections with compact
support. For $C^2$ maps this is the  Euler--Lagrange equation
  \[
 \operatorname{div}_{g,h}(|du|^{k-2}du)=0,
\]
 see \cite[\S2, p.~218]{Takeuchi}.
Here $\operatorname{div}_{g,h}$ is the \defterm{covariant divergence}, which
for a   $C^1$  one-form $A$ with $A(X)\in T_{u(x)}M$ when $X\in T_xL$, is defined as
\[
 \operatorname{div}_{g,h}A:=
 \sum_{i=1}^k\left\{\nabla^u_{e_i}\bigl(A(e_i)\bigr)
                    -A(\nabla^g_{e_i}e_i)\right\},
\]
where the frame $(e_i)$ is $g$-orthonormal, $\nabla^g$ is the Levi--Civita
connection of $g$, and $\nabla^u$ differentiates vector fields along $u$
using the Levi--Civita connection $\nabla^M$ of $h$.
For Euclidean metrics this is the usual componentwise divergence,
written simply as $\operatorname{div}$ below.

For a   $C^2$  immersion $u:L\to(M,h)$, let $U\subset L$ be an
open set such that $u|_U$ is an embedding. For   $C^1$  tangent
vector fields $V,W$ on the submanifold $u(U)$, the
\defterm{second fundamental form} is
\[
 \mathrm{II}(V,W):=(\nabla^M_VW)^\perp.
\]
Here $\perp$ is the $h$-orthogonal projection to the normal space
of $u(U)$, and the norm of $\mathrm{II}$ is taken with the metric
induced by $h$ on $u(U)$.
A  $C^2$  {\color{black}weakly-conformal} $k$-harmonic immersion has a \defterm{minimal image},
meaning that the trace of this second fundamental form vanishes
\cite[Corollary~4, pp.~221--222]{Takeuchi}.

\subsection{\texorpdfstring{Calibrations and conformal $\phi$-curves}{Calibrations and conformal phi-curves}}
\label{sec:smith-scope}
A \defterm{simple $k$-vector} is a wedge of $k$ vectors.
A smooth $k$-form $\phi$ on $(M^n,h)$ is a \defterm{calibration} if
$\dd\phi=0$ and $\phi(\xi)\leq1$ for every oriented unit simple
$k$-vector $\xi$. An oriented $k$-plane is
\defterm{$\phi$-calibrated} when equality holds on its oriented unit
$k$-vector. When considering conformal $\phi$-curves, we fix an orientation on $L$.
A $C^1$ map $u:(L^k,g)\to(M^n,h)$ is a
\defterm{conformal $\phi$-curve} if
\[
 u^*h=\lambda^2g,\qquad u^*\phi=\lambda^k\vol_g
\]
for some function $\lambda:L\to[0,\infty)$. The first identity implies that
\[
 \lambda=\frac{|du|}{\sqrt{k}}.
\]
Conversely, with $\lambda=|du|/\sqrt{k}$, the second identity implies the first.
We use the terminology of Ikonen--Pankka \cite[\S1]{IP}.
The literature also uses \defterm{Smith map}, in the vector-cross-product
setting of Cheng--Karigiannis--Madnick \cite[\S1.1 and \S2]{CKM},
and \defterm{Smith immersion}, for general calibrations in
Iliashenko--Karigiannis \cite[Definition~3.3]{IK}.
The definition allows critical points. We do not require $du$ to have full
rank everywhere \cite[Definition~3.3 and the remark following it]{IK}.
We keep using the name conformal $\phi$-curve throughout.
At each noncritical point $p$, the image tangent plane $du_p(T_pL)$, with orientation induced by $L$, is $\phi$-calibrated.

Closedness of $\phi$ implies that every conformal $\phi$-curve is
$k$-harmonic. Smith established the energy-minimizing property in
the vector-cross-product setting \cite[Corollary~6.8]{SmithOriginal}.
Iliashenko--Karigiannis treat general closed calibrations
\cite[Theorems~3.4 and~3.9]{IK}.  Even though originally the maps were considered to be smooth, the same variational argument applies to
$C^1$ maps. The variational approach originates in Smith's work
\cite[\S6]{SmithOriginal}. Related developments include regularity and bubble-tree theory of
Cheng--Karigiannis--Madnick \cite{CKM} and their variational treatment
of calibrated submanifolds \cite{CKMVar}.
Pankka \cite[\S1]{Pankka} studies the broader class of quasiregular curves,
and Ikonen--Pankka \cite[\S1]{IP} study conformal curves and their
Liouville-type rigidity.

Our interest in conformal $\phi$-curves also comes from hyperbolicity:
Broder--Iliashenko--Madnick use maps from the Poincar\'e $k$-ball
in place of holomorphic disks to define calibrated analogues of
Kobayashi--Royden hyperbolicity \cite[\S5.1 and \S6]{BIM}.
Iliashenko--Karigiannis--Madnick relate these notions to equicontinuity
and compactness of families of conformal $\phi$-curves, providing calibrated analogues of Royden and Montel theorems
\cite[Theorem~1.2 and Corollary~1.3]{IKM}.
These results motivate our study of the regularity and critical points
of such maps.

\subsection{Relation to the main results}\label{sec:main-hypotheses}
An integral current calibrated by a closed form is locally area minimizing
\cite[Theorem~II.4.2]{HL}. We recall the definitions and the argument in
Section~\ref{sec:area-minimizing-currents}.

  For a smooth conformal $\phi$-curve with  locally conformally flat  domain,
  the induced metric on each submanifold described in
Section~\ref{sec:energy} is locally conformally flat.
This explains the current theorem's application to suitable Euclidean
image currents of conformal $\phi$-curves.  One must still verify that an image current exists,
has locally finite mass and integer multiplicity, and no boundary.  An arbitrary conformal $\phi$-curve need not have these properties.
The roles of the necessary assumptions can be summarized as follows.
\begin{enumerate}[label=(\roman*),leftmargin=2em]
\item \emph{Branch obstruction, Theorem~\ref{thm:intro-branch}.}
  The object is a locally area-minimizing integral current $T$ of dimension
$k\geq3$ and codimension at least two in Euclidean space, without boundary, and   with locally conformally flat induced
metric on $\reg T$. The pointwise and dimension conclusions have the additional
branch-point and frequency restrictions stated in that theorem.
It is not necessary that the current comes from a conformal $\phi$-curve.
For a conformal $\phi$-curve image current, calibration supplies area minimization once
the current is well defined, but one must still verify that it has locally finite mass and no boundary.   Smooth weakly-conformal $k$-harmonic maps   locally parametrize minimal submanifolds away from their critical sets, but minimality alone is not sufficient for this
area-minimizing hypothesis.
\item \emph{Taylor flatness, Theorem~\ref{thm:intro-flatness}.}
The object is a $C^j$ {\color{black}weakly-conformal} map, $j\geq2$, between smooth
Riemannian manifolds of dimensions $n\geq k\geq3$.
We do not assume the existence of an image current, minimality,
calibration, $k$-harmonicity, or flatness of the metric. Thus the result applies to conformal $\phi$-curves with this
regularity on arbitrary smooth Riemannian domains and targets.
\item \emph{Conditional continuation,
Theorem~\ref{thm:intro-continuation}.}
The map is smooth and {\color{black}weakly-conformal}, its domain is connected and has
dimension $k\geq3$, and the stated curvature bound holds.
For   every smooth weakly-conformal map, the contracted Gauss equation
 makes the bound on $r_u\lambda|\mathrm{II}|_\gamma$ sufficient.
  This includes smooth conformal $\phi$-curves and weakly-conformal
$k$-harmonic maps.
\end{enumerate}
Proposition~\ref{prop:curved} concerns a smooth ambient metric but assumes
an explicit graph expansion. We do not extend the current theorem to
general curved targets. The local graph obstruction in
Theorem~\ref{thm:local}, by contrast, uses the Cotton and Weyl tensors and the
stated convergence of the graph functions, not minimality or calibration.

\section{Currents and tangent functions}\label{sec:setup}
{\color{black}In the current-theoretic sections, $m\geq2$
denotes the codimension of the $k$-dimensional current, so its ambient
Euclidean space is $\mathbb R^{k+m}$.}
We now introduce the current-theoretic notation used in
Theorem~\ref{thm:intro-branch}. Throughout this section, $Z$ denotes a
point of the ambient Euclidean space. In graph formulas, $X$ lies in the
tangent plane and $Y$ in its orthogonal complement. We reserve $p$ for
a point in the domain of the original map.
\subsection{Area-minimizing currents, tangent cones, and branch points}
\label{sec:area-minimizing-currents}

The image of a smooth immersion need not be an embedded
submanifold. Integral currents provide the appropriate setting
because they retain the orientations and integer multiplicities.
We recall below the definitions and conventions used in this paper,
following Simon \cite[Chapters~1, 3, 6, and~7]{SimonGMT}.

Let $\mathcal U\subset\R^{k+m}$ be open. All currents in the following
definition are defined in $\mathcal U$.  We use
$B_r(Z):=\{W:|W-Z|<r\}$ for an ambient Euclidean ball and
$B_r^P(X):=\{Y\in P:|Y-X|<r\}$ for a ball in the plane $P$.
For sets $A,B$, the notation $A\Subset B$ means that $\overline A$
is compact and contained in $B$. For a measure $\mu$ and a measurable
set $E$, its restriction is $(\mu\llcorner E)(B):=\mu(B\cap E)$.

We use $\mathcal H^k$ for \defterm{$k$-dimensional Hausdorff measure}, normalized
to equal Lebesgue measure on a $k$-plane. For a set $E$ in Euclidean space, its \defterm{Hausdorff dimension} is
$\dim_{\HH}E:=\inf\{s\geq0:\HH^s(E)=0\}$, with
$\dim_{\HH}\varnothing:=-\infty$.

\begin{definition}[Integral current and local area minimization]\label{def:integral-current}
  An $\mathcal H^k$-measurable  set $E \subset \mathcal U$ is \defterm{countably $k$-rectifiable} if, outside a set of
zero $k$-dimensional Hausdorff measure, it is covered by countably many
Lipschitz images of subsets of $\R^k$.
A $k$-dimensional \defterm{integer-rectifiable current} $T$ has,
for a countably $k$-rectifiable set $E\subset\mathcal U$,  the form
\[
 T(\omega):=\int_E \theta(X)\langle\omega(X),\xi(X)\rangle
                      \,\dd\mathcal H^k(X).
\]
Here $\omega$ is a smooth compactly supported differential $k$-form,
$\theta:E\to\{1,2,\ldots\}$ is locally integrable in $\mathcal U$ with
respect to $\mathcal H^k\llcorner E$, meaning
$\int_{E\cap K}\theta\,\dd\mathcal H^k<\infty$ for every compact
$K\subset\mathcal U$, and $\xi(X)$ is a
measurable oriented unit simple $k$-vector spanning the approximate tangent
plane of $E$ at almost every $X$.
At $\mathcal H^k$-almost every such $X$, the \defterm{approximate tangent plane}
$P_X$ is characterized by
\[
 \lim_{\rho\downarrow0}\rho^{-k}
    \int_E f((Y-X)/\rho)\,\dd\mathcal H^k(Y)
       =\int_{P_X}f(Y)\,\dd\mathcal H^k(Y)
       \quad(f\in C_c(\R^{k+m})).
\]
Here $C_c$ denotes continuous compactly supported functions.
The \defterm{boundary} is the current
$\partial T(\eta):=T(\dd\eta)$.
The current is \defterm{integral} if its boundary is also integer rectifiable
with locally finite mass.
Its \defterm{mass measure} and \defterm{mass} are
\[
 \|T\|(B):=\int_{E\cap B}\theta\,\dd\mathcal H^k,
 \qquad \mathbf M(T):=\|T\|(\mathcal U).
\]
Thus $T$ has \defterm{locally finite mass}: $\|T\|(K)<\infty$ for every
compact $K\subset\mathcal U$, by the local integrability of $\theta$.
Its \defterm{support} is
$\supp T:=\{X\in\mathcal U:\|T\|(B_r(X))>0\text{ for every }r>0
\text{ with }B_r(X)\Subset\mathcal U\}$.
  For a Borel set $B$, the  \defterm{restriction} $T\llcorner B$ is defined by restricting the
integral to $E\cap B$.
The symbol $\llbracket E\rrbracket$ denotes integration with multiplicity
one and the specified orientation.

The current is \defterm{locally area minimizing} in $\mathcal U$ if, for every
open $W\Subset\mathcal U$ and every integral $(k+1)$-current $S$ with
$\supp S\Subset W$,
\[
 \mathbf M(T\llcorner W)
 \leq \mathbf M((T+\partial S)\llcorner W).
\]
Adding $\partial S$, with $W$ and $S$ as above, preserves the
boundary of $T$ because $\partial^2S=0$. Local area minimization
means that this addition cannot decrease its mass in $W$.
\end{definition}

The mass measure $\|T\|$ is a \defterm{Radon measure}: a nonnegative Borel
measure finite on compact sets, inner regular and outer regular.
Radon measures $\mu_i$ \defterm{converge weakly locally} to $\mu$ in $\mathcal U$,
written $\mu_i\to\mu$, if
\[
 \int_{\mathcal U} f\,\dd\mu_i\longrightarrow
 \int_{\mathcal U} f\,\dd\mu
 \qquad\text{for every }f\in C_c(\mathcal U),
\]
where $C_c(\mathcal U)$ denotes the continuous real-valued functions
with compact support in $\mathcal U$. Convergence of mass measures is understood in this sense.

An integral current $T$ is \defterm{$\phi$-calibrated} if
$\phi(X)(\xi(X))=1$ for $\|T\|$-almost every $X$.
{We claim that such a current is locally area minimizing when
$\dd\phi=0$, as in \cite[Theorem~II.4.2]{HL}.
Take $W$ and $S$ as in the preceding definition.} Since $\supp S\Subset W$, closedness gives
$(\partial S)(\phi)=S(\dd\phi)=0$, where we insert a smooth compactly
supported cutoff equal to one near $\supp S$ when evaluating these currents.
Calibration and the pointwise bound on $\phi$ then give
\[
 \mathbf M(T\llcorner W)
 =(T\llcorner W)(\phi)
 =((T+\partial S)\llcorner W)(\phi)
 \leq\mathbf M((T+\partial S)\llcorner W).
\]
Here evaluation on $\phi$ for the restricted currents means integration
over $W$ in the representation above.

\begin{definition}[Magnification and tangent cone]\label{def:magnification}
The \defterm{pushforward} under a smooth map $f$ is defined by
$(f_\#T)(\omega):=T(f^*\omega)$   when $f$ is proper on $\supp T$,
meaning that $\supp T\cap f^{-1}(K)$ is compact for every compact
set $K$ in the target. Here $T(f^*\omega)$ means $T(\chi f^*\omega)$,
where $\chi$ is smooth, compactly supported, and equal to $1$ near
$\supp T\cap f^{-1}(\supp\omega)$. The value is independent of the
choice of $\chi$.
 For a Lipschitz map proper on $\supp T$, we use the pushforward defined in
\cite[Chapter~6, \S3, p.~174]{SimonGMT}.

For $Z\in\R^{k+m}$ and $\rho>0$, define the \defterm{magnification map}
 by
\[
 \eta_{Z,\rho}(X):=\frac{X-Z}{\rho}.
\]

  The magnification map
$\eta_{Z,\rho}:\mathcal U\to\eta_{Z,\rho}(\mathcal U)$ is proper on $\supp T$.
  Define  the \defterm{magnified current} by
\[
 T_{Z,\rho}:=(\eta_{Z,\rho})_\#T.
\]

A sequence of $k$-currents $S_i$ \defterm{converges weakly} to $S$ in
$\mathcal U$ if $S_i(\omega)\to S(\omega)$ for every smooth compactly
supported $k$-form $\omega$ in $\mathcal U$.

For an open set $W\Subset\mathcal U$ and a $k$-current $S$, define its
\defterm{flat seminorm in $W$} by
\[
 \mathcal F_W(S):=\inf\bigl\{\mathbf M(R\llcorner W)
  +\mathbf M(Q\llcorner W): S=R+\partial Q\text{ in }W\bigr\},
\]
where the infimum is over rectifiable $k$-currents $R$ and
$(k+1)$-currents $Q$ of locally finite mass, with the equality tested
against forms compactly supported in $W$.
A sequence $S_i$ \defterm{converges locally in the flat topology} to $S$ if
\[
 \mathcal F_W(S_i-S)\longrightarrow0
 \qquad\text{for every open }W\Subset\mathcal U.
\]
Local flat convergence implies weak convergence of currents.
Unless otherwise stated, convergence of currents \(S_i\to S\) means local flat convergence.
A current $C$ is a \defterm{tangent cone of $T$ at $Z$} if
there is a sequence $\rho_i\downarrow0$ such that
$T_{Z,\rho_i}\to C$.
\end{definition}

If $T$ is locally area minimizing and
$Z\in\supp T\setminus\supp(\partial T)$, every tangent cone $C$ of $T$
at $Z$ is invariant under every dilation $\delta_\lambda(X):=\lambda X$
\cite[Chapter~7, Theorem~3.1(2), p.~217]{SimonGMT}:
\[
 (\delta_\lambda)_\#C=C\qquad(\lambda>0).
\]
Different sequences $\rho_i\downarrow0$ may produce different tangent cones.
The tangent cone is \defterm{unique} if every such sequence produces the same
current $C$.

\paragraph{Monotonicity.}
Let $T$ be a locally area-minimizing integral $k$-current in
$\mathcal U\subset\mathbb R^{k+m}$. If
$\overline{B_R(Z)}\subset\mathcal U$ and $\partial T=0$ in $B_R(Z)$,
then
\begin{equation}\label{eq:current-monotonicity}
 \frac{\|T\|(B_\sigma(Z))}{\omega_k\sigma^k}
 \leq\frac{\|T\|(B_\rho(Z))}{\omega_k\rho^k}
 \qquad(0<\sigma<\rho<R),
\end{equation}
where $\omega_k$ denotes the $k$-dimensional Lebesgue measure of the unit
ball in $\mathbb R^k$. Thus the mass ratio is nondecreasing with radius.
This is the Euclidean monotonicity formula for locally area-minimizing
currents~\cite[Chapter~4, equations~(3.8)--(3.10), p.~109, and
Chapter~7, Theorem~3.1(1), p.~217]{SimonGMT}.

\begin{definition}[  Density, singular point, branch point, and
flat singular point]\label{def:density-branch}
Let $T$ be a locally area-minimizing $k$-current and let
$Z\in\supp T\setminus\supp(\partial T)$.  The \defterm{density of $T$ at $Z$}
is
\[
 \Theta_k(\|T\|,Z):=\lim_{\rho\downarrow0}
 \frac{\|T\|(B_\rho(Z))}{\omega_k\rho^k},
\]
and the limit exists and is finite by \eqref{eq:current-monotonicity}
and local finiteness of the mass. At every such support point the density
is at least one, so monotonicity also gives the lower mass bound
\[
 \|T\|(B_\rho(Z))\geq\omega_k\rho^k\qquad(0<\rho<R)
\]
whenever $R$ satisfies the hypotheses above. The density bound follows
from integer multiplicity and upper semicontinuity of the density
\cite[Chapter~7, Theorem~3.1(1), p.~217]{SimonGMT}.

  Define $\reg T$ to be the set of points
$Z\in\supp T\setminus\supp(\partial T)$ for which there exist an open
neighborhood $W$ of $Z$, an oriented smooth embedded minimal
$k$-submanifold $\Sigma\subset W$ containing $Z$, and an integer $q\geq1$
such that $T\llcorner W=q\llbracket\Sigma\rrbracket$.
The points of $\reg T$ are called \defterm{regular points of $T$}. Set
\[
 \sing T:=(\supp T\setminus\supp\partial T)\setminus\reg T.
\]
A point in $\sing T$ is a \defterm{singular point}.  Following the terminology
of Krummel--Wickramasekera, a singular point $Z$ is a \defterm{branch point} if
at least one tangent cone of $T$ at $Z$ is supported on a $k$-plane \cite[p.~3]{KWIII}.  Such a tangent
cone has the form $q\llbracket P\rrbracket$ for an integer $q\geq2$.

{\color{black}By \cite[Chapter~7, Remark~3.2(1), p.~217]{SimonGMT},
such a tangent cone $C$ satisfies $C=q\llbracket P\rrbracket$ for
an integer $q\geq1$, after orienting $P$.}

{\color{black}  Multiplicity one would imply $Z\in\reg T$ by
Allard's theorem \cite[Chapter~7, Theorem~4.3, p.~224]{SimonGMT},
contradicting its singularity.}

A singular point $Z$ is a \defterm{flat singular point of density $q$} if its
tangent cone is unique and equals
\[
 q\llbracket P_Z\rrbracket,
\]
where $P_Z$ is an oriented $k$-plane through the origin and
$\llbracket P_Z\rrbracket$ denotes integration over that plane with
multiplicity one.  The plane $P_Z$ is the \defterm{tangent plane of $T$ at $Z$}.
Note that a flat singular point is a branch point with a unique tangent cone.
\end{definition}

By Almgren's theorem~\cite[Chapter~7, \S4, p.~223]{SimonGMT}, the singular set, and hence
the branch set, has Hausdorff dimension at most $k-2$.

For a nonempty set $E\subset\R^{k+m}$, the \defterm{Euclidean distance}
from a point $W$ to $E$ is
\[
 \dist(W,E):=\inf_{Y\in E}|W-Y|.
\]
Define the \defterm{affine tangent plane} by $Z+P_Z:=\{Z+V:V\in P_Z\}$. Then
\[
 \dist(W,Z+P_Z)=\inf_{V\in P_Z}|W-Z-V|.
\]
This is the perpendicular distance to that plane: if $W-Z=X+Y$ with
$X\in P_Z$ and $Y\in P_Z^\perp$, then $\dist(W,Z+P_Z)=|Y|$.

The pointwise obstruction and dimension estimate in
Theorem~\ref{thm:intro-branch}(i)--(ii) apply to the following rapid-decay
branch set $\mathcal R_q(T)$. Part~(iii) gives a separate almost-everywhere
conclusion on the full branch set.

\begin{definition}[Rapid-decay branch set]\label{def:rapid}
For an integer $q\geq2$, define $\mathcal R_q(T)$ to be the set of flat
singular points $Z$ of density $q$ for which there exist constants
$\mu\in(0,1)$, $C<\infty$, and $r_0>0$ such that
\begin{equation}\label{eq:rapid}
 \int_{B_r(Z)}\dist^2(W,Z+P_Z)\,\dd\|T\|(W)
       \leq C r^{k+2+2\mu}\qquad(0<r<r_0).
\end{equation}
We call $\mathcal R_q(T)$ the \defterm{rapid-decay branch set of density $q$}.
The constants may depend on $Z$.
\end{definition}

Later we integrate over cylinders that restrict the coordinate along
the tangent plane but not the normal coordinate. Such cylinders may
contain distant parts of the current. The following lemma allows us to
retain only the part near the point under study, without introducing
boundary inside the cylinder or losing local area minimization.

\begin{lemma}[Normal localization]\label{lem:localization}
Let \(T\) be a locally area-minimizing integral \(k\)-current in an open set \(\mathcal U\subset\mathbb R^{k+m}\), and let
\[
Z\in\operatorname{supp}T\setminus\operatorname{supp}(\partial T).
\]Suppose \(T\) has the unique tangent cone \(q\llbracket P\rrbracket\) at \(Z\).
After translation of $Z$ to $0$, there is $R>0$ such that the part of $T$
in $B_R^P(0)\times B_{2R}^{P^\perp}(0)$, extended by zero in the normal
direction, is boundaryless and locally area minimizing in
$B_R^P(0)\times P^\perp$.  It has bounded normal support and agrees with
$T$ near $0$.
\end{lemma}
\begin{proof}

First, we need to choose $R>0$. Since we assumed
$Z=0\in\supp T\setminus\supp(\partial T)$, there exists $r_{\star}>0$
such that $\overline{B_{r_{\star}}(0)}\subset\mathcal U$ and
$\partial T=0$ in $B_{r_{\star}}(0)$.
For any $r>0$, define
\[
 \eta_r(W):=W/r\quad(W\in\mathbb R^{k+m}),
 \qquad T_r:=(\eta_r)_\#T.
\]
We claim that, as $r\downarrow0$,
\begin{equation}\label{eq:localization-tangent-convergence}
 T_r\to q\llbracket P\rrbracket\quad\text{locally flatly},
 \qquad \|T_r\|\to q\mathcal H^k\llcorner P\quad\text{weakly locally}.
\end{equation}

Fix any sequence $r_j\downarrow0$ and $L>0$. Then $T_{r_j}$ is locally
area minimizing in $\mathcal U_{r_j}:=\eta_{r_j}(\mathcal U)$, since a
competitor pulls back to a competitor for $T$, and magnification scales
both masses by $r_j^{-k}$.
For sufficiently large $j$, $2Lr_j<r_{\star}$, so
$B_{2L}(0)\subset\mathcal U_{r_j}$ and $\partial T_{r_j}=0$ in $B_{2L}(0)$.
We can discard the first finitely many terms for which this might not hold,
and hence regard all the remaining currents as currents in the same domain
$B_{2L}(0)$, where they are also locally area minimizing.

Now, for any open $V\Subset B_{2L}(0)$, choose $0<b<2L$ with
$V\subset B_b(0)$. Then
\[
 \|T_{r_j}\|(V)\leq r_j^{-k}\|T\|(B_{br_j}(0))
 =b^k\frac{\|T\|(B_{br_j}(0))}{(br_j)^k},
\]
where, by monotonicity, the rightmost term is bounded for sufficiently
small $r_j$. Since $\partial T_{r_j}=0$ in $V$, we get
\[
 \sup_j\bigl(\|T_{r_j}\|(V)+\|\partial T_{r_j}\|(V)\bigr)<\infty.
\]
Then, by the Federer--Fleming compactness theorem
\cite[Chapter~6, Theorem~8.2]{SimonGMT}, a subsequence of $T_{r_j}$
converges weakly to an integral current $C$. Applying this to successively
larger balls $B_{2L}(0)$ and taking a diagonal subsequence, we get weak
convergence in $\mathbb R^{k+m}$. By
\cite[Chapter~6, Theorem~7.2, p.~197]{SimonGMT}, this convergence is also
locally flat. Thus $C$ is a tangent cone, and uniqueness gives
$C=q\llbracket P\rrbracket$.

This actually implies that the whole original sequence $T_{r_j}\to C$
locally flatly. If not, we can find an open $V\Subset\mathbb R^{k+m}$,
$\varepsilon>0$, and a subsequence with
$\mathcal F_V(T_{r_{j_\ell}}-C)\geq\varepsilon$ for every $\ell$.
Applying compactness and the weak-to-flat convergence theorem again gives
a further subsequence converging locally flatly to $q\llbracket P\rrbracket$,
contradicting this lower bound.
Then we apply \cite[Chapter~7, Theorem~2.4, p.~214]{SimonGMT} on each fixed
ball, using local area minimization as established above, to conclude that
\[
 \|T_{r_j}\|\to\|q\llbracket P\rrbracket\|
 =q\mathcal H^k\llcorner P
 \quad\text{weakly locally}.
\]
Since the sequence $r_j\downarrow0$ was arbitrary, this proves the claimed
convergences as $r\downarrow0$.

Now, we claim that for all sufficiently small $R>0$ we have
\begin{equation}\label{eq:localization-support-strip}
 \supp T\cap\overline{B_{3R}(0)}
 \subset\{W\in\mathbb R^{k+m}:|\pi_{P^\perp}W|<R/2\},
\end{equation}
where $\pi_{P^\perp}$ denotes orthogonal projection onto $P^\perp$.
If this were not true, there would exist a sequence $R_j\downarrow0$
and points $W_j\in\supp T\cap\overline{B_{3R_j}(0)}$ with
$|\pi_{P^\perp}W_j|\geq R_j/2$. Setting $A_j:=W_j/R_j$, we get
$A_j\in\supp T_{R_j}\cap\overline{B_3(0)}$ and
$\dist(A_j,P)\geq1/2$. After passing to a subsequence,
$A_j\to A\in\overline{B_3(0)}$ with $\dist(A,P)\geq1/2$.

Set $s:=1/8$. For large $j$, $B_s(A_j)\subset B_{2s}(A)$,
and $\overline{B_{2s}(A)}$ is disjoint from $P$.
Now, we claim that
\[
 \|T_{R_j}\|(B_{2s}(A))\geq\omega_k(s/2)^k>0.
\]
This follows from the monotonicity formula, since
$A_j\in\supp T_{R_j}$, we have
$\|T_{R_j}\|(B_{s/2}(A_j))>0$. At $\|T_{R_j}\|$-almost every point,
the density equals its positive integer multiplicity, so there exists
$X_j\in B_{s/2}(A_j)$ with
\[
 \Theta_k(\|T_{R_j}\|,X_j)
 =\lim_{t\downarrow0}\frac{\|T_{R_j}\|(B_t(X_j))}{\omega_k t^k}
 \geq1.
\]
For sufficiently large $j$, $T_{R_j}$ is locally area minimizing and
has no boundary in $B_4(0)$, which contains all the balls used here.
Hence, by \eqref{eq:current-monotonicity},
\[
 \frac{\|T_{R_j}\|(B_{s/2}(X_j))}{\omega_k(s/2)^k}
 \geq\Theta_k(\|T_{R_j}\|,X_j)\geq1.
\]
Finally, $|X_j-A_j|<s/2$ implies
\[
 B_{s/2}(X_j)\subset B_s(A_j)\subset B_{2s}(A).
\]
Taking the $\|T_{R_j}\|$-measure of these balls gives the required claim.

Now choose a continuous compactly supported test function $f$ with
$0\leq f\leq1$, equal to one on $K:=\overline{B_{2s}(A)}$, and
supported away from $P$. The measure convergence proved above gives
\[
 0<\omega_k(s/2)^k
 \leq\|T_{R_j}\|(B_{2s}(A))
 \leq\|T_{R_j}\|(K)
 \leq\int f\,\dd\|T_{R_j}\|
 \longrightarrow\int f\,\dd(q\mathcal H^k\llcorner P)=0,
\]
which is a contradiction. This proves \eqref{eq:localization-support-strip}.

Fix such an $R$ with $4R<r_{\star}$. If $X\in P$, $Y\in P^\perp$,
$|X|<R$, and $R\leq|Y|\leq2R$, then orthogonality gives
$|X+Y|<\sqrt{5}R<3R$. Equation~\eqref{eq:localization-support-strip}
therefore gives
\[
 \supp T\cap\bigl(B_R^P(0)\times
 \{Y\in P^\perp:R\leq|Y|\leq2R\}\bigr)=\varnothing.
\]
The region in which we retain $T$ is
$B_R^P(0)\times B_{2R}^{P^\perp}(0)$, which lies in $B_{3R}(0)$.
By \eqref{eq:localization-support-strip}, the support of $T$ in this
region actually satisfies $|Y|<R/2$. Thus, within the cylinder
$B_R^P(0)\times P^\perp$, there is a neighborhood of the surface $|Y|=2R$ containing no points of $\supp T$.
Set
\[
 \mathcal C:=B_R^P(0)\times P^\perp,
 \qquad D:=B_R^P(0)\times B_{2R}^{P^\perp}(0).
\]
Using the representation of $T$ in Definition~\ref{def:integral-current}, define the
\defterm{retained current} $\widetilde T$ in $\mathcal C$ by
\[
 \widetilde T(\omega):=
 \int_{E\cap D}\theta(W)\langle\omega(W),\xi(W)\rangle
 \,\dd\mathcal H^k(W)
\]
for every smooth compactly supported $k$-form $\omega$ in $\mathcal C$.

We now show that $\partial\widetilde T=0$ in $\mathcal C$.
Choose a smooth function $\chi:P^\perp\to[0,1]$ equal to one on
$B_R^{P^\perp}(0)$ and compactly supported in $B_{2R}^{P^\perp}(0)$.
We also write $\chi$ for the function $X+Y\mapsto\chi(Y)$ on
$\mathcal C$. Its differential is supported in
\[
 B_R^P(0)\times\{Y\in P^\perp:R\leq|Y|\leq2R\},
\]
which contains no points of $\supp T$.
For a smooth compactly supported $(k-1)$-form $\beta$ in $\mathcal C$,
the form $\chi\beta$ has compact support in
$B_R^P(0)\times B_{2R}^{P^\perp}(0)\subset B_{3R}(0)\subset\mathcal U$.
We extend it by zero to $\mathcal U$ when applying $T$ or $\partial T$.
Consequently,
\begin{align*}
 \partial\widetilde T(\beta)
 &=\widetilde T(\dd\beta)=T(\chi\,\dd\beta)\\
 &=T(\dd(\chi\beta))-T(\dd\chi\wedge\beta)\\
 &=\partial T(\chi\beta)-T(\dd\chi\wedge\beta)=0.
\end{align*}
The first term vanishes because the original boundary is absent in
$B_{3R}(0)$. The second vanishes because $T$ has no support where
$\dd\chi\ne0$ inside the cylinder.

The restriction is integer rectifiable with finite mass, and its boundary
in $\mathcal C$ is zero, so it is an integral current there. Moreover,
\[
 \mathbf M(\widetilde T)\leq\|T\|(B_{3R}(0))<\infty,
 \qquad\supp_{\mathcal C}\widetilde T
 \subset\{X+Y\in\mathcal C:|Y|\leq R/2\},
\]
where $\supp_{\mathcal C}$ denotes support relative to $\mathcal C$.
It agrees with $T$ in $B_{R/2}(0)$ and hence has the same tangent cones at $0$.

Finally, we show that $\widetilde T$ is locally area minimizing in
$\mathcal C$. Let $W\Subset\mathcal C$ be open, let $S$ be an integral
$(k+1)$-current with $\supp S\Subset W$, and set
$A:=\widetilde T+\partial S$. The support of $S$ need not lie in
$\mathcal U$, so we first return it to the region where the original
current is defined. Define
\[
 p_R(Y):=\begin{cases}
 Y,&|Y|\leq R,\\
 RY/|Y|,&|Y|>R,
 \end{cases}
 \qquad F(X+Y):=X+p_R(Y)\quad(X\in P,\ Y\in P^\perp).
\]

\begingroup\color{black}
The map $p_R$ is nearest-point projection onto the closed convex ball
$\overline{B_R^{P^\perp}(0)}$.

Since $p_R(Y)$ minimizes the distance to $Y$ within the ball,
differentiating $|Y-p_R(Y)-t(V-p_R(Y))|^2$ at $t=0^+$ gives
\[
 \langle Y-p_R(Y),V-p_R(Y)\rangle\leq0
 \quad\bigl(V\in\overline{B_R^{P^\perp}(0)}\bigr).
\]
Apply this with $V=p_R(Y')$ and interchange $Y,Y'$ to obtain
\[
 |p_R(Y)-p_R(Y')|^2
 \leq\langle p_R(Y)-p_R(Y'),Y-Y'\rangle
 \leq |p_R(Y)-p_R(Y')|\,|Y-Y'|.
\]
Thus $p_R$ is $1$-Lipschitz. The orthogonal splitting
$P\oplus P^\perp$ gives the same bound for $F$.
\endgroup
Since $F$ is the identity on $|Y|<R$ and
$\supp_{\mathcal C}\widetilde T\subset\{|Y|\leq R/2\}$,
we have $F_\#\widetilde T=\widetilde T$.
Put $S':=F_\#S$. Pushforward commutes with boundary, preserves integral
currents, and satisfies the mass estimate for Lipschitz maps
\cite[Chapter~6, \S\S2--3]{SimonGMT}. Thus
\[
 F_\#A=\widetilde T+\partial S',
 \qquad\mathbf M(F_\#A)\leq\mathbf M(A).
\]
These pushforwards are well-defined even though $F$ is not proper on
the entire cylinder: $S$ has compact support, while $A$ has bounded
normal support. Since $F$ leaves $X$ unchanged, its restriction to
either support has compact inverse images of compact subsets of
$\mathcal C$.

Choose $0<a<a'<R$ such that $|X|<a$ on $\supp S$. Then
\[
 \supp S'\subset\{X+Y:|X|<a,\ |Y|\leq R\}
 \Subset O:=B_{a'}^P(0)\times B_{3R/2}^{P^\perp}(0).
\]
Here $O\Subset\mathcal U$ because $\overline O\subset B_{3R}(0)$,
and $T=\widetilde T$ in $O$ because $O\subset D$.
We extend the compactly supported current $S'$ by zero to
$\mathcal U$. Local area minimization of $T$ now gives
\begin{align*}
 \mathbf M(\widetilde T\llcorner O)
 &=\mathbf M(T\llcorner O)\\
 &\leq\mathbf M((T+\partial S')\llcorner O)\\
 &=\mathbf M((\widetilde T+\partial S')\llcorner O).
\end{align*}
Outside $O$, the currents $\widetilde T$ and
$\widetilde T+\partial S'$ agree. Adding their equal masses on
$\mathcal C\setminus O$ gives
\[
 \mathbf M(\widetilde T)
 \leq\mathbf M(\widetilde T+\partial S')
 =\mathbf M(F_\#A)
 \leq\mathbf M(A)
 =\mathbf M(\widetilde T+\partial S).
\]
All these masses in $\mathcal C$ are finite, since $\widetilde T$ has
finite mass and $S$ is an integral current with compact support.
Finally, $A=\widetilde T$ outside $W$, so subtracting their equal
masses on $\mathcal C\setminus W$ yields
\[
 \mathbf M(\widetilde T\llcorner W)
 \leq\mathbf M((\widetilde T+\partial S)\llcorner W).
\]
This proves local area minimization in $\mathcal C$ and completes the proof.

\end{proof}

\subsection{Planar frequency and tangent functions}\label{sec:planar-frequency}

We follow the planar-frequency and tangent-function framework of
Krummel--Wickramasekera \cite[\S3]{KWIII}. In particular, the cylinder,
projected radius, cutoff, and height, tilt, and frequency quantities below
are those of \cite[Definition~3.1, pp.~22--23]{KWIII}, with our notation.
We also recall the standard definitions of associated varifolds and
$q$-valued functions needed to state this framework.

Fix $Z\in\R^{k+m}$, $\rho_0>0$, and a $k$-plane $P$ through the origin.
Let $\pi_P$ denote orthogonal projection onto $P$.  The \defterm{projected
cylinder of radius $\rho_0$ about $Z$ relative to $P$} is
\[
 \mathbf C_{\rho_0}(Z,P)
 :=\{Z+X+Y:X\in P,\ |X|<\rho_0,\ Y\in P^\perp\}.
\]
To apply this setup at a flat singular point $Z$ of the original current,
first use Lemma~\ref{lem:localization} to obtain a localized current in
$\mathbf C_R(Z,P)$, where $P:=P_Z$. It is locally area minimizing,
has no boundary in this cylinder, and has bounded normal support.
It agrees with the original current near $Z$, so the local singularity
and tangent cones at $Z$ are unchanged. We take $\rho_0=R$ and continue
to denote the localized current by $T$.

Thus the cylindrical setup is as follows: $T$ is a
$k$-dimensional locally area-minimizing integral current
in this cylinder, with
$(\partial T)\llcorner\mathbf C_{\rho_0}(Z,P)=0$ and
\[
 \sup_{W\in\supp T\cap\mathbf C_{\rho_0}(Z,P)}
 \dist(W,Z+P)<\infty.
\]
The \defterm{projected radius from $Z$ relative to $P$} is
\[
 \mathcal R_{P,Z}(W):=|\pi_P(W-Z)|.
\]
Define the \defterm{cutoff function} $\vartheta:[0,\infty)\to[0,1]$ by
\[
 \vartheta(s):=\begin{cases}
 1,&0\leq s<1/2,\\
 2-2s,&1/2\leq s<1,\\
 0,&s\geq1.
 \end{cases}
\]
Set $\vartheta'(1/2):=\vartheta'(1):=0$.

The \defterm{Grassmann bundle} over the cylinder is the set of pairs
$(W,S)$, where $W\in\mathbf C_{\rho_0}(Z,P)$ and $S$ is an unoriented
$k$-dimensional linear subspace of $\mathbb R^{k+m}$.
Write $T$ in the integral representation of Definition~\ref{def:integral-current}, with
rectifiable set $E$, multiplicity $\theta$, and approximate tangent plane
$P_W$ at $\mathcal H^k$-almost every $W\in E$.
The \defterm{varifold associated with $T$} is the nonnegative Radon measure
$|T|$ on this bundle defined by
\[
 \int f(W,S)\,\dd|T|(W,S)
 :=\int_{E\cap\mathbf C_{\rho_0}(Z,P)}
 f(W,P_W)\,\theta(W)\,\dd\mathcal H^k(W)
\]
for every nonnegative Borel function $f$ on the bundle
\cite[Chapter~8, \S1]{SimonGMT}.
Thus the mass at $W$ is assigned to the pair $(W,P_W)$. In particular, $S=P_W$ at $|T|$-almost every pair $(W,S)$. The integrals below are integrals of scalar functions against this
measure, so they require no orientation.

Let $\pi_S$ be orthogonal projection onto $S$.
{\color{black}Where
$\mathcal R_{P,Z}(W)>0$, the ordinary gradient is defined, and we define}
the \defterm{tangential gradient of $\mathcal R_{P,Z}$ along $S$}
by
\[
 \nabla_S\mathcal R_{P,Z}:=\pi_S(\nabla\mathcal R_{P,Z}).
\]
The symbol $\|\pi_S-\pi_P\|$ denotes the Hilbert--Schmidt norm of
the difference of the two projection matrices.

We include the formulas below to define planar frequency explicitly.
Our proofs use its cited properties, not calculations with these integrals.
For $0<\rho\leq\rho_0$, define the \defterm{height quantity}
$H_{T,P,Z}(\rho)$, the \defterm{tilt quantity} $D_{T,P,Z}(\rho)$,
and the \defterm{frequency ratio} $N_{T,P,Z}(\rho)$ by
\begin{align*}
 H_{T,P,Z}(\rho)&:=-\rho^{1-k}\int
 \dist^2(W,Z+P)|\nabla_S\mathcal R_{P,Z}|^2
 \frac{1}{\mathcal R_{P,Z}(W)}
 \vartheta'(\mathcal R_{P,Z}(W)/\rho)
 \,\dd|T|(W,S),\\
 D_{T,P,Z}(\rho)&:=\frac{1}{2}\rho^{2-k}\int\|\pi_S-\pi_P\|^2
 \vartheta(\mathcal R_{P,Z}(W)/\rho)\,\dd|T|(W,S),\\
 N_{T,P,Z}(\rho)&:=\frac{D_{T,P,Z}(\rho)}{H_{T,P,Z}(\rho)}
 \quad\text{when }H_{T,P,Z}(\rho)>0.
\end{align*}
The integrand defining $H_{T,P,Z}$ is set equal to zero where
$\mathcal R_{P,Z}=0$. This is consistent with the fact that $\vartheta'$ is
zero near the origin.

The height quantity $H_{T,P,Z}(\rho)$ integrates the squared distance
from each point $W$ of the current's support to the affine plane $Z+P$,
with the weights in the formula above and over the region
$\rho/2<\mathcal R_{P,Z}(W)<\rho$.
The tilt quantity $D_{T,P,Z}(\rho)$ measures the weighted deviation of
the current's tangent planes from $P$ within the cylinder.
Under the hypotheses below, the limit of their ratio $N_{T,P,Z}(\rho)$
gives the homogeneity degree of the normalized tangent functions.

At a branch point $Z$ whose tangent cone is uniquely
$q\llbracket P_Z\rrbracket$, the \defterm{planar frequency of $T$ at $Z$} is,
when the limit exists,
\[
 N_{T,\mathrm{Pl}}(Z):=\lim_{\rho\downarrow0}N_{T,P_Z,Z}(\rho).
\]
We now justify existence of this limit for $Z\in\mathcal R_q(T)$.
Fix the exponent $\mu$ in Definition~\ref{def:rapid}, and apply
Lemma~\ref{lem:localization}. In the cylindrical formulas below, $T$
denotes the retained current. It is locally area minimizing and has
no boundary in a cylinder about $Z$, and it agrees with the original
current near $Z$. Define its \defterm{cylindrical height excess} by
\begin{equation}\label{eq:cylindrical-height-excess}
 e_Z(r):=
 \left(\frac{1}{\omega_k r^{k+2}}
 \int_{\mathbf C_r(Z,P_Z)}
       \dist^2(W,Z+P_Z)\,\dd\|T\|(W)\right)^{1/2}.
\end{equation}
This is the quantity in \cite[(2.1), p.~21]{KWIII}.
The estimate \eqref{eq:rapid} concerns ambient balls, whereas the
frequency theorem requires decay of $e_Z$ on cylinders. The support estimate \eqref{eq:localization-support-strip} bridges
these two conditions: translate $Z$ to $0$ and write
$W=X+Y$, with $X\in P_Z$ and $Y\in P_Z^\perp$.
Choose the localization radius $R$ so that the support estimate holds
at every radius at most $R$. On the retained support we have
$|X|<R$ and $|Y|\leq R/2$, so $|W|/3<R$.
For each nonzero retained support point, apply the estimate to the
original current at radius $|W|/3$. It gives $|Y|<|W|/6$, and hence
$|W|<6|X|/\sqrt{35}$. Consequently,
\begin{equation}\label{eq:localized-cylinder-ball}
 \supp T\cap\mathbf C_{7\rho/4}(Z,P_Z)
 \subset B_{2\rho}(Z)
 \qquad\text{for all sufficiently small }\rho>0,
\end{equation}
where we used $(6/\sqrt{35})(7/4)<2$ and returned to the original
coordinates. The mass measure of the retained current is the restriction
of the original current's mass measure to the retained region
\[
 Z+\bigl(B_R^{P_Z}(0)\times B_{2R}^{P_Z^\perp}(0)\bigr).
\]
Therefore \eqref{eq:rapid} gives
\begin{equation}\label{eq:cylindrical-height-decay}
 e_Z(7\rho/4)^2
 \leq \frac{C(2\rho)^{k+2+2\mu}}{\omega_k(7\rho/4)^{k+2}}
 \leq C'\rho^{2\mu},
\end{equation}
with $C'$ independent of sufficiently small $\rho$.
Together with local area minimization, absence of boundary, and density
$q$ at $Z$, this verifies condition (3.4) of
Krummel--Wickramasekera. Their Corollary~3.3 gives
\begin{equation}\label{eq:planar-frequency-lower-bound}
 1<1+\mu\leq N_{T,\mathrm{Pl}}(Z)<\infty
\end{equation}
and existence of the defining limit
\cite[Corollary~3.3, p.~24]{KWIII}.

A branch point $Z$ is a \defterm{frequency-$2$ branch point} if
$N_{T,\mathrm{Pl}}(Z)=2$. The stratification result
\cite[Theorem~3.7(i),(iii), pp.~31--32]{KWIII} and the
almost-everywhere existence and uniqueness result for
cylindrical normalized tangent functions
\cite[Theorem~1.8(ii)--(iv), pp.~15--16]{KWIII}
used below assume $N_{T,\mathrm{Pl}}(Z)\ne2$.
Separately, the local graph obstruction in dimension three
(Theorem~\ref{thm:local}) requires the leading degree to differ from $2$:
Example~\ref{ex:boruvka-quadratic} gives a conformally flat minimal
graph with a nonzero quadratic leading term.

We next recall $q$-valued functions and define normalized tangent functions.

\begin{definition}[$q$-valued functions]\label{def:q-valued}
For a finite-dimensional real inner-product space $V$ and an integer
$q\geq1$, define
\[
 \mathcal A_q(V):=V^q/\operatorname{Sym}(q),\qquad
 \sum_{l=1}^q\llbracket w_l\rrbracket\in\mathcal A_q(V),
\]
where $\operatorname{Sym}(q)$ is the group of permutations of $\{1,\ldots,q\}$.
An element $W$ of $\mathcal A_q(V)$ is an
\defterm{unordered $q$-tuple with multiplicity}.
The symbol $\llbracket w\rrbracket$ records one occurrence of the value
$w$. It is distinct from the notation for an integration current.
The \defterm{distance between unordered tuples} and their \defterm{norm} are
\[
 \mathcal G\!\left(\sum_l\llbracket w_l\rrbracket,
                    \sum_l\llbracket w'_l\rrbracket\right)
 :=\min_{\sigma\in\operatorname{Sym}(q)}
       \left(\sum_l|w_l-w'_{\sigma(l)}|^2\right)^{1/2},
 \qquad |W|:=\mathcal G(W,q\llbracket0\rrbracket).
\]
A \defterm{$q$-valued function} is a map into $\mathcal A_q(V)$.
Its Lipschitz constant is defined using $\mathcal G$.
Strong $L^2$ convergence $F_i\to F$ on $B$ means
$\int_B\mathcal G(F_i,F)^2\to0$.

The \defterm{sum} of tuples is their concatenation, so it adds their numbers
of values. A \defterm{repetition} $aW$, for an integer $a\geq1$, denotes the
concatenation of $a$ identical copies.
For $s\in\R$, \defterm{scalar multiplication} acts on
each value: $s\cdot W:=\sum_l\llbracket sw_l\rrbracket$.
Thus $m_jF_j$ in a component decomposition denotes repetition,
whereas $\varepsilon_i^{-1}\cdot u_i$ denotes scalar multiplication and
does not change the number of values. The \defterm{average} is
$\operatorname{av}(W):=q^{-1}\sum_lw_l \in V$.
A function is \defterm{nonzero} if it is not identically $q\llbracket0\rrbracket$.
It is \defterm{homogeneous of degree $\nu$} if
$F(sX)=s^\nu\cdot F(X)$ for every $s>0$.

{\color{black}Let $B$ be a bounded open ball in a Euclidean
space, and let $F:B\to\mathcal A_q(V)$ be measurable.}
The \defterm{scalar Sobolev space} $W^{1,2}(B)$ consists of
square-integrable real-valued functions on $B$ with square-integrable
distributional first derivatives.
We define the \defterm{Sobolev space} $W^{1,2}(B,\mathcal A_q(V))$
as follows: $F\in
W^{1,2}(B,\mathcal A_q(V))$ if $|F|\in L^2(B)$,
{\color{black}for every fixed $Q\in\mathcal A_q(V)$ the
scalar function $\mathcal G(F,Q)$ belongs to $W^{1,2}(B)$, and for
each coordinate $i$ there is a nonnegative $g_i\in L^2(B)$,
independent of $Q$, such that
$|\partial_i\mathcal G(F,Q)|\leq g_i$ almost everywhere.}
For a countable dense set $\mathcal D\subset\mathcal A_q(V)$ put
\[
 |\partial_iF|:=\mathop{\operatorname{ess\,sup}}_{Q\in\mathcal D}
      |\partial_i\mathcal G(F,Q)|,\qquad
 |dF|^2:=\sum_i|\partial_iF|^2.
\]
Here $|\partial_iF|$ denotes the nonnegative scalar function defined
above, not the norm of a separately defined vector-valued derivative
$\partial_iF$.
One can check that this definition is independent of $\mathcal D$.
{\color{black}The \defterm{Dirichlet energy} is
$\operatorname{Dir}(F;B):=\int_B|dF|^2$.
On an open domain, $F\in W^{1,2}_{\mathrm{loc}}$ means that its
restriction to each relatively compact ball belongs to $W^{1,2}$.}
{\color{black}A function
$F\in W^{1,2}_{\mathrm{loc}}$} is \defterm{locally Dirichlet minimizing} if on every relatively
compact ball its energy is no greater than that of any $W^{1,2}$
competitor with the same Sobolev boundary trace
 \cite[Definition~0.7 and Proposition~2.10]{DLSQ}.
   These conventions are those used in \cite[\S2 and \S3.2]{KWIII}.

Let $E,V$ be finite-dimensional real inner-product spaces and let
$B\subset E$ be open. The \defterm{graph} of a single-valued function
$f:B\to V$ is
\[
 \operatorname{graph}(f):=\{(X,f(X)):X\in B\}\subset E\times V.
\]
We call $f$ the \defterm{height function}. For a $q$-valued function $f:B\to\mathcal A_q(V)$,
the definition becomes
\[
 \operatorname{graph}(f):=
 \{(X,Y):X\in B,\ Y\text{ is a value of }f(X)\}\subset E\times V.
\]
In the current setting, take $E=P_Z$ and $V=P_Z^\perp$ and identify
$(X,Y)$ with $X+Y\in\R^{k+m}$. The translated graph is then
$Z+\operatorname{graph}(f)$. For single-valued $f$, its height function is still $f$. If $f$ is Lipschitz, its graph is rectifiable. After choosing an
orientation of $E$, the \defterm{graph current} $\mathbf G_f$ is the current of integration
over this graph: each value is counted with its repetition multiplicity,
and each approximate tangent plane is oriented so that projection onto
the oriented space $E$ preserves orientation. {\color{black}Rectifiability and the graph-current construction
for Lipschitz $q$-valued functions are recalled in
\cite[\S3.5, p.~29]{KWI}.}
\end{definition}

\paragraph{Construction of normalized tangent functions.}
Let $Z\in\mathcal R_q(T)$, and use the localized current from
Lemma~\ref{lem:localization} and the cylindrical height excess $e_Z$
defined in \eqref{eq:cylindrical-height-excess}.

{\color{black}We use the construction of
\cite[\S3.2, pp.~25--27]{KWIII}, checking how its estimates apply
in our notation.}
Choose $\rho_i\downarrow0$ and put
\begin{equation}\label{eq:tangent-scales}
 T_i:=T_{Z,\rho_i},\qquad \varepsilon_i:=e_Z(\rho_i).
\end{equation}
{\color{black}{We claim that $\varepsilon_i>0$ for all sufficiently
large $i$. If $e_Z(\rho_i)=0$,} the nonnegative integrand in
\eqref{eq:cylindrical-height-excess} forces the support of the current
in that cylinder to lie in $Z+P_Z$. Since the current has no boundary
there, the constancy theorem on this plane
\cite[Chapter~6, Theorem~2.41]{SimonGMT} makes it a constant integer
multiple of the plane near $Z$.   This would imply $Z\in\reg T$,
contrary to $Z\in\mathcal R_q(T)$. Thus the divisions by
$\varepsilon_i$ below are well-defined. Moreover, $\varepsilon_i\to0$
by the cylindrical decay estimate {\color{black}\eqref{eq:cylindrical-height-decay}}.}
If the localized current is defined in $\mathbf C_R(Z,P_Z)$, then
$T_i$ is defined in $\mathbf C_{R/\rho_i}(0,P_Z)$.
Fix a base ball $B\subset P_Z$ and choose $\sigma\geq1$ so that
$\overline B\subset B_\sigma(0,P_Z)$, where
$B_\sigma(0,P_Z):=\{X\in P_Z:|X|<\sigma\}$.
Here $R$ is the localization radius in the original coordinates,
whereas $\sigma$ is fixed in the magnified coordinates.
Since $R/\rho_i\to\infty$, the domain of $T_i$ contains
$\mathbf C_{4\sigma}(0,P_Z)$ for all sufficiently large $i$.
In particular,
\[
 \mathbf C_\sigma(0,P_Z)=B_\sigma(0,P_Z)\times P_Z^\perp.
\]
This cylinder need not equal $B\times P_Z^\perp$.

For all sufficiently large $i$, the construction in
\cite[(3.6)--(3.9), p.~25]{KWIII} provides a Lipschitz $q$-valued
function $u_{\sigma,i}:B_\sigma(0,P_Z)\to\mathcal A_q(P_Z^\perp)$
and a relatively closed set $K_{\sigma,i}\subset B_\sigma(0,P_Z)$. Set
\[
 u_i:=u_{\sigma,i}|_B,\qquad K_i:=K_{\sigma,i}\cap B.
\]
Write $\mathbf G_{u_i}$ for the graph current associated with
$u_i:B\to\mathcal A_q(P_Z^\perp)$, as defined in
Definition~\ref{def:q-valued}. Thus $K_i$ is closed relative to $B$. Restricting the equality of
currents in \cite[(3.6), p.~25]{KWIII} to
$K_i\times P_Z^\perp$ gives
\begin{equation}\label{eq:tangent-graph-agreement}
 \mathbf G_{u_i}\llcorner(K_i\times P_Z^\perp)
 =T_i\llcorner(K_i\times P_Z^\perp).
\end{equation}
This equality means that the graph current and $T_i$ agree above $K_i$,
including their multiplicities and orientations. Their omitted-region
estimate concerns
\[
 \mathbf C_\sigma(0,P_Z)\setminus\pi_{P_Z}^{-1}(K_{\sigma,i})
 =\bigl(B_\sigma(0,P_Z)\setminus K_{\sigma,i}\bigr)\times P_Z^\perp.
\]
Our set $(B\setminus K_i)\times P_Z^\perp$ is contained in this set.
Restricting their estimates to $B$ therefore gives
\begin{equation}\label{eq:tangent-approximation-bounds}
 \begin{gathered}
 |B\setminus K_i|+
 \|T_i\|((B\setminus K_i)\times P_Z^\perp)
       =O(\varepsilon_i^{2+2\varsigma}),\\
 \sup_B|u_i|=O(\varepsilon_i),\qquad
 \operatorname{Lip}(u_i)=O(\varepsilon_i^{2\varsigma})
 \end{gathered}
\end{equation}
for a fixed $\varsigma>0$, where $|B\setminus K_i|$ is Lebesgue measure.
The constants in these $O$-bounds may depend on the fixed $\sigma$,
but not on $i$.
For the Lipschitz bound, \cite[(3.9), p.~25]{KWIII} gives
$|du_i|\leq C_\sigma\varepsilon_i^{2\varsigma}$ almost everywhere.
Each scalar distance function $X\mapsto\mathcal G(u_i(X),Q)$,
with $Q\in\mathcal A_q(P_Z^\perp)$ fixed, has weak gradient bounded
by $|du_i|$. {We claim the following elementary estimate: if a real-valued}
Lipschitz function $h$ on an open convex set $B$ satisfies
$|\nabla h|\leq M$ almost everywhere, then
$|h(X)-h(X')|\leq M|X-X'|$ for every $X,X'\in B$.
For smooth $h$, this follows by integrating its gradient along the
segment joining $X'$ to $X$. For Lipschitz $h$, smooth by convolution
in a neighborhood of that segment. The gradient bound is preserved,
and passing to the limit gives the same inequality for $h$.
Applying this fact to $h(X):=\mathcal G(u_i(X),Q)$ with
$M:=C_\sigma\varepsilon_i^{2\varsigma}$ gives
\[
 |\mathcal G(u_i(X),Q)-\mathcal G(u_i(X'),Q)|
 \leq C_\sigma\varepsilon_i^{2\varsigma}|X-X'|
 \qquad(X,X'\in B).
\]
For each pair, take $Q=u_i(X')$ to obtain the stated bound on
$\operatorname{Lip}(u_i)$.

By \cite[(3.10)--(3.11), p.~25, and Theorem~3.4(i), p.~26]{KWIII},
after passing to a subsequence and relabelling, there exists
$\widehat\Phi\in W^{1,2}(B,\mathcal A_q(P_Z^\perp))$ such that
\begin{equation}\label{eq:tangent-convergence}
 \begin{gathered}
 \varepsilon_i^{-1}\cdot u_i\longrightarrow\widehat\Phi
       \quad\hbox{strongly in }L^2(B)\hbox{ and almost everywhere},\\
 \varepsilon_i^{-1}|du_i|\longrightarrow|d\widehat\Phi|
       \quad\hbox{strongly in }L^2(B).
 \end{gathered}
\end{equation}
After diagonal subsequence selection, convergence holds along one
common subsequence on all balls centred at $0$ with positive integer
radii. Along this subsequence, the limits agree on overlaps because the
approximations coincide outside sets of vanishing measure and their
normalized heights are uniformly bounded. These compatible limits
therefore define one function $\widehat\Phi$ on all of $P_Z$.

\begingroup\color{black}
\begin{definition}[Normalized tangent function]\label{def:tangent-function}
A \defterm{normalized tangent function of $T$ at $Z$} is any function
$\widehat\Phi:P_Z\to\mathcal A_q(P_Z^\perp)$ obtained as a limit by
the preceding construction, for some sequence of radii tending to zero
and a common subsequence on all bounded base balls.
\end{definition}
\endgroup

The normalized tangent functions constructed above have the following
additional properties by \cite[Theorem~3.4(ii), pp.~26--27]{KWIII}:
each is nonzero,
homogeneous of degree $\nu:=N_{T,\mathrm{Pl}}(Z)>1$
by \eqref{eq:planar-frequency-lower-bound}, and locally
Dirichlet minimizing, with
\[
 \frac{1}{\omega_k}\int_{B_1^{P_Z}(0)}|\widehat\Phi|^2=1.
\]
Note that at the same fixed point $Z$, a different sequence or
subsequence of scales may produce a different normalized tangent function.

{\color{black}Each normalized tangent function, viewed as a
Sobolev equivalence class, has a unique locally H\"older-continuous
representative. Existence follows from \cite[Theorem~0.9]{DLSQ}.
Uniqueness follows because two continuous representatives of the same
Sobolev class agree almost everywhere and hence everywhere.
We always use this unique representative, including when stating
homogeneity and defining the spine.}

For a homogeneous function $F:P\to\mathcal A_q(P^\perp)$ define its
\defterm{spine} by
\[
 \operatorname{spine}F
 :=\{a\in P:F(X+a)=F(X)\text{ for every }X\in P\}.
\]
For the homogeneous functions used here it is a linear
subspace. {\color{black}The defining translation invariance is
closed under addition and negatives. If $a$ belongs to the spine,
homogeneity gives, for $s>0$,
$F(X+sa)=s^\nu\cdot F(X/s+a)=s^\nu\cdot F(X/s)=F(X)$,
so the spine is also closed under scalar multiplication.}
We call $F$ \defterm{cylindrical} if its spine has dimension $k-2$.
For a nonzero homogeneous Dirichlet minimizer $F$ of degree $\nu>0$,
  let $N_F(a)$ denote its Dirichlet frequency at $a\in P$,
as defined in \cite[(3.12)--(3.14), pp.~25--26]{KWIII}.
Homogeneity gives  $N_F(0)=\nu$. Krummel--Wickramasekera define the spine
by the frequency condition $N_F(a)=\nu$
\cite[\S3.3, p.~31]{KWIII}. The equivalence with our definition is
justified in Part~I, \cite[\S3.6, p.~30]{KWI}, where they state,
as a consequence of Almgren's Theorem~2.14, that
\[
 N_F(a)=\nu\quad\Longrightarrow\quad
 F(X+a)=F(X)\qquad\text{for every }X\in P.
\]

{\color{black}Conversely, translation invariance along $a$
identifies the frequency integrals centered at $a$ with those centered
at $0$, so $N_F(a)=N_F(0)=\nu$.}

\paragraph{Almost-everywhere uniqueness.}
We need this setup to apply Krummel--Wickramasekera's
Theorem~1.8 \cite[pp.~15--16]{KWIII}, which gives almost-everywhere uniqueness of tangent
functions and assumes Hypothesis~1.4.
Fix $q\geq2$. We first justify the unit-cylinder hypotheses near a
flat singular point of density $q$. Translate that point to $0$, and
let $T$ denote the retained current from Lemma~\ref{lem:localization}
in $\mathbf C_R(0,P)$. {We claim that the support-convergence argument in that lemma gives}
\begin{equation}\label{eq:vanishing-relative-normal-height}
 \sup_{W\in\supp T\cap B_{2r}(0)}\dist(W,P)=o(r).
\end{equation}
{For every $\epsilon>0$, that argument applies with strip width
$\epsilon$ and gives}
$\dist(W,P)<\epsilon r$ for every
$W\in\supp T\cap B_{2r}(0)$ once $r$ is sufficiently small
(depending on $\epsilon$). Thus the displayed supremum divided by
$r$ tends to zero, which is the meaning of $o(r)$.
The ball--cylinder comparison \eqref{eq:localized-cylinder-ball} puts the support in
$\mathbf C_r(0,P)$ inside $B_{2r}(0)$ for small $r$. Thus, after
magnification by $1/r$, its normal support in the unit cylinder is
arbitrarily small.

Projection onto $P$ is proper on the retained support because that
support is normally bounded. The projected current $\pi_{P\#}T$ has no boundary
in $B_R^P(0)$, so constancy gives
$\pi_{P\#}T=d\llbracket B_R^P(0)\rrbracket$ for an integer $d$
\cite[Chapter~6, Theorem~2.41]{SimonGMT}.
We now show that magnification leaves this integer unchanged, while
$T_r\to q\llbracket P\rrbracket$ forces $d=q$.
{Let $\alpha$ be a smooth compactly supported $k$-form on $P$.}
The support bound \eqref{eq:vanishing-relative-normal-height}, together
with \eqref{eq:localized-cylinder-ball}, allows us to choose a smooth compactly
supported function $\chi:P^\perp\to[0,1]$, independent of sufficiently
small $r$, equal to one near $0$ and at every normal coordinate of
$\supp T_r$ above $\supp\alpha$. The form
$\chi(Y)\pi_P^*\alpha$ is compactly supported, so current convergence gives
\begin{equation}\label{eq:projected-current-limit}
 (\pi_{P\#}T_r)(\alpha)
 =T_r\bigl(\chi(Y)\pi_P^*\alpha\bigr)
 \longrightarrow q\llbracket P\rrbracket(\alpha).
\end{equation}
Since $\alpha$ has compact support and $R/r\to\infty$ as $r\downarrow0$,
$\supp\alpha\subset B_{R/r}^P(0)$ for all sufficiently small $r$.
For every such $r$, the identity
$\pi_{P\#}T_r=d\llbracket B_{R/r}^P(0)\rrbracket$ shows that the
left-hand side of \eqref{eq:projected-current-limit} equals
$d\int_P\alpha$. Choosing $\alpha$ with
$\int_P\alpha=1$ gives $d=q$. Thus
\begin{equation}\label{eq:localized-projection-identity}
 (\pi_P)_\#T=q\llbracket B_R^P(0)\rrbracket.
\end{equation}

Consequently, restricting to a sufficiently small cylinder, rotating
$P$ to the coordinate plane, and rescaling verifies all of
\cite[Hypothesis~1.4, p.~14]{KWIII}: local area minimization,
absence of boundary, small normal support, and projection equal to
$q\llbracket B_1\rrbracket$.

{\color{black}In the following almost-everywhere statement
and normalization formula, we continue to write $T$ for this restricted,
rotated, and rescaled current. The point $Z$ and plane $P_Z$ are also
expressed in these transformed coordinates.}

At $\HH^{k-2}$-almost every
branch point $Z$ of density $q$ in
$\mathbf C_{1/2}(0,\mathbb R^k\times\{0\})
=B_{1/2}^k(0)\times\mathbb R^m$
with $\nu:=N_{T,\mathrm{Pl}}(Z)\ne2$, Theorem~1.8(ii)--(iv)
provides a unique nonzero, cylindrical, homogeneous degree-$\nu$
Dirichlet-minimizing profile
$\Phi:P_Z\to\mathcal A_q(P_Z^\perp)$ satisfying its approximation
and multiplicity conclusions \cite[pp.~15--16]{KWIII}.
We call this the \defterm{physical leading profile}.

At each such point $Z$, Krummel--Wickramasekera identify the unit
mean-square normalization of $\Phi$ as the unique normalized tangent
function obtained by the blow-up procedure recalled in
Definition~\ref{def:tangent-function}
\cite[proof of Theorem~11.1, Step~4, p.~143]{KWIII}. Therefore,
\[
 \widehat\Phi=
 \frac{\sqrt{\omega_k}}{\|\Phi\|_{L^2(B_1^{P_Z}(0))}}\cdot\Phi.
\]
Thus, at each point $Z$ covered by this almost-everywhere conclusion,
every sequence producing a normalized tangent function at that fixed
$Z$ gives the same $q$-valued function $\widehat\Phi$.
The formula above makes explicit the normalization distinction left
implicit in the notation of Theorem~1.8(iii)--(iv).

Our branch theorem, Theorem~\ref{thm:intro-branch}, uses $\widehat\Phi$.
Its pointwise conclusion in part~(i)
applies to every cylindrical normalized tangent function, whether or
not uniqueness holds. Since positive scalar multiplication preserves
component repetition multiplicities, the same multiplicity
conclusion holds for $\Phi$ at the points described here.

\subsection{Coordinates, cyclic components, and exact sheets}\label{sec:cyclic-components}
Let $Z\in\mathcal R_q(T)$, where $\mathcal R_q(T)$ is defined in Definition~\ref{def:rapid}. The unique tangent cone at $Z$ is $q\llbracket P_Z\rrbracket$.
Fix a cylindrical normalized tangent function $F:=\widehat\Phi$ at $Z$ from Definition~\ref{def:tangent-function}, so that
\[
 F:P_Z\longrightarrow\mathcal A_q(P_Z^\perp).
\]
Its \defterm{branch axis} is
\[
 \mathcal L_Z:=\operatorname{spine}F.
\]
Since $F$ is cylindrical, $\dim\mathcal L_Z=k-2$. Note that this branch axis does not have to be the branch set of the current.

Choose orthonormal coordinates identifying
\[
 P_Z=\R^2\times\R^{k-2},\qquad
 \mathcal L_Z=\{(0,0,t):t\in\R^{k-2}\}.
\]
Write $X=(x,y,t)$, where $t\in\R^{k-2}$, and put
\[
 z:=x+iy,\qquad r:=|z|.
\]
Thus $\mathcal L_Z=\{z=0\}$, and $F$ is independent of $t$.

The function $F$ is homogeneous of degree
$\nu=N_{T,\mathrm{Pl}}(Z)>1$. In this case, the classification of cylindrical homogeneous Dirichlet
minimizers \cite[\S3.3, pp.~32--33, and Definition~5.1 and Remark~5.2]{KWIII}
gives the following: $\nu$ is a rational number of the form $p_0/q_0$,  for some relatively prime positive integers $p_0, q_0$,
\[
 F=m_0\llbracket0\rrbracket+\sum_{j=1}^Jm_jF_j,\qquad
 q=m_0+q_0\sum_{j=1}^Jm_j,\qquad J\geq1,
\]
where $m_0\geq0$ and each $m_j\geq1$ is an integer. When $m_0=0$,
we do not write the zero term. The functions $F_j$ are distinct nonzero
$q_0$-valued functions, which we call \defterm{cyclic components}.
On an open set \(A\subset P_Z\setminus\mathcal L_Z\) admitting a continuous choice of \(\theta=\arg z\), the restriction $F_j|_A$ is the unordered collection of the following $q_0$ values:
\begin{equation}\label{eq:branches}
 \varphi_{j,l}(z,t):=
 \operatorname{Re}\!\left(c_j r^\nu e^{i\nu(\theta+2\pi l)}\right),
 \qquad l=0,\ldots,q_0-1,\qquad
 c_j\in P_Z^\perp\otimes\C\setminus\{0\}.
\end{equation}
The integer $q_0$ is the \defterm{cyclic sheet number}, and $m_j$ is the
\defterm{repetition multiplicity}. The \defterm{zero component} is
$\llbracket0\rrbracket$. For a complex normal vector $c=(c^\alpha)$ we use the
\defterm{complex-bilinear product}
\begin{equation}\label{eq:complex-bilinear-product}
 c\cdot c:=\sum_\alpha(c^\alpha)^2
\end{equation}
and the \defterm{Hermitian squared norm} $|c|^2:=\sum_\alpha|c^\alpha|^2$.

\begin{definition}[Annular sheet chart and sheet height]\label{def:chart}
Let $\Omega$ be an open subset of $\{X\in P_Z:|X|=1\}$ such that
$\overline\Omega\cap\mathcal L_Z=\varnothing$ and there is a continuous
function $\theta:\Omega\to\R$ satisfying $z=|z|e^{i\theta}$.
Define its \defterm{conical sector} by
\[
 C_\Omega:=\{sX:s>0,\ X\in\Omega\}.
\]
Extend the angle to $C_\Omega$ by $\theta(sX):=\theta(X)$.
An \defterm{annular sheet chart} is a nonempty open set
$A\Subset C_\Omega\cap\{X:1/2<|X|<2\}$.
For $\rho>0$, the set $\rho A:=\{\rho X:X\in A\}$ is the \defterm{dilation} of $A$ by $\rho$.
We call the functions $\varphi_{j,l}:C_\Omega\to P_Z^\perp$ defined by
\eqref{eq:branches} the \defterm{local branches} of $F_j$ on $C_\Omega$.

An \defterm{exact sheet-height function} on $\rho A$ is a smooth single-valued map
$v:\rho A\to P_Z^\perp$ whose
translated graph $Z+\operatorname{graph}(v)$
is an open subset of $\reg T$.
Lemma~\ref{lem:extract} establishes the existence of exact sheet-height
functions under the assumptions stated there.
\end{definition}

\subsection{The induced metric of a graph}\label{sec:graph-metric}
Let $E,V$ be finite-dimensional real inner-product spaces, let
$A\subset E$ be open, and let $w:A\to V$ be smooth.
We use the graph defined in Definition~\ref{def:q-valued}.
The map
\[
 \Gamma:A\to E\times V,\qquad \Gamma(X):=(X,w(X)),
\]
parametrizes $\operatorname{graph}(w)$.
Write $\delta_{\mathrm{amb}}$ for the product metric on $E\times V$
and $\delta$ for the metric of $E$, restricted to $A$.
The \defterm{induced metric of $\operatorname{graph}(w)$, expressed on $A$},
is $\Gamma^*\delta_{\mathrm{amb}}$.
For $a,b\in E$, identified with tangent vectors at $X\in A$,
$d\Gamma_X(a)=(a,dw_X(a))$. The Euclidean product metric gives
\[
 (\Gamma^*\delta_{\mathrm{amb}})_X(a,b)
 =\langle a,b\rangle+\langle dw_X(a),dw_X(b)\rangle.
\]
Thus
\begin{equation}\label{eq:induced-graph-metric}
 \Gamma^*\delta_{\mathrm{amb}}=\delta+(dw)^Tdw,
 \qquad ((dw)^Tdw)_{ab}:=\langle\partial_a w,\partial_b w\rangle
\end{equation}
in orthonormal coordinates on $E$ and $V$.
This applies with $E=\R^k$ and $V=\R^m$ in
Section~\ref{sec:local}. In the current setting, for
$A\subset P_Z$ and $w:A\to P_Z^\perp$, the parametrization
$X\mapsto Z+X+w(X)$ of $Z+\operatorname{graph}(w)$ induces the same
formula \eqref{eq:induced-graph-metric}, with $\delta$ the Euclidean
metric on $P_Z$ restricted to $A$.

\subsection{Normalized heights}\label{sec:normalized-heights}
We use two normalizations: division by $\rho^\nu$ for an exact
sheet-height function with a prescribed leading order $\nu$, and
division by the cylindrical height excess for an exact sheet-height
function in a magnified current. The first is convenient for explicit
examples and the fixed-function formulation of Theorem~\ref{thm:local}.
The proof of Theorem~\ref{thm:intro-branch} uses the second and does not
assume that $e_Z(\rho)/\rho^{\nu-1}$ converges to a finite positive number.

Fix an annular sheet chart $A$ as in Definition~\ref{def:chart}.

Let $\nu>1$ and let $v$ be a smooth single-valued function on
$\bigcup_{0<\rho<\rho_0}\rho A$ whose restriction to each $\rho A$
is an exact sheet-height function. Its \defterm{height normalized by $\rho^\nu$}
is
\[
 f_\rho:A\to P_Z^\perp,\qquad f_\rho(X):=\rho^{-\nu}v(\rho X).
\]
Put $\varepsilon_\rho:=\rho^{\nu-1}$.
The exact sheet-height function in the magnified current is
\[
 w_\rho(X):=\rho^{-1}v(\rho X)=\varepsilon_\rho f_\rho(X).
\]
By \eqref{eq:induced-graph-metric}, the induced metric of
$\operatorname{graph}(w_\rho)$, expressed on $A$, is
\[
 G_\rho:=\delta+(dw_\rho)^Tdw_\rho
        =\delta+\varepsilon_\rho^2(df_\rho)^Tdf_\rho.
\]
Writing $G:=\delta+(dv)^Tdv$ and $d_\rho(X):=\rho X$, we have
$G_\rho=\rho^{-2}d_\rho^*G$.
Thus $G_\rho$ is locally conformally flat on $A$ if and only if $G$
is locally conformally flat on $\rho A$.

For the magnified currents $T_i$ and numbers $\varepsilon_i=e_Z(\rho_i)$
in Definition~\ref{def:tangent-function}, suppose that an exact
sheet-height function $w_i:A\to P_Z^\perp$ for $T_i$ exists.
Its \defterm{height normalized by the cylindrical height excess}
is $f_i:A\to P_Z^\perp$, and the induced metric of
$\operatorname{graph}(w_i)$ is $G_i$, where
\[
 f_i:=\varepsilon_i^{-1}w_i,\qquad
 G_i:=\delta+\varepsilon_i^2(df_i)^Tdf_i.
\]
The corresponding exact sheet-height function for $T$ on $\rho_iA$ is
$v_i(X):=\rho_iw_i(X/\rho_i)$.

\subsection{A two-valued example}

Let $k=3$, $m=2$, and use coordinates $(z,t,w)\in\C\times\R\times\C = \mathbb R^5$.
Let $T$ be integration with multiplicity one   over the smooth submanifold
$\{(z,t,w)\in\mathbb R^5:(w-z^2)^2=z^3,\ (z,w)\ne(0,0)\}$, with the complex-curve
orientation determined by $(e, i e)$, for any nonzero tangent vector $e$ to the complex curve in the $(z,w)$ variables, followed by the positive $t$-direction. This current is calibrated
by the standard K\"ahler form wedged with $dt$, i.e. by $(dx \wedge dy + du \wedge dv) \wedge dt$, where $z = x+iy, w = u+iv$, and hence is locally
area minimizing. Put $P:=\{w=0\}$ and
$\mathcal L:=\{z=w=0\}$. At $Z=0$ one can check that the tangent cone is $2\llbracket P\rrbracket$.

On $(\C\setminus(-\infty,0])\times\R$, choose $\arg z\in(-\pi,\pi)$
and set $z^{3/2}:=|z|^{3/2}e^{(3i/2)\arg z}$.
The functions $v_\pm(z,t):=z^2\pm z^{3/2}$ restrict to exact
sheet-height functions on annular sheet charts in this domain.
The physical leading profile and its
normalized tangent function are
\[
 \Phi(z,t):=\llbracket z^{3/2}\rrbracket+
                  \llbracket-z^{3/2}\rrbracket,
 \qquad
 \widehat\Phi=\frac{4}{\sqrt{3\pi}}\cdot\Phi,
\]
since $\omega_3^{-1}\int_{B_1^P}|\Phi|^2=3\pi/16$. Both are homogeneous of degree $\nu = 3/2$.

For $\rho>0$, the magnification map and magnified current at $0$
from Definition~\ref{def:magnification} are
\[
 \eta_{0,\rho}(z,t,w):=(z/\rho,t/\rho,w/\rho),
 \qquad T_\rho:=(\eta_{0,\rho})_\#T.
\]
On an annular sheet chart $A$ in the domain chosen above, the exact
sheet-height functions for $T_\rho$ are
\[
 w_{\rho,\pm}(z,t):=\rho^{-1}v_\pm(\rho z,\rho t)
            =\rho z^2\pm\rho^{1/2}z^{3/2}.
\]
For $\nu=3/2$, the heights normalized by $\rho^\nu$ are
\[
 f_{\rho,\pm}(z,t):=\rho^{-3/2}v_\pm(\rho z,\rho t)
 =\sqrt\rho\,z^2\pm z^{3/2}
 \longrightarrow\pm z^{3/2}.
\]
The cylindrical height excess satisfies
$e_0(\rho)=(\sqrt{3\pi}/4)\rho^{1/2}(1+o(1))$.
For a sequence $\rho_i\downarrow0$, the heights normalized by the cylindrical height excess satisfy
\[
 \frac{w_{\rho_i,\pm}(z,t)}{e_0(\rho_i)}
 =\frac{\sqrt{\rho_i}}{e_0(\rho_i)}
       \bigl(\sqrt{\rho_i}\,z^2\pm z^{3/2}\bigr)
 \longrightarrow\pm\frac{4}{\sqrt{3\pi}}z^{3/2}.
\]
The limits $\pm z^{3/2}$ of $f_{\rho,\pm}$ are local branches of $\Phi$, and the second limits $\pm\frac{4}{\sqrt{3\pi}}z^{3/2}$ are
local branches of $\widehat\Phi$. Both convergences hold in every fixed
$C^b(A)$ norm. To summarize, the function $\widehat\Phi$ has homogeneous degree $\nu=3/2$,
one cyclic component, cyclic sheet number two, and repetition
multiplicity one. Theorem~\ref{thm:local} therefore rules out local
conformal flatness of the induced metric on $\reg T$ near the origin.

\subsection{Curvature conventions}\label{sec:curvature-conventions}

We recall the curvature objects on a smooth $k$-dimensional Riemannian
manifold $(L^k,g)$, with $k\geq2$. The Schouten, Cotton, and Weyl
tensors below are defined for $k\geq3$. Here $\nabla$ is the
Levi--Civita connection of $g$, and $X,Y,V,W$ are smooth vector fields
on $L$.
The \defterm{curvature operator} is
\[
 R^g(X,Y)V:=
 \nabla_X\nabla_YV-\nabla_Y\nabla_XV-\nabla_{[X,Y]}V,
\]
where   the \defterm{Lie bracket} is the vector field characterized by
$[X,Y]f=X(Yf)-Y(Xf)$ for every smooth scalar function $f$.
The \defterm{Riemann curvature tensor} is
\[
 \operatorname{Rm}_g(X,Y,V,W):=g\bigl(R^g(X,Y)V,W\bigr).
\]
Choose any smooth local frame $e_1,\ldots,e_k$ of $TL$. Write
\[
 g_{ij}:=g(e_i,e_j),\qquad
 (g^{ij}):=(g_{ij})^{-1},\qquad
 \operatorname{Rm}_{ijkl}(g):=\operatorname{Rm}_g(e_i,e_j,e_k,e_l).
\]
Repeated contraction indices are summed from $1$ to $k$.
The \defterm{Ricci tensor} and \defterm{scalar curvature} are defined by
\[
 \operatorname{Ric}_{jk}(g):=g^{il}\operatorname{Rm}_{ijkl}(g),
 \qquad
 \operatorname{Scal}_g:=g^{jk}\operatorname{Ric}_{jk}(g).
\]
In an orthonormal frame, $g_{ij}=g^{ij}=\delta_{ij}$.
For later calculations in local coordinates $(x^1,\ldots,x^k)$,
we take $e_i=\partial_i:=\partial/\partial x^i$.

For linearly independent vectors $v,w\in T_pL$, put
$\sigma:=\operatorname{span}\{v,w\}$. The \defterm{sectional curvature}
of $\sigma$ is
\begin{equation}\label{eq:sectional-curvature}
 K_g(\sigma)=K_g(v\wedge w):=
 \frac{\operatorname{Rm}_g(v,w,w,v)}{g(v,v)g(w,w)-g(v,w)^2}.
\end{equation}
This value is independent of the chosen basis $v,w$ of $\sigma$.
In particular, if $e_i,e_j$ are orthonormal and $i\ne j$, then
\[
 K_g(e_i\wedge e_j)=\operatorname{Rm}_{ijji}(g)
 \qquad\text{(no sum over $i,j$)}.
\]

For an immersion $F:L\to(M,h)$ with induced metric $g:=F^*h$ on $L$,
use any smooth local frame $e_1,\ldots,e_k$ on $L$ and write
\[
 E_i:=dF(e_i),\qquad
 \mathrm{II}_{ij}:=\mathrm{II}(e_i,e_j)
       =\bigl(\nabla^M_{e_i}(dF(e_j))\bigr)^\perp,
\]
where $\nabla^M$ here denotes the connection along $F$ induced by
the Levi--Civita connection of $h$, and $\perp$ is the
$h$-orthogonal projection onto the normal space.
With the curvature convention above, the \defterm{Gauss equation}
\cite[Theorem~8.5, p.~230]{LeeRM} is
\begin{equation}\label{eq:gauss-equation}
 \begin{aligned}
 \operatorname{Rm}_{ijkl}(g)
 &=h\bigl(R^h(E_i,E_j)E_k,E_l\bigr)\\
 &\quad+h(\mathrm{II}_{il},\mathrm{II}_{jk})
       -h(\mathrm{II}_{ik},\mathrm{II}_{jl}).
 \end{aligned}
\end{equation}
For a Euclidean ambient metric $h$, $R^h=0$, so the first term vanishes.

For $k\geq3$, the \defterm{Schouten tensor} is
\[
 S_g:=\frac{1}{k-2}\left(\operatorname{Ric}_g-
                   \frac{\operatorname{Scal}_g}{2(k-1)}g\right).
\]
We use $S_g$ rather than $P_g$ to distinguish it from the tangent plane $P$.
For symmetric covariant $2$-tensors $a,b$ on $L$, the
\defterm{Kulkarni--Nomizu product} is defined by
\[
 (a\owedge b)_{ijkl}:=
 a_{il}b_{jk}+a_{jk}b_{il}-a_{ik}b_{jl}-a_{jl}b_{ik}.
\]
Write $\nabla_k:=\nabla_{e_k}$ and
$\nabla_k(S_g)_{ij}:=(\nabla_{e_k}S_g)(e_i,e_j)$.
The \defterm{Cotton tensor} and \defterm{Weyl tensor} are
\[
 C_{ijk}(g):=(\nabla_{e_k}S_g)(e_i,e_j)
             -(\nabla_{e_j}S_g)(e_i,e_k),\qquad
 W(g):=\operatorname{Rm}_g-g\owedge S_g.
\]
Here the covariant derivatives act on the tensor $S_g$:
if $\nabla_{e_k}e_i=\Gamma^\ell_{ki}e_\ell$, then
\[
 (\nabla_{e_k}S_g)(e_i,e_j)
 =e_k\bigl((S_g)_{ij}\bigr)
  -\Gamma^\ell_{ki}(S_g)_{\ell j}
  -\Gamma^\ell_{kj}(S_g)_{i\ell}.
\]
In dimension three, local conformal flatness is equivalent to $C=0$.
In dimensions at least four it is equivalent to $W=0$
\cite[\S1.J]{Besse}.

\begin{samepage}
\subsection{The branch theorem}\label{sec:branch-statement}
With these definitions, we can state the main branch-obstruction theorem.
The ambient metric in this theorem is Euclidean.

\begin{theorem}[Branch obstruction and dimension bound]\label{thm:intro-branch}
Let $T$ be a $k$-dimensional locally area-minimizing integral current in an
open set $\mathcal U\subset\R^{k+m}$, with $k\geq3$, $m\geq2$, and no
boundary in $\mathcal U$. Suppose that the induced metric on $\reg T$
is locally conformally flat. For every integer $q\geq2$, the following
statements hold.
\begin{enumerate}[label=(\roman*),leftmargin=2em]
\item Let $Z\in\mathcal R_q(T)$ and
$\nu:=N_{T,\mathrm{Pl}}(Z)\ne2$. Every cylindrical normalized tangent
function
\[
 \widehat\Phi=m_0\llbracket0\rrbracket+
                     \sum_{j=1}^Jm_j\widehat\Phi_j
\]
has $m_j\geq2$ for every nonzero component. If its rational degree is
$\nu=p_0/q_0$ in lowest terms, then
\[
 2q_0J\leq q-m_0.
\]
For $k\geq4$ this pointwise statement also holds at $\nu=2$ whenever
the normalized tangent function is cylindrical.
\item Define
\[
 \mathcal E_q:=
 \{Z\in\mathcal R_q(T):\nu:=N_{T,\mathrm{Pl}}(Z)\notin\mathbb Q
 \ \text{or}\;
 \nu=p_0/q_0\notin\mathbb Z,\ \gcd(p_0,q_0)=1,\ 2q_0>q\}.
\]
Then $\dim_{\HH}\mathcal E_q\leq k-3$. In particular,
\[
 \dim_{\HH}\{Z\in\mathcal R_2(T)\cup\mathcal R_3(T):
                    N_{T,\mathrm{Pl}}(Z)\notin\mathbb Z\}\leq k-3.
\]
\item At $\HH^{k-2}$-almost every branch point of density $q$ and planar
frequency $\nu\ne2$, the unique normalized tangent function $\widehat\Phi$ has repetition
multiplicity $m_j \geq 2$ for every nonzero component and satisfies
$2q_0J\leq q-m_0$. In particular, the nonintegral-frequency
points of densities $2$ and $3$ have zero $\HH^{k-2}$ measure.
\end{enumerate}
\end{theorem}
\end{samepage}

\begin{remark}
We observe several things.
Part (i) applies to every cylindrical normalized tangent function at the
specified point, without asserting that all tangent functions are
cylindrical. Part (ii) follows by combining this conclusion with
the stratification discussed in
Section~\ref{sec:stratification}. Its dimension bound concerns only the rapid-decay set  $\mathcal R_q(T)$:
the remaining branch points form an $\HH^{k-2}$-null set, and we do not have a stronger dimension estimate for that.
Finally, when both $\Phi$ and $\widehat\Phi$ exist, they have the same cyclic sheet numbers and
repetition multiplicities, so their different scalar normalizations do
not affect the conclusions.
\end{remark}

\section{The local conformal-curvature obstruction}\label{sec:local}
This section is independent of the minimal-surface equation and of
current theory. We use the splitting
$\R^k=\R^2\times\R^{k-2}$, with coordinates $(z,t)$,
where $z=x+iy$, $r=|z|$, and $t\in\R^{k-2}$.
For smooth functions $f_i:A\to\R^m$ on an open set $A\subset\R^k$
and positive numbers $\varepsilon_i$, consider the graphs
$\operatorname{graph}(\varepsilon_i f_i)\subset\R^k\times\R^m$.
Their induced metrics, expressed on $A$, are
$G_i:=\delta+\varepsilon_i^2(df_i)^Tdf_i$, as in
Section~\ref{sec:graph-metric}. In Theorem~\ref{thm:local}, $A$ is an
annular sheet chart, $\varepsilon_i\to0$, and $f_i$ converge to
$\varphi(z,t):=\operatorname{Re}(cz^\nu)$, with the power defined using
the angle in Definition~\ref{def:chart}. We do not assume that $\varphi$
is a local branch of a cyclic component of a cylindrical normalized
tangent function or of a physical leading profile.
In the application to a cylindrical normalized tangent function,
the two factors of $\R^k$ are $\mathcal L_Z^\perp\cap P_Z$ and
$\mathcal L_Z$, of dimensions $2$ and $k-2$, respectively.
Indices $1,2$ denote $x,y$. Indices $\ell=3,\ldots,k$ denote the
coordinates in the second factor. For $k=3$, $t$ is a single real coordinate.

We first express the Riemann curvature tensor of the induced metric of
$\operatorname{graph}(\varepsilon f)$ in terms of the derivatives of $f$
and identify its limit after division by $\varepsilon^2$. This is a
standard consequence of the Gauss equation. We include the
calculation to fix conventions and track the derivatives required below.
\begin{lemma}[Gauss equation for a small graph]\label{lem:gauss}
Let $A\subset\R^k$ be open, let $f:A\to\R^m$ be smooth, and let
$\varepsilon>0$. The induced metric of $\operatorname{graph}(\varepsilon f)$,
expressed on $A$, is $G:=\delta+\varepsilon^2(df)^Tdf$. Define
$B_\varepsilon:=(I_m+\varepsilon^2df(df)^T)^{-1}$.
Then
\begin{equation}\label{eq:gauss-exact}
 \operatorname{Rm}_{ijkl}(G)
 =\varepsilon^2\bigl(
 \langle B_\varepsilon f_{il},f_{jk}\rangle
 -\langle B_\varepsilon f_{ik},f_{jl}\rangle\bigr).
\end{equation}
Consequently, let $\varphi:A\to\R^m$ be a $C^2$ function with
$\|\varphi\|_{C^2(A)}<\infty$, let
$f_\rho:A\to\R^m$ be smooth functions converging to $\varphi$ in
$C^2(A)$, and let $\varepsilon_\rho>0$ tend to zero.
For the induced metrics
$G_\rho:=\delta+\varepsilon_\rho^2(df_\rho)^Tdf_\rho$ of
$\operatorname{graph}(\varepsilon_\rho f_\rho)$, expressed on $A$, we have
\begin{equation}\label{eq:gauss-limit}
 \varepsilon_\rho^{-2}\operatorname{Rm}_{ijkl}(G_\rho)
 \longrightarrow
 \langle\varphi_{il},\varphi_{jk}\rangle
 -\langle\varphi_{ik},\varphi_{jl}\rangle.
\end{equation}
The convergence in \eqref{eq:gauss-limit} is in $C^0(A)$.
It is in $C^1(A)$ if $\varphi\in C^3(A,\R^m)$,
$\|\varphi\|_{C^3(A)}<\infty$, and $f_\rho\to\varphi$ in $C^3(A)$.

\end{lemma}
\begin{proof}
Let $e_1,\ldots,e_k$ be the standard basis vectors of $\R^k$.
For the graph immersion $F(X):=(X,\varepsilon f(X))$, the vectors
$dF_X(e_i)=(e_i,\varepsilon\partial_i f(X))$, $i=1,\ldots,k$,
form a basis of the tangent space to the graph at $F(X)$.
The symbol ${}^\perp$ denotes orthogonal projection onto the normal
space at $F(X)$. For $a\in\R^m$, write $n_a:=(0,a)^\perp$.
The normal vectors have the form $(-\varepsilon(df_X)^T\eta,\eta)$.
Requiring $(0,a)-n_a$ to be tangent gives
$(I_m+\varepsilon^2df_X(df_X)^T)\eta=a$, so
\[
(0,a)^\perp = n_a=(-\varepsilon(df_X)^TB_\varepsilon a,B_\varepsilon a).
\]
For $b\in\R^m$, the vector $(0,b)-n_b$ is tangent. Hence
\begin{equation}\label{eq:normal-projection-inner-product}
 \langle n_a,n_b\rangle=\langle n_a,(0,b)\rangle
 =\langle B_\varepsilon a,b\rangle.
\end{equation}
The second fundamental form has coordinate components
$\mathrm{II}_{ik}:=(\partial_i\partial_kF)^\perp$.
We apply the Gauss equation \eqref{eq:gauss-equation} from
Section~\ref{sec:curvature-conventions} with the Euclidean ambient
metric.
Since $\partial_i\partial_kF=(0,\varepsilon f_{ik})$, where
$f_{ik}:=\partial_i\partial_k f$, we have
$\mathrm{II}_{ik}=n_{\varepsilon f_{ik}}$.
Substituting this into \eqref{eq:gauss-equation} and using
\eqref{eq:normal-projection-inner-product} gives \eqref{eq:gauss-exact}.
For the family $f_\rho$, the matrices
$B_{\varepsilon_\rho}:=(I_m+\varepsilon_\rho^2df_\rho(df_\rho)^T)^{-1}$
converge to $I_m$ in $C^1(A)$:
the  matrices being inverted converge to $I_m$ in $C^1(A)$ by the assumed $C^2(A)$ convergence
and the bound $\|\varphi\|_{C^2(A)}<\infty$.
Applying \eqref{eq:gauss-exact} gives the convergence in
\eqref{eq:gauss-limit} in $C^0(A)$. If $\|\varphi\|_{C^3(A)}<\infty$
and $f_\rho\to\varphi$ in $C^3(A)$, the derivatives of $f_\rho$
through order three are uniformly bounded on $A$. Differentiating
\eqref{eq:gauss-exact} once therefore gives convergence in $C^1(A)$.
\end{proof}

The next result gives the obstruction used in
Theorem~\ref{thm:intro-branch}: if $\varepsilon_i\to0$ and
$f_i\to\operatorname{Re}(cz^\nu)$ in the norms specified
below, with $c\ne0$, then the induced metrics of
$\operatorname{graph}(\varepsilon_i f_i)$ are not locally conformally
flat for all sufficiently large $i$, except possibly when $k=3$ and
$\nu=2$.
\begin{theorem}[Local graph obstruction]\label{thm:local}
Let $k\geq3$, $\nu>1$, and $c\in\C^m\setminus\{0\}$.
Let $A$ be an annular sheet chart as in Definition~\ref{def:chart},
with the base plane identified with $\R^k=\R^2\times\R^{k-2}$.
Using the continuous angle $\theta$ on the conical sector $C_\Omega$
from that definition, set
\[
 z^\nu:=r^\nu e^{i\nu\theta},\qquad
 \varphi(z,t):=\operatorname{Re}(c z^\nu)\quad\text{on }A.
\]
Let $\varepsilon_i>0$ tend to zero and let $f_i:A\to\R^m$ be smooth functions with
\[
 f_i\longrightarrow\varphi
 \quad\begin{cases}
       \text{in }C^3(A),&k=3,\\
       \text{in }C^2(A),&k\geq4.
      \end{cases}
\]
If $k=3$ and $\nu\ne2$, or if $k\geq4$ and $\nu>1$ is arbitrary,
the induced metrics $G_i:=\delta+\varepsilon_i^2(df_i)^Tdf_i$ of
$\operatorname{graph}(\varepsilon_i f_i)$, expressed on $A$, are not
locally conformally flat for all sufficiently large $i$.

In particular, fix $\rho_0 > 0$ and let
$v:\bigcup_{0<\rho<\rho_0}\rho A\to\R^m$ be a smooth function
satisfying
\begin{equation}\label{eq:power-height-convergence}
 \rho^{-\nu}v(\rho\,\cdot)\longrightarrow\varphi
 \quad\begin{cases}
       \text{in }C^3(A),&k=3,\\
       \text{in }C^2(A),&k\geq4,
      \end{cases}
\end{equation}
as $\rho\downarrow0$. Under the same conditions on $k$ and $\nu$, the
induced metric $G:=\delta+(dv)^Tdv$ of $\operatorname{graph}(v)$,
expressed on its domain, is not locally conformally flat on the
dilation $\rho A$ of $A$ for all sufficiently small $\rho>0$.
This applies, in particular, when each restriction $v|_{\rho A}$ is
an exact sheet-height function as in Definition~\ref{def:chart}.
\end{theorem}
\begin{proof}
 The function $\varphi$ is smooth on a neighborhood of the compact set $\overline A$,
so its derivatives through order three are bounded on $A$.
Hence Lemma~\ref{lem:gauss} applies.
 We first compute:
\begin{align*}
\varphi_{xx} &= \operatorname{Re}(c\nu (\nu-1) z^{\nu - 2}),\\
\varphi_{xy} &= -\operatorname{Im}(c\nu (\nu-1) z^{\nu - 2}),\\
\varphi_{yy} &= -\varphi_{xx},\\
H:&=|\varphi_{xx}|^2+|\varphi_{xy}|^2
   =\nu^2(\nu-1)^2|c|^2r^{2\nu-4}.
\end{align*}
Hence, \eqref{eq:gauss-limit} gives
\[
 \varepsilon_i^{-2}\operatorname{Rm}_{1212}(G_i)\longrightarrow H.
\]
Because $\varphi$ is independent of $t$, the limit of $\varepsilon_i^{-2}\operatorname{Rm}_{abcd}(G_i)$ is zero
whenever any of $a,b,c,d$ belongs to $\{3,\ldots,k\}$.
Since $G_i^{-1}\to\delta^{-1}$, contraction of indices gives
\[
 \varepsilon_i^{-2}\operatorname{Ric}_{ab}(G_i)\longrightarrow
 \begin{cases}-H,&a=b\in\{1,2\},\\0,&\text{otherwise},\end{cases}
 \qquad
 \varepsilon_i^{-2}\operatorname{Scal}_{G_i}\longrightarrow-2H.
\]

For $k=3$, these limits give
\begin{equation}\label{eq:schouten-limit-k3}
 \varepsilon_i^{-2}(S_{G_i})_{tt}\longrightarrow H/2,\qquad
 \varepsilon_i^{-2}(S_{G_i})_{tx}\longrightarrow0
 \quad\text{in }C^1(A).
\end{equation}
Since $G_i\to\delta$ in $C^1(A)$, the Christoffel symbols
$\Gamma^p_{ab}(G_i)$ tend to zero in $C^0(A)$.
The covariant derivative of the Schouten tensor is
\[
 (\nabla_c S_{G_i})_{ab}
 =\partial_c(S_{G_i})_{ab}
 -\Gamma^p_{ca}(G_i)(S_{G_i})_{pb}
 -\Gamma^p_{cb}(G_i)(S_{G_i})_{ap}.
\]
The tensors $\varepsilon_i^{-2}S_{G_i}$ are bounded in $C^0(A)$
by the preceding curvature limits. Hence, after multiplication by
$\varepsilon_i^{-2}$, terms containing Christoffel symbols
tend to zero in $C^0(A)$.
Using the $C^1(A)$ convergence in \eqref{eq:schouten-limit-k3}
and the definition of the Cotton tensor, we obtain
\begin{align}
 \varepsilon_i^{-2}C_{ttx}(G_i)
 &\longrightarrow
\frac{1}{2} \partial_x H =  \nu^2(\nu-1)^2(\nu-2)|c|^2r^{2\nu-6}x,\label{eq:cotton-new}\\
 \varepsilon_i^{-2}C_{tty}(G_i)
 &\longrightarrow
\frac{1}{2} \partial_y H =  \nu^2(\nu-1)^2(\nu-2)|c|^2r^{2\nu-6}y.\nonumber
\end{align}
For $\nu\ne2$ these two components cannot both vanish at any $r>0$.

For $k\geq4$ and any fixed $\ell\in\{3,\ldots,k\}$, the definition of the
Weyl tensor gives
\begin{equation}\label{eq:weyl-new}
 \varepsilon_i^{-2}W_{1\ell\ell1}(G_i)
 \longrightarrow
 \frac{k-3}{(k-2)(k-1)} H = \frac{k-3}{(k-2)(k-1)}
       \nu^2(\nu-1)^2|c|^2r^{2\nu-4},
\end{equation}
which is nonzero. If $G_i$ were locally conformally flat on $A$, its
Cotton tensor would vanish identically for $k=3$, and its Weyl tensor
would vanish identically for $k\geq4$. However, by the limits above, this is impossible for $i$ large enough.

For the last claim take
$f_i(X):=\rho_i^{-\nu}v(\rho_iX)$ and
$\varepsilon_i:=\rho_i^{\nu-1}$ along any sequence $\rho_i\downarrow0$.
Writing $d_\rho(X):=\rho X$, the induced metric $G_\rho$ of
$\operatorname{graph}(\rho^{-1}v(\rho\,\cdot))$, expressed on $A$, satisfies
\[
 G_\rho=\delta+\rho^{2\nu-2}
 \bigl(d(\rho^{-\nu}v(\rho\,\cdot))\bigr)^T
 d(\rho^{-\nu}v(\rho\,\cdot))
 =\rho^{-2}d_\rho^*G.
\]
Pullback by the diffeomorphism $d_\rho:A\to\rho A$ and multiplication
by a positive constant preserve local conformal flatness. Thus
$G$ is locally conformally flat on $\rho A$ if and only if $G_\rho$
is locally conformally flat on $A$.
If $G$ were locally conformally flat on $\rho A$ for arbitrarily small
$\rho>0$, such a sequence would contradict the first part of the theorem.
\end{proof}

\section{A curved-ambient extension and the quadratic case}\label{sec:extensions}
We first extend the conclusion of Theorem~\ref{thm:local} to graphs in
a smooth ambient metric, under the same convergence assumption and
with $1<\nu<2$. We then discuss the case $\nu=2$.

\subsection{A local graph obstruction in a smooth ambient metric}
We now consider the induced metric of $\operatorname{graph}(v)$ when
the ambient metric is a smooth Riemannian metric $h$.
We assume that $\rho^{-\nu}v(\rho\,\cdot)$ converges on a fixed annular
sheet chart as in \eqref{eq:power-height-convergence}.
For $1<\nu<2$, the contribution of the ambient Riemann curvature
tensor vanishes after the normalization used in
Theorem~\ref{thm:local}, so the same obstruction applies.
\begin{proposition}[Local graph obstruction in a smooth ambient metric]
\label{prop:curved}
Let $k\geq3$, and let $h$ be a smooth Riemannian metric near
$0\in\R^{k+m}$. Choose normal coordinates centered at $0$, so that
\[
 h_{\alpha\beta}(0)=\delta_{\alpha\beta},\qquad
 \partial_\gamma h_{\alpha\beta}(0)=0
 \quad\text{for all }\alpha,\beta,\gamma\in\{1,\ldots,k+m\}.
\]
In these coordinates, use the splitting
$\R^{k+m}=\R^k\times\R^m$ and let $A\subset\R^k$ be an annular
sheet chart as in Definition~\ref{def:chart}.
Fix $\rho_0>0$, and let
$v:\bigcup_{0<\rho<\rho_0}\rho A\to\R^m$ be a smooth function
whose graph is contained in the coordinate neighborhood.
Suppose that \eqref{eq:power-height-convergence} holds with
$\varphi(z,t)=\operatorname{Re}(cz^\nu)$, where
$c\in\C^m\setminus\{0\}$ and $1<\nu<2$, using the angle from
Definition~\ref{def:chart}.
For the parametrization $\Gamma(X):=(X,v(X))$, the induced metric
$g:=\Gamma^*h$ of $\operatorname{graph}(v)$, expressed on its domain,
is not locally conformally flat on the dilation $\rho A$ of $A$
for all sufficiently small $\rho>0$.
\end{proposition}
\begin{proof}
Let $U\subset\R^{k+m}$ be the normal coordinate neighborhood of $0$
for $h$. For $D_\rho(Y):=\rho Y$, where $Y \in \mathbb{R}^{k+m}$, set
$h_\rho:=\rho^{-2}D_\rho^*h$, which is defined on
$\rho^{-1}U:=\{Y\in\R^{k+m}:\rho Y\in U\}$.
Its coefficients are
\[
 (h_\rho)_{\alpha\beta}(Y)=h_{\alpha\beta}(\rho Y).
\]
For every fixed compact set $K\subset\R^{k+m}$, we have
$\rho K\subset U$ for all sufficiently small $\rho>0$.
Taylor's theorem for $h_{\alpha \beta}(\rho \cdot)$ and the normal-coordinate identities then give that for every fixed integer $j\geq0$,
\begin{equation}\label{eq:ambient-metric-estimate}
 \|h_\rho-\delta_{\mathrm{amb}}\|_{C^j(K)} = O(\rho^2).
\end{equation}
Put $f_\rho(X):=\rho^{-\nu}v(\rho X)$ and
$\varepsilon_\rho:=\rho^{\nu-1}$. The map
$\Gamma_\rho(X):=(X,\varepsilon_\rho f_\rho(X))$ parametrizes
$\operatorname{graph}(\rho^{-1}v(\rho\,\cdot))$.
Since $A$ is bounded, $f_\rho\to\varphi$ in $C^0(A)$, and
$\varepsilon_\rho\to0$, the sets $\Gamma_\rho(A)$ lie in a single
fixed compact ball $K\subset\R^{k+m}$ for all sufficiently small
$\rho>0$. For these $\rho$, after decreasing them further if necessary,
$\rho K\subset U$, so the preceding estimates hold on $\Gamma_\rho(A).$
Its induced metrics for the ambient metrics $h_\rho$ and
$\delta_{\mathrm{amb}}$, respectively, are
\[
 g_\rho:=\Gamma_\rho^*h_\rho,\qquad
 G_\rho:=\Gamma_\rho^*\delta_{\mathrm{amb}}
 =\delta+\varepsilon_\rho^2(df_\rho)^Tdf_\rho.
\]
We apply the Gauss equation \eqref{eq:gauss-equation} from
Section~\ref{sec:curvature-conventions} to $\Gamma_\rho$ with ambient
metric $h_\rho$. Its first term is the ambient Riemann curvature
tensor evaluated on the tangent vectors
\[
 E_a^\rho:=d\Gamma_\rho(\partial_a)
 =(e_a,\varepsilon_\rho\partial_a f_\rho),
\]
where $e_a$ is the $a$th standard basis vector of $\R^k$.
By \eqref{eq:ambient-metric-estimate} with $j=3$, we have
$\operatorname{Rm}(h_\rho)=O(\rho^2)$ in $C^1(K)$.
The uniform $C^2(A)$ bounds on $f_\rho$ give uniform $C^1(A)$
bounds on $E_a^\rho$, since
\[
 \partial_bE_a^\rho
 =(0,\varepsilon_\rho\partial_b\partial_a f_\rho).
\]
Since $\Gamma_\rho(A)\subset K$ and
$d\Gamma_\rho(\partial_a)=E_a^\rho$, we obtain
\[
 h_\rho\bigl(R^{h_\rho}(E_a^\rho,E_b^\rho)E_c^\rho,E_d^\rho\bigr)
 =O(\rho^2)\quad\text{in }C^1(A).
\]
Next, we justify why the coordinate components of the second fundamental form for
$h_\rho$ differ from those for $\delta_{\mathrm{amb}}$ by
$O(\rho^2)$ in $C^0(A)$ when
$f_\rho\to\varphi$ in $C^2(A)$, and in $C^1(A)$ when
$f_\rho\to\varphi$ in $C^3(A)$.
Let $\Pi_\rho$ and $\Pi_0$ denote the orthogonal projections onto
the normal space of $\Gamma_\rho(A)$ for $h_\rho$ and
$\delta_{\mathrm{amb}}$, respectively. In ambient coordinates
$(Y^1,\ldots,Y^{k+m})$, define the vector $B^\rho_{ab}$ by
\[
 (B^\rho_{ab})^\alpha
 :=\sum_{\beta,\gamma=1}^{k+m}
 \Gamma^\alpha_{\beta\gamma}(h_\rho)(\Gamma_\rho(X))
 (E_a^\rho)^\beta(E_b^\rho)^\gamma,
\]
where $\Gamma^\alpha_{\beta\gamma}(h_\rho)$ are the Christoffel
symbols of $h_\rho$. The second fundamental forms are
\[
 \mathrm{II}^{h_\rho}_{ab}
 =\Pi_\rho\bigl(\partial_a\partial_b\Gamma_\rho+B^\rho_{ab}\bigr),
 \qquad
 \mathrm{II}^{\delta_{\mathrm{amb}}}_{ab}
 =\Pi_0\bigl(\partial_a\partial_b\Gamma_\rho\bigr).
\]
Consequently,
\[
 \mathrm{II}^{h_\rho}_{ab}-\mathrm{II}^{\delta_{\mathrm{amb}}}_{ab}
 =(\Pi_\rho-\Pi_0)\partial_a\partial_b\Gamma_\rho
   +\Pi_\rho B^\rho_{ab}.
\]
  For explicit projection formulas, let $E_\rho$ be the matrix
with columns $E_a^\rho$ and put
$H_\rho:=h_\rho(\Gamma_\rho(X))$. Then
\[
 \begin{aligned}
 \Pi_\rho&=I-E_\rho(E_\rho^TH_\rho E_\rho)^{-1}E_\rho^TH_\rho,\\
 \Pi_0&=I-E_\rho(E_\rho^TE_\rho)^{-1}E_\rho^T.
 \end{aligned}
\]
The matrices $E_\rho^T H_\rho E_\rho$ and $E_\rho^T E_\rho$ are uniformly positive definite.
Indeed, the first $k$
rows of $E_\rho$ form the identity and, by
\eqref{eq:ambient-metric-estimate}, $H_\rho$ is uniformly close
to the identity. The $C^2$ bounds on $f_\rho$ give $C^1$ bounds
on $E_\rho$. Differentiating these formulas once and using
$H_\rho-I=O(\rho^2)$  in $C^1(A)$  gives
$\Pi_\rho-\Pi_0=O(\rho^2)$ in $C^1(A)$. The second coordinate derivatives of $\Gamma_\rho$ are uniformly
bounded in $C^0(A)$ under the $C^2(A)$ convergence of $f_\rho$,
and in $C^1(A)$ under its $C^3(A)$ convergence. Also, \eqref{eq:ambient-metric-estimate} with $j=2$ gives
Christoffel symbols of $h_\rho$ that are $O(\rho^2)$ in $C^1(K)$. Combining these proves the
stated estimates for the difference.
Thus \eqref{eq:gauss-equation} gives
\[
 \operatorname{Rm}(g_\rho)-\operatorname{Rm}(G_\rho)=O(\rho^2)
\]
in $C^0(A)$ for $k\geq4$, and in $C^1(A)$ for $k=3$. Since
\[
 \frac{\rho^2}{\varepsilon_\rho^2}
 =\rho^{4-2\nu}\longrightarrow0
 \qquad(1<\nu<2),
\]
division by $\varepsilon_\rho^2$ gives the same Riemann curvature
tensor limits as in the Euclidean case.
Moreover,
\[
 \begin{aligned}
 g_\rho-\delta
 &=\varepsilon_\rho^2(df_\rho)^Tdf_\rho
   +\Gamma_\rho^*(h_\rho-\delta_{\mathrm{amb}})\\
 &=O(\varepsilon_\rho^2+\rho^2)\quad\text{in }C^1(A).
 \end{aligned}
\]
Here the first term is controlled by the uniform $C^2(A)$ bounds
on $f_\rho$, and the second by
\eqref{eq:ambient-metric-estimate} and the uniform $C^2(A)$
bounds on $\Gamma_\rho$. Thus $g_\rho\to\delta$ in $C^1(A)$.
Therefore,
the contractions and covariant derivatives in the
proof of Theorem~\ref{thm:local} give the Cotton tensor limit
\eqref{eq:cotton-new} for $k=3$ and the Weyl tensor limit
\eqref{eq:weyl-new} for $k\geq4$.

These nonzero limits rule out local conformal flatness of $g_\rho$
on $A$ for all sufficiently small $\rho>0$.
Finally, with $d_\rho(X):=\rho X$ on the base $\R^k$, the identity
$D_\rho\circ\Gamma_\rho=\Gamma\circ d_\rho$ gives
\[
 g_\rho=\rho^{-2}d_\rho^*g.
\]
Pullback by $d_\rho:A\to\rho A$ and multiplication by a positive
constant preserve local conformal flatness. Hence, the conclusion follows.
\end{proof}

\begin{remark}
To extend Theorem~\ref{thm:intro-branch} to locally area-minimizing
integral currents in a smooth ambient metric, one would also need
analogues of the approximation by normalized tangent functions
and the estimates on the support of the current used in
Section~\ref{sec:current}, which we do not have. At $\nu=2$, the contribution of the ambient Riemann curvature tensor
can have the same order as the terms involving the second fundamental
form in the Gauss equation, so division by $\varepsilon_\rho^2$
does not force that contribution to vanish.
\end{remark}

\subsection{Spherical harmonics and the quadratic exception}
\label{sec:spherical-definitions}
We recall the definition of spherical harmonics and their Laplacian
eigenvalues, which are also used in the proof of
Theorem~\ref{thm:homogeneous-polynomial}. We then look at an example of a minimal
graph whose induced metric is locally conformally flat and whose
height function satisfies
$\rho^{-2}v(\rho\,\cdot)\to\operatorname{Re}(cz^2)$ with $c\ne0$.
This shows why Theorem~\ref{thm:local} excludes $\nu=2$ when $k=3$.

For an integer $d_{\rm sph}\geq1$, define the \defterm{unit sphere} by
\[
 S^{d_{\rm sph}}:=\{\theta\in\R^{d_{\rm sph}+1}:|\theta|=1\}.
\]
Its \defterm{unit round metric} $g_{S^{d_{\rm sph}}}$ is the induced
metric from $\R^{d_{\rm sph}+1}$.
We write $\Delta_{S^{d_{\rm sph}}}$ for the Laplacian of this metric.

For an integer $h\geq0$, a
\defterm{spherical harmonic of degree $h$ on $S^{d_{\rm sph}}$}
is the restriction of a real homogeneous polynomial $H_h$ of degree
$h$ satisfying $\Delta_{\R^{d_{\rm sph}+1}}H_h=0$.
Write $\mathcal H_h(S^{d_{\rm sph}})$ for the space of these
restrictions, as in \cite[\S1.4]{Gallier}.
Every $Y\in\mathcal H_h(S^{d_{\rm sph}})$ satisfies
\begin{equation}\label{eq:spherical-harmonic-eigenvalue}
 -\Delta_{S^{d_{\rm sph}}}Y
 =h(h+d_{\rm sph}-1)Y.
\end{equation}
Moreover, these are all the eigenvalues of
$-\Delta_{S^{d_{\rm sph}}}$, and $\mathcal H_h(S^{d_{\rm sph}})$
is the eigenspace for $h(h+d_{\rm sph}-1)$
\cite[Proposition~1.7 and Theorem~1.9]{Gallier}.

\paragraph{The Bor\r{u}vka sphere and its cone.}
On $S^2$, write $dS$ for the Riemannian volume measure of $g_{S^2}$.
The space $\mathcal H_3(S^2)$ has dimension $7$, since
$\dim\mathcal H_h(S^2)=2h+1$ \cite[Corollary~1.5]{Gallier}.
Choose a basis $Y_1,\ldots,Y_7$ orthonormal for
$\langle f,g\rangle_{L^2}:=\int_{S^2}fg\,dS$, and define
\[
 \beta:S^2\longrightarrow S^6\subset\R^7,\qquad
 \beta(\theta):=\sqrt{\frac{4\pi}{7}}\,
       (Y_1(\theta),\ldots,Y_7(\theta)).
\]
This is the \defterm{Bor\r{u}vka sphere in $S^6$}. Its image is an
embedded two-sphere, and
\begin{equation}\label{eq:boruvka-sphere}
 |\beta|=1,\qquad \Delta_{S^2}\beta=-12\beta,\qquad
 \beta^*g_{S^6}=6g_{S^2}.
\end{equation}
More about the construction can be found in~\cite[\S2.1]{BM}. We verify the normalization and the induced metric.
Orthogonal transformations of $\R^3$ act by orthogonal changes
of the basis $Y_1,\ldots,Y_7$. Thus $\sum_iY_i^2$ is constant on $S^2$.
Since its integral is $7$ and $S^2$ has area $4\pi$, this constant
equals $7/(4\pi)$.
This proves $|\beta|=1$. Invariance under rotations of $S^2$ also
gives $\beta^*g_{S^6}=\kappa g_{S^2}$ for a constant $\kappa$, so $|d\beta|^2=2\kappa$.
Applying $\Delta_{S^2}$ to $|\beta|^2=1$ gives
\[
 0=2\langle\Delta_{S^2}\beta,\beta\rangle+2|d\beta|^2
   =-24+2|d\beta|^2 = -24 + 4 \kappa.
\]
Hence, $\kappa=6$.

Define the \defterm{Bor\r{u}vka cone}, including its vertex, by
\[
 \mathcal C_\beta:=\{R\beta(\theta):R\geq0,\ \theta\in S^2\}.
\]
We choose the orthonormal basis \(Y_1,\ldots,Y_7\) so that \(\mathcal C_\beta\setminus\{0\}\), with a suitable orientation, is calibrated by the standard associative \(3\)-form \(\varphi_0\) on $\mathbb R^7$, which is given by
\[
 \varphi_0:=dy^{123}+dy^{145}+dy^{167}+dy^{246}
                -dy^{257}-dy^{347}-dy^{356},
 \qquad dy^{ijk}:=dy^i\wedge dy^j\wedge dy^k
\]
in orthonormal coordinates $y^1,\ldots,y^7$ on $\R^7$.
This associativity property is recalled in
\cite[\S1, the discussion preceding Theorem~1.4]{BM}.
  We choose $\beta$ so that the parametrization $R\theta\mapsto R\beta(\theta)$
induces the calibrated orientation of the cone.
We fix this parametrization, these target coordinates,  and the
calibrated orientation from now on.

\begin{example}[The Bor\r{u}vka cone and the quadratic exception]
\label{ex:boruvka-quadratic}
We claim that the Bor\r{u}vka cone $\mathcal C_\beta$ has the following properties:
\begin{itemize}
\item The submanifold $\mathcal C_\beta\setminus\{0\}$ is minimal.
\item Its induced metric is locally conformally flat.
\item For any $\theta_0\in S^2$, set $Z:=\beta(\theta_0)$.
A neighborhood of $Z$ in this submanifold is
$Z+\operatorname{graph}(v)$ for a smooth function
\[
 v:B_\eta^{T_Z\mathcal C_\beta}(0)\longrightarrow
       (T_Z\mathcal C_\beta)^\perp
\]
and some $\eta>0$, with $v(0)=0$ and $dv(0)=0$.
Choose orthonormal coordinates $X=(x,y,t)$ on $T_Z\mathcal C_\beta$,
with $t$ in the radial direction, and orthonormal normal vectors
$n_1,n_2$ as in the proof below. Put $z:=x+iy$ and
$c:=\frac{\sqrt{15}}{12}(n_1-i n_2)$.
The height function $v$ satisfies
\[
 v(X)=\varphi(X)+O(|X|^3),\qquad
 \varphi(x,y,t):=\frac{\sqrt{15}}{12}
       \bigl((x^2-y^2)n_1+2xy\,n_2\bigr)
       =\operatorname{Re}(cz^2).
\]
In particular, $c\ne0$ and $c\cdot c=0$, where $\cdot$ is defined in \eqref{eq:complex-bilinear-product}.
\end{itemize}
\end{example}
\begin{proof}
Associativity implies minimality. By \eqref{eq:boruvka-sphere}, the
induced metric of $\mathcal C_\beta\setminus\{0\}$, expressed on
$(0,\infty)\times S^2$ by $(R,\theta)\mapsto R\beta(\theta)$, is
\[
 g_{\mathcal C_\beta}=\dd R^2+6R^2g_{S^2}.
\]
With $\varrho:=R^{1/\sqrt6}$, this is
$6R^2\varrho^{-2}(\dd\varrho^2+\varrho^2g_{S^2})$.
The metric in parentheses is Euclidean in polar coordinates $(\varrho,\theta)$,
so $g_{\mathcal C_\beta}$ is locally conformally flat for $R>0$.

Let $e_1,e_2$ be an orthonormal basis of
$T_Z\beta(S^2) \subset T_Z\mathcal C_\beta$, and let $e_3:=Z \in T_Z\mathcal C_\beta$ be the unit radial vector.

Write $\mathrm{II}_Z$ for the second fundamental form of
$\mathcal C_\beta\setminus\{0\}\subset\R^7$ at $Z$.
In local coordinates $(\theta^1,\theta^2)$ on $S^2$, use the
parametrization $F(R,\theta^1,\theta^2):=R\beta(\theta^1,\theta^2)$.
Its derivatives satisfy
\[
 \frac{\partial^2F}{\partial R^2}=0,\qquad
 \frac{\partial^2F}{\partial R\,\partial\theta^a}
 =\frac{1}{R}\frac{\partial F}{\partial\theta^a}
 \quad(a=1,2).
\]
The mixed derivatives are tangent to the cone because
$\partial F/\partial\theta^a$ is tangent to the cone. Thus $\mathrm{II}_Z(e_3,\cdot)=0$. Minimality gives
$\mathrm{II}_Z(e_2,e_2)=-\mathrm{II}_Z(e_1,e_1)$.
Put
\[
 U:=\mathrm{II}_Z(e_1,e_1),\qquad
 V:=\mathrm{II}_Z(e_1,e_2).
\]
We now use the rotational symmetry of the cone to prove that
$U$ and $V$ have equal length and are orthogonal.
Let $Q\in\mathrm{SO}(3)$ satisfy $Q\theta_0=\theta_0$.

The action on the spherical harmonics $Y_1,\ldots,Y_7$, which form
an orthonormal basis of $\mathcal H_3(S^2)$, determines
a unique orthogonal transformation $A_Q$ of $\R^7$ such that
\begin{equation}\label{eq:boruvka-equivariance}
 A_Q\beta(\theta)=\beta(Q\theta)
 \qquad(\theta\in S^2).
\end{equation}
In particular, $A_QZ=\beta(Q\theta_0)=Z$, and
$A_Q(R\beta(\theta))=R\beta(Q\theta)$, so $A_Q$ preserves
$\mathcal C_\beta$. Differentiating \eqref{eq:boruvka-equivariance} at $\theta_0$ gives
\begin{equation}\label{eq:boruvka-differentiated-equivariance}
 A_Q\circ d\beta_{\theta_0}
 =d\beta_{\theta_0}\circ Q|_{T_{\theta_0}S^2}.
\end{equation}
Since $\beta^*g_{S^6}=6g_{S^2}$ by \eqref{eq:boruvka-sphere},
$|d\beta_{\theta_0}(v)|=\sqrt6\,|v|$ for every
$v\in T_{\theta_0}S^2$. Thus
\[
 O:=\frac{1}{\sqrt6}\,d\beta_{\theta_0}:
 T_{\theta_0}S^2\longrightarrow T_Z\beta(S^2)
\]
is a linear isometry. Dividing
\eqref{eq:boruvka-differentiated-equivariance} by $\sqrt6$ and
composing with $O^{-1}$ gives
\[
 A_Q|_{T_Z\beta(S^2)}
 =O\circ Q|_{T_{\theta_0}S^2}\circ O^{-1}.
\]
Hence as $Q$ ranges over rotations about the axis $\R\theta_0$, its
restriction ranges over all rotations of $T_{\theta_0}S^2$.
The induced transformations therefore give all plane rotations of
$T_Z\beta(S^2)=\operatorname{span}(e_1,e_2)$.

For any angle $\alpha$, choose $Q$ so that its induced action satisfies
\[
 A_Qe_1=\cos\alpha\,e_1+\sin\alpha\,e_2,\qquad
 A_Qe_2=-\sin\alpha\,e_1+\cos\alpha\,e_2,\qquad
 A_Qe_3=e_3.
\]
Since $A_Q$ is an ambient isometry fixing $Z$ and preserving the cone,
the second fundamental form is equivariant under $A_Q$, meaning that
\[
 A_Q\bigl(\mathrm{II}_Z(v,w)\bigr)
 =\mathrm{II}_Z(A_Qv,A_Qw)
 \qquad(v,w\in T_Z\mathcal C_\beta).
\]
By bilinearity and
$\mathrm{II}_Z(e_2,e_2)=-U$, we obtain
\[
 \begin{aligned}
 A_QU
 &=\mathrm{II}_Z(A_Qe_1,A_Qe_1)\\
 &=\cos^2\alpha\,U
   +2\sin\alpha\cos\alpha\,V-\sin^2\alpha\,U\\
 &=\cos(2\alpha)U+\sin(2\alpha)V.
 \end{aligned}
\]
Since $A_Q$ preserves norms, choosing $\alpha=\pi/4$ gives
$|U|=|V|$. Then choosing $\alpha=\pi/8$ gives
$\langle U,V\rangle=0$.
Next, we compute the sectional curvature of
$\operatorname{span}(e_1,e_2)$ by applying the Gauss equation
\eqref{eq:gauss-equation} from Section~\ref{sec:curvature-conventions}
to the two inclusions
\[
 \Sigma:=\beta(S^2)
 \ \subset\ \mathcal C_\beta\setminus\{0\}
 \ \subset\ \R^7.
\]
For the first inclusion, the ambient manifold is the cone with metric
$g_{\mathcal C_\beta}$, and we write $B$ for the second fundamental
form of $\Sigma$ inside the cone. Its values are tangent to the cone
and normal to $\Sigma$ within the cone.

To compute $B$, consider the cross-sections
$\Sigma_R:=\{R\beta(\theta):\theta\in S^2\}$.
Their unit normal within the cone is
$N:=dF(\partial_R)=\beta(\theta)$, which at $Z$ equals $e_3$.
In the coordinates $(R,\theta^1,\theta^2)$, the cone metric has
coefficients
\[
 g_{RR}=1,\qquad g_{Ra}=0,\qquad
 g_{ab}=6R^2(g_{S^2})_{ab}\quad(a,b=1,2).
\]
The radial Christoffel symbols are therefore
\[
 \Gamma^R_{ab}
 =-\frac{1}{2}\partial_Rg_{ab}
 =-6R(g_{S^2})_{ab}.
\]
For $E_a:=dF(\partial_{\theta^a})$,
\[
 \begin{aligned}
 B^{\Sigma_R\subset\mathcal C_\beta}(E_a,E_b)
 &:=g_{\mathcal C_\beta}
       \bigl(\nabla^{\mathcal C_\beta}_{E_a}E_b,N\bigr)N\\
 &=\Gamma^R_{ab}N\\
 &=-\frac{1}{R} g_{\mathcal C_\beta}(E_a,E_b)N.
 \end{aligned}
\]
At $R=1$, $\Sigma_R=\Sigma$, so in our orthonormal basis
\[
 B(e_1,e_1)=B(e_2,e_2)=-e_3,\qquad B(e_1,e_2)=0.
\]

Let $g_\Sigma$ be the induced metric on $\Sigma$.
By \eqref{eq:boruvka-sphere}, the parametrization
$\beta:S^2\to\Sigma$ satisfies $\beta^*g_\Sigma=6g_{S^2}$,
so the sectional curvature of $g_\Sigma$ is $1/6$.
Hence, applying~\eqref{eq:gauss-equation} to $\Sigma\subset\mathcal C_\beta$ and recalling that $K(e_1\wedge e_2)=\operatorname{Rm}(e_1,e_2,e_2,e_1)$, gives
\[
 \begin{aligned}
 \frac{1}{6}
 &=K_{\mathcal C_\beta}(e_1\wedge e_2)
   +\langle B(e_1,e_1),B(e_2,e_2)\rangle
   -|B(e_1,e_2)|^2\\
 &=K_{\mathcal C_\beta}(e_1\wedge e_2)+1.
 \end{aligned}
\]
Thus $K_{\mathcal C_\beta}(e_1\wedge e_2)=-5/6$.

We now apply \eqref{eq:gauss-equation} to the second inclusion,
$\mathcal C_\beta\setminus\{0\}\subset\R^7$.
Its second fundamental form is $\mathrm{II}_Z$, with $U$ and $V$ as above.
The ambient curvature term vanishes because $\R^7$ is Euclidean.
Using the same indices $(1,2,2,1)$, we obtain
\[
 \begin{aligned}
 -\frac{5}{6}
 &=\langle\mathrm{II}_Z(e_1,e_1),\mathrm{II}_Z(e_2,e_2)\rangle
       -|\mathrm{II}_Z(e_1,e_2)|^2\\
 &=\langle U,-U\rangle-|V|^2
 =-|U|^2-|V|^2.
 \end{aligned}
\]
Since $|U|=|V|$, this gives $2|U|^2=5/6$, and hence
$|U|=|V|=\sqrt{5/12}$. Define
$n_1:=U/|U|$ and $n_2:=V/|V|$.

Let $\pi_Z:\R^7\to T_Z\mathcal C_\beta$ be the orthogonal projection.
The map $q\mapsto\pi_Z(q-Z)$, restricted to
$\mathcal C_\beta\setminus\{0\}$ near $Z$, has differential
$\pi_Z|_{T_Z\mathcal C_\beta}=\operatorname{id}_{T_Z\mathcal C_\beta}$
at $Z$. It is therefore invertible nearby. The inverse function theorem gives
a parametrization $G(X):=Z+X+v(X)$ of an open neighborhood of $Z$
in $\mathcal C_\beta\setminus\{0\}$, with
$X\in B_\eta^{T_Z\mathcal C_\beta}(0)$ for some $\eta>0$ and $v(0)=0$.
For every $\xi\in T_Z\mathcal C_\beta$, the vector
$dG_0(\xi)=\xi+dv(0)\xi$ is tangent to the cone at $Z$.
Since $\xi$ is also tangent but $dv(0)\xi$ is normal, we obtain $dv(0)=0$.

Write $X=X^1e_1+X^2e_2+X^3e_3$, where
$(X^1,X^2,X^3)=(x,y,t)$, and put $\partial_i:=\partial/\partial X^i$.
Since $G(X)=Z+X+v(X)$ and $dv(0)=0$,
\[
 \partial_iG(0)=e_i,\qquad
 \partial_i\partial_jG(0)=\partial_i\partial_jv(0).
\]
Apply the definition of the second fundamental form from
Section~\ref{sec:curvature-conventions} to $G$ and the coordinate
frame $\partial_1,\partial_2,\partial_3$ to get
\[
 \begin{aligned}
 \mathrm{II}_Z(e_i,e_j)
 &=\mathrm{II}_Z\bigl(\partial_iG(0),\partial_jG(0)\bigr)\\
 &=\bigl(\partial_i\partial_jG(0)\bigr)^\perp\\
 &=\partial_i\partial_jv(0).
 \end{aligned}
\]
Thus  $D^2v(0)$ from Section~\ref{sec:derivative-notation}
consists of the components $\mathrm{II}_Z(e_i,e_j)$ in these coordinates.
Taylor's theorem gives
\[
 \begin{aligned}
 v(X)
 &=\frac{1}{2}\sum_{i,j=1}^3X^iX^j\partial_i\partial_jv(0)+O(|X|^3)\\
 &=\frac{1}{2}\mathrm{II}_Z(X,X)+O(|X|^3)\\
 &=\frac{1}{2}\bigl((x^2-y^2)U+2xyV\bigr)+O(|X|^3),
 \end{aligned}
\]
which is the expression we claimed for $\varphi$.
Fix an annular sheet chart $A$ in $T_Z\mathcal C_\beta$ with
$\{z=0\}=\operatorname{span}(e_3)$.
Smoothness of $v$ and Taylor's theorem give
$\rho^{-2}v(\rho\,\cdot)\to\varphi$ in $C^3(A)$ as $\rho\downarrow0$,
so \eqref{eq:power-height-convergence} holds with $\nu=2$.
\end{proof}

\begin{remark}
Example~\ref{ex:boruvka-quadratic} concerns a smooth point of
$\mathcal C_\beta\setminus\{0\}$, not a branch point of an integral
current. It shows that the conclusion of Theorem~\ref{thm:local}
fails when $k=3$ and $\nu=2$, even when the graph is associative
and $c\cdot c=0$.
\end{remark}

\begin{remark}[The cone over the Veronese surface]
The cone over the Veronese   surface gives the same quadratic exception:
near each nonzero point it is a minimal graph with locally conformally
flat induced metric and a nonzero quadratic leading term of the form
$\operatorname{Re}(cz^2)$. See \cite[\S2.1]{BM} for the construction
of the Veronese   surface.
\end{remark}

\section{From tangent functions to the branch theorem}\label{sec:current}
We prove Theorem~\ref{thm:intro-branch} by applying the local graph
obstruction in Theorem~\ref{thm:local} to magnified currents.
The main step is to construct exact sheet-height functions when a
nonzero cyclic component has repetition multiplicity one.
In the cases covered by Theorem~\ref{thm:local}, this gives the
pointwise obstruction. Stratification then gives the Hausdorff
dimension and measure conclusions.

Throughout this section, $T$ satisfies the hypotheses of
Theorem~\ref{thm:intro-branch}.

Let $P$ be an oriented $k$-plane and $B\subset P$ a connected open set.
For a $k$-dimensional integral current $S$ in $B\times P^\perp$ with
bounded normal support, its \defterm{projection multiplicity} over $B$
is the integer $d$, whenever
\[
 (\pi_P)_\#S=d\llbracket B\rrbracket,
 \qquad \pi_P(X,Y):=X.
\]
Here $\llbracket B\rrbracket$ denotes
integration over $B$ with multiplicity one and the orientation of $P$. The graph current of a single-valued Lipschitz function has projection
multiplicity one. Note that, conversely, projection multiplicity one alone
does not imply that the current is the graph current of a smooth function.

\subsection{\texorpdfstring{From convergence of $q$-valued functions to currents with disjoint supports}{From convergence of q-valued functions to currents with disjoint supports}}
We first write the magnified current as a sum of currents with pairwise
disjoint supports, one near each distinct local branch of a cylindrical
normalized tangent function. Lemma~\ref{lem:matching} shows that each
of these currents has projection multiplicity equal to the repetition
multiplicity of the corresponding local branch.

Fix $Z\in\mathcal R_q(T)$, as defined in Definition~\ref{def:rapid},
and a cylindrical normalized tangent function $\widehat\Phi$.
Use the retained current from Lemma~\ref{lem:localization}, translate
$Z$ to $0$, and write $P:=P_Z$ and $\nu:=N_{T,\mathrm{Pl}}(Z)$.
Choose a sequence of radii $\rho_i\downarrow0$ realizing
$\widehat\Phi$ in Definition~\ref{def:tangent-function}, and use
$T_i=T_{Z,\rho_i}$ and $\varepsilon_i=e_Z(\rho_i)$ from
\eqref{eq:tangent-scales}. In particular, $\varepsilon_i>0$ and
$\varepsilon_i\to0$. By \eqref{eq:localized-cylinder-ball}, for each fixed $L>0$ and
all sufficiently large $i$,
\[
 \supp T_i\cap\bigl(B_L^P(0)\times P^\perp\bigr)
 \subset B_{8L/7}^{\R^{k+m}}(0).
\]
Together with the local mass convergence in
\eqref{eq:localization-tangent-convergence}, this gives, for each fixed
$L>0$ and all sufficiently large $i$,
\begin{equation}\label{eq:magnified-local-bounds}
 \begin{aligned}
 \|T_i\|\bigl(B_L^P(0)\times P^\perp\bigr)&\leq C_L,\\
 |Y|&\leq\frac{8L}{7}
 \quad\bigl((X,Y)\in\supp T_i,\ |X|<L\bigr),
 \end{aligned}
\end{equation}
where $C_L<\infty$ is independent of $i$.
By \eqref{eq:localized-projection-identity}, after magnification,
$(\pi_P)_\#T_i=q\llbracket B\rrbracket$ on each fixed base ball $B$
for all sufficiently large $i$.

Use the branch axis $\mathcal L_Z$, cyclic components, and local
branches from Section~\ref{sec:cyclic-components}. Choose a nonempty
connected conical sector $C_\Omega$ with its continuous angle as in
Definition~\ref{def:chart}.
List the distinct local branches as $\psi_1,\ldots,\psi_s$, taking
values in $P^\perp$, including the zero function if a zero component
is present. Identical repetitions are listed only once.
An annular sheet chart $A\Subset C_\Omega$ is
\defterm{profile-separated} if there exists $\delta>0$ such that
\[
 |\psi_a(X)-\psi_b(X)|\geq\delta
 \quad\text{for every }X\in\overline A
 \text{ and }1\leq a<b\leq s.
\]
If $s=1$, the condition is automatic.
Such annular sheet charts can be chosen to be open balls away from
the branch axis and the
\defterm{coincidence sets} $E_{ab}:=\{X:\psi_a(X)=\psi_b(X)\}$,
where $X$ ranges over the common domain of the branches and $a<b$.
The branches are real analytic on $C_\Omega$ by
formula~\eqref{eq:branches}. Since $C_\Omega$ is connected, the identity
principle implies that two such branches agreeing on a nonempty open
subset agree throughout $C_\Omega$.
Hence distinct branches have coincidence sets with empty interior.
Choose nested open balls $A\Subset A'\Subset A^+$ that are all
annular sheet charts in $C_\Omega$, with
$\overline{A^+}\subset\{X:1/2<|X|<2\}$ and $A^+$ profile-separated.
These nested annular sheet charts leave room for interior estimates.
Estimates on $A^+$ allow us to write the current over $A'$ as a sum of
currents with pairwise disjoint supports, one near each distinct local branch.
When one of these currents has projection multiplicity one, we obtain an exact
sheet-height function and derivative estimates on smaller annular
sheet charts, ultimately on $A$. Each smaller chart has positive distance
from the boundary of the next. All three charts are fixed as $i\to\infty$.
Write $n_1,\ldots,n_s$ for the repetition multiplicities of the branches,
so that
\[
 \widehat\Phi|_{A^+}=\sum_{a=1}^s n_a\llbracket\psi_a\rrbracket,
 \qquad \sum_a n_a=q.
\]
If $\psi_a$ is a local branch of the $j$th nonzero cyclic component,
then $n_a=m_j$, with $m_j$ as in Section~\ref{sec:cyclic-components}.
 This is because formula~\eqref{eq:branches} shows that the branches
within each cyclic component are distinct, and that branches from
distinct components cannot agree on an open set.
 For the zero function, if present, $n_a=m_0$.

Lemma~\ref{lem:matching} isolates the part of the current near each
distinct local branch and identifies its projection multiplicity
with that local branch's repetition multiplicity. It is a consequence of
Krummel--Wickramasekera's support estimate and standard multiplicity
arguments. We include the proof to verify their application in our setting.
\begin{lemma}[Support estimate and multiplicity matching]\label{lem:matching}
For all sufficiently large $i$, there are locally area-minimizing
integral currents $S_{i,a}$ with pairwise disjoint supports and no
boundary in $A'\times P^\perp$ such that
\[
 T_i\llcorner(A'\times P^\perp)=\sum_{a=1}^s S_{i,a}.
\]
Each $S_{i,a}$ satisfies
\begin{equation}\label{eq:multiplicity-matching}
 \sup_{(X,Y)\in\supp S_{i,a},\,X\in A'}
       |Y-\varepsilon_i\psi_a(X)|=o(\varepsilon_i),
 \qquad
 (\pi_P)_\#S_{i,a}=n_a\llbracket A'\rrbracket.
\end{equation}
Each $S_{i,a}$ has bounded normal support and is interpreted as a
current in the cylinder $A'\times P^\perp$.
\end{lemma}
\begin{proof}
We first check the hypotheses of \cite[Theorem~4.19]{KWIII} and
apply it. We then estimate the fine excess in its conclusion and use
the resulting bound to construct the currents $S_{i,a}$ and determine
their projection multiplicities.

\emph{Step 1: verify the hypotheses of \cite[Theorem~4.19]{KWIII}.}
Let $r_+>0$ and $X_+\in P$ denote the radius and center of $A^+$,
respectively, so that $A^+=B_{r_+}^P(X_+)$.
Choose $\gamma\in(0,1)$ such that
$\overline{A'}\subset B_{\gamma r_+}^P(X_+)$.
Then for all sufficiently large $i$, the magnified current $T_i$ is locally
area minimizing and has no boundary in $A^+\times P^\perp$.
Equation~\eqref{eq:magnified-local-bounds}, with $L$ chosen so that
$A^+\subset B_L^P(0)$, gives bounded normal support there,
and \eqref{eq:localized-projection-identity}, after magnification
and restriction to $A^+$, gives
\[
 (\pi_P)_\#T_i=q\llbracket A^+\rrbracket.
\]

Let $u_i$ and $K_i$ be the Lipschitz $q$-valued function and relatively
closed set from the construction preceding
Definition~\ref{def:tangent-function}, on a fixed ball containing
$\overline{A^+}$. For a fixed open ball $B\subseteq A^+$, define the
nonnegative \defterm{oriented tilt excess} $t_i(B)$ by
\begin{equation}\label{eq:oriented-tilt-excess}
 t_i(B)^2:=
 \frac{1}{2|B|}\int_{B\times P^\perp}
       |\vec T_i-\vec P|^2\,\dd\|T_i\|,
\end{equation}
where $|B|$ is the $k$-dimensional Lebesgue measure of $B$.
$\vec T_i$ is the oriented unit $k$-vector representing $T_i$ as in
Definition~\ref{def:integral-current}, and $\vec P$ represents the
oriented plane $P$. Since the projection multiplicity is $q$,
this is the excess in \cite[(2.2), p.~21]{KWIII}.
Taking $B=A^+$ in \eqref{eq:oriented-tilt-excess}, we first estimate
the part of its integral over $(K_i\cap A^+)\times P^\perp$.
Choose a positively oriented orthonormal basis $e_1,\ldots,e_k$ of
$P$, so that $\vec P=e_1\wedge\cdots\wedge e_k$.
For a linear map $M:P\to P^\perp$, the vectors
$v_a:=(e_a,Me_a)$ span its graph in $P\times P^\perp$.
Their matrix of inner products is $I+M^TM$, as in
\eqref{eq:induced-graph-metric}. Hence
\begin{equation}\label{eq:tilt-graph-jacobian}
 J_M:=|v_1\wedge\cdots\wedge v_k|
 =\sqrt{\det(\langle v_a,v_b\rangle)_{a,b=1}^k}
 =\sqrt{\det(I+M^TM)}.
\end{equation}
The positively oriented unit tangent $k$-vector is
$\tau_M:=(v_1\wedge\cdots\wedge v_k)/J_M$.
Since $\langle v_a,(e_b,0)\rangle=\delta_{ab}$, we get
\[
 \langle\tau_M,\vec P\rangle
 =\frac{\det(\langle v_a,(e_b,0)\rangle)_{a,b=1}^k}{J_M}
 =\frac{1}{J_M}.
\]
Both $\tau_M$ and $\vec P$ have norm one, so
\begin{equation}\label{eq:tilt-graph-identity}
 |\tau_M-\vec P|^2J_M=2(J_M-1).
\end{equation}
Note that every eigenvalue of $M^TM$ is at most $|M|^2$. Then for $|M|\leq1$, \eqref{eq:tilt-graph-jacobian} gives
\[
 0\leq J_M-1
 \leq(1+|M|^2)^{k/2}-1
 \leq C|M|^2,
\]
for some fixed $C > 0$.
By \eqref{eq:tilt-graph-identity},
\[
 |\tau_M-\vec P|^2J_M=2(J_M-1)\leq 2C|M|^2.
\]

By \cite[Theorem~1.13]{DLSQ}, at almost every base point $X\in A^+$
there exist linear maps $M_{i,1}(X),\ldots,M_{i,q}(X):P\to P^\perp$
that are the derivatives of the $q$ affine functions in the first-order
approximation of $u_i$ at $X$, counting repetitions.
By \cite[Proposition~2.17]{DLSQ},
\begin{equation}\label{eq:q-valued-derivative-norm}
 \sum_{\ell=1}^q|M_{i,\ell}(X)|^2=|du_i(X)|^2,
\end{equation}
with $|du_i|$ from Definition~\ref{def:q-valued}. The Lipschitz bound in \eqref{eq:tangent-approximation-bounds}
gives $|M_{i,\ell}|\leq\sqrt{k}\operatorname{Lip}(u_i)\leq1$
for all sufficiently large $i$.
Using the area formula for the graph current
\cite[(3.29)--(3.30), p.~36]{KWI}, the equality \eqref{eq:tangent-graph-agreement},
and \eqref{eq:tilt-graph-identity}--\eqref{eq:q-valued-derivative-norm},
we obtain
\begin{equation}\label{eq:matching-good-set-tilt}
 \begin{aligned}
 &\int_{(K_i\cap A^+)\times P^\perp}
       |\vec T_i-\vec P|^2\,\dd\|T_i\|\\
 &\quad=\int_{K_i\cap A^+}\sum_{\ell=1}^q
       |\tau_{M_{i,\ell}(X)}-\vec P|^2
       J_{M_{i,\ell}(X)}\,\dd X\\
 &\quad=2\int_{K_i\cap A^+}
       \sum_{\ell=1}^q(J_{M_{i,\ell}(X)}-1)\,\dd X\\
 &\quad\leq C\int_{K_i\cap A^+}|du_i|^2\,\dd X
 \leq C\int_{A^+}|du_i|^2\,\dd X.
 \end{aligned}
\end{equation}
Finally, the strong $L^2$ convergence of
$\varepsilon_i^{-1}|du_i|$ in \eqref{eq:tangent-convergence}
gives a uniform bound for its $L^2(A^+)$ norm. Therefore
\[
 \int_{A^+}|du_i|^2\,\dd X
 =\varepsilon_i^2\int_{A^+}
       \bigl(\varepsilon_i^{-1}|du_i|\bigr)^2\,\dd X
 =O(\varepsilon_i^2).
\]
Thus \eqref{eq:matching-good-set-tilt} bounds the integral over
$(K_i\cap A^+)\times P^\perp$ in the definition of the oriented tilt
excess \eqref{eq:oriented-tilt-excess} by $O(\varepsilon_i^2)$. Off $K_i$, the integrand is bounded
and the mass is $O(\varepsilon_i^{2+2\varsigma})$ by
\eqref{eq:tangent-approximation-bounds}. Consequently,
\begin{equation}\label{eq:matching-tilt-decay}
 t_i(A^+)^2=O(\varepsilon_i^2).
\end{equation}

We first define functions on $A^+$ by
\begin{equation}\label{eq:matching-height-functions}
 h_{i,a}:=\varepsilon_i\psi_a,\qquad
 h_i:=\sum_{a=1}^s\llbracket h_{i,a}\rrbracket.
\end{equation}
We will verify that, after translation and rescaling, $h_i$ satisfies
\cite[Hypothesis~4.16]{KWIII}, as required by
\cite[Theorem~4.19]{KWIII}.
Note that each $h_{i,a}$ is harmonic and
\begin{equation}\label{eq:matching-height-properties}
 \begin{aligned}
 \operatorname{graph}(h_i)
 &=\operatorname{graph}((\varepsilon_i\cdot\widehat\Phi)|_{A^+}),\\
 \|dh_{i,a}\|_{C^0(A^+)}&=O(\varepsilon_i).
 \end{aligned}
\end{equation}
The local branches $\psi_a$ are smooth on a neighborhood of $\overline{A^+}$.
Their derivatives are therefore bounded there, and profile separation on $A^+$
gives a fixed $\kappa<\infty$ such that
\[
 \sup_{A^+}|d\psi_a-d\psi_b|
 \leq\kappa\inf_{A^+}|\psi_a-\psi_b|
 \qquad(a\ne b).
\]
The same inequality holds after multiplication by $\varepsilon_i$.
If $s=1$, there is no pair to check.
After translating and dilating $A^+$ to the unit ball and rotating
$P$ to the coordinate plane, these bounds verify
\cite[Hypothesis~4.16]{KWIII}.
The fixed dilation only changes $\kappa$ by a fixed factor.
Together with \eqref{eq:matching-tilt-decay} and the current properties
above, this verifies all the hypotheses in
\cite[Theorem~4.19]{KWIII} for sufficiently large $i$.

\emph{Step 2: apply the support estimate.}
As in \cite[Definition~4.11, (4.47), p.~49]{KWIII}, for
$\sigma\in(0,1)$ and a fixed open ball $B\subseteq A^+$ of radius
$r_B$, define the nonnegative \defterm{fine excess}
$E_{\sigma,i}(B)$ by
\begin{equation}\label{eq:comparison-fine-excess}
 \begin{split}
 E_{\sigma,i}(B)^2:={}&
 \frac{1}{|B|r_B^2}\int_{B\times P^\perp}
       \dist^2(W,\operatorname{graph}(h_i|_B))\,\dd\|T_i\|(W)\\
 &+\sum_a\|dh_{i,a}\|_{C^0(B)}\,t_i(B)^2
   +t_i(B)^{2+\sigma}.
 \end{split}
\end{equation}
The Laplacian term in \cite[Definition~4.11, (4.47)]{KWIII}
vanishes because $\Delta h_{i,a}=0$ for every $a$.
By \cite[Theorem~4.19, (4.77)]{KWIII},
the hypotheses verified in Step~1 imply that there exists
$\sigma\in(0,1)$, independent of $i$, such that, after translating
and rescaling back to $A^+$, for all sufficiently large $i$ we have
\begin{equation}\label{eq:matching-support-bound}
 \sup_{\substack{W\in\supp T_i\\ \pi_P(W)\in A'}}
       \dist(W,\operatorname{graph}(h_i))
 \leq C r_+ E_{\sigma,i}(A^+),
\end{equation}
where $C$ is independent of $i$.
It remains to estimate the right-hand side. In Step~3 we prove
$E_{\sigma,i}(A^+)=o(\varepsilon_i)$. Since $r_+$ is fixed,
\eqref{eq:matching-support-bound} will then give a distance bound
of $o(\varepsilon_i)$ for every point of $\supp T_i$ over $A'$.

\emph{Step 3: estimate the fine excess.} We claim that $E_{\sigma,i}(A^+)=o(\varepsilon_i)$. We estimate the three terms in \eqref{eq:comparison-fine-excess},
starting with the distance integral. We split its domain into
$(K_i\cap A^+)\times P^\perp$ and
$(A^+\setminus K_i)\times P^\perp$.

For $X\in A^+$, write
$u_i(X)=\sum_{\ell=1}^q\llbracket u_{i,\ell}(X)\rrbracket$,
counting repetitions. Each point $(X,\varepsilon_i\psi_a(X))$
belongs to $\operatorname{graph}(h_i)$ by
\eqref{eq:matching-height-functions}. Therefore
\begin{equation}\label{eq:matching-pointwise-distance}
 \begin{aligned}
 &\sum_{\ell=1}^q
   \dist^2((X,u_{i,\ell}(X)),\operatorname{graph}(h_i))\\
 &\quad\leq
   \sum_{\ell=1}^q\min_{1\leq a\leq s}
       |u_{i,\ell}(X)-\varepsilon_i\psi_a(X)|^2\\
 &\quad\leq
   \mathcal G(u_i(X),\varepsilon_i\cdot\widehat\Phi(X))^2\\
 &\quad=\varepsilon_i^2
   \mathcal G(\varepsilon_i^{-1}\cdot u_i(X),\widehat\Phi(X))^2.
 \end{aligned}
\end{equation}
On $(K_i\cap A^+)\times P^\perp$, the equality of graph currents
\eqref{eq:tangent-graph-agreement} allows us to use the area formula
\cite[(3.29), p.~36]{KWI}. By Step~1, almost everywhere and for all sufficiently large $i$,
\begin{equation}\label{eq:matching-jacobian-bound}
 J_{M_{i,\ell}}
 \leq(1+|M_{i,\ell}|^2)^{k/2}\leq 2^{k/2},
 \qquad
 |M_{i,\ell}|\leq\sqrt{k}\operatorname{Lip}(u_i)\leq1.
\end{equation}
The area formula, \eqref{eq:matching-jacobian-bound}, and
\eqref{eq:matching-pointwise-distance} now yield
\begin{equation}\label{eq:matching-good-set-distance}
 \begin{aligned}
 &\int_{(K_i\cap A^+)\times P^\perp}
       \dist^2(W,\operatorname{graph}(h_i))\,\dd\|T_i\|(W)\\
 &\quad=\int_{K_i\cap A^+}\sum_{\ell=1}^q
       \dist^2((X,u_{i,\ell}(X)),\operatorname{graph}(h_i))
       J_{M_{i,\ell}(X)}\,\dd X\\
 &\quad\leq 2^{k/2}\varepsilon_i^2
       \int_{K_i\cap A^+}
       \mathcal G(\varepsilon_i^{-1}\cdot u_i,\widehat\Phi)^2\,\dd X\\
 &\quad\leq 2^{k/2}\varepsilon_i^2
       \int_{A^+}
       \mathcal G(\varepsilon_i^{-1}\cdot u_i,\widehat\Phi)^2\,\dd X
       =o(\varepsilon_i^2),
 \end{aligned}
\end{equation}
because the first line of \eqref{eq:tangent-convergence}, with
strong $L^2$ convergence understood as in Definition~\ref{def:q-valued}
and restricted to $A^+$, says precisely that the last integral tends to zero.

For the integral over $(A^+\setminus K_i)\times P^\perp$, choose
a fixed $L$ with $A^+\subset B_L^P(0)$. The functions $\psi_a$
are smooth on a neighborhood of $\overline{A^+}$, so
$M_0:=\max_a\|\psi_a\|_{C^0(A^+)}<\infty$.
For $(X,Y)\in\supp T_i$ with $X\in A^+$,
\eqref{eq:magnified-local-bounds} and
\eqref{eq:matching-height-functions} give
\begin{equation}\label{eq:matching-distance-uniform-bound}
 \dist((X,Y),\operatorname{graph}(h_i))
 \leq |Y-h_{i,1}(X)|
 \leq \frac{8L}{7}+\varepsilon_i M_0
 \leq \frac{8L}{7}+M_0
\end{equation}
for all sufficiently large $i$. Multiplying the square of this
bound by the mass estimate in
\eqref{eq:tangent-approximation-bounds} gives
\begin{equation}\label{eq:matching-bad-set-distance}
 \begin{aligned}
 &\int_{(A^+\setminus K_i)\times P^\perp}
       \dist^2(W,\operatorname{graph}(h_i))\,\dd\|T_i\|(W)\\
 &\quad\leq
       \left(\frac{8L}{7}+M_0\right)^2
       \|T_i\|((A^+\setminus K_i)\times P^\perp)\\
 &\quad=O(\varepsilon_i^{2+2\varsigma})
       =o(\varepsilon_i^2).
 \end{aligned}
\end{equation}
Adding \eqref{eq:matching-good-set-distance} and
\eqref{eq:matching-bad-set-distance} proves
\begin{equation}\label{eq:matching-current-l2}
 \int_{A^+\times P^\perp}
       \dist^2(W,\operatorname{graph}(h_i))\,\dd\|T_i\|(W)
 =o(\varepsilon_i^2).
\end{equation}

For the second term in \eqref{eq:comparison-fine-excess}, the
derivative bound in \eqref{eq:matching-height-properties} and
\eqref{eq:matching-tilt-decay}, together with the fixed number $s$
of distinct local branches, give
\begin{equation}\label{eq:matching-derivative-excess-term}
 \sum_{a=1}^s\|dh_{i,a}\|_{C^0(A^+)}\,t_i(A^+)^2
 \leq C\varepsilon_i\,C\varepsilon_i^2
 =O(\varepsilon_i^3).
\end{equation}
For the third term, \eqref{eq:matching-tilt-decay} gives
\begin{equation}\label{eq:matching-excess-power}
 t_i(A^+)^{2+\sigma}
 =\bigl(t_i(A^+)^2\bigr)^{1+\sigma/2}
 =O(\varepsilon_i^{2+\sigma}).
\end{equation}
Substituting \eqref{eq:matching-current-l2},
\eqref{eq:matching-derivative-excess-term}, and
\eqref{eq:matching-excess-power} into
\eqref{eq:comparison-fine-excess}, and using that
$|A^+|$ and $r_+$ are fixed and positive, we obtain
\begin{equation}\label{eq:matching-fine-excess-squared}
 \begin{aligned}
 E_{\sigma,i}(A^+)^2
 &\leq o(\varepsilon_i^2)
       +C\varepsilon_i^3+C\varepsilon_i^{2+\sigma} = o(\varepsilon_i^2).
 \end{aligned}
\end{equation}
Hence,
\begin{equation}\label{eq:matching-fine-excess-decay}
 E_{\sigma,i}(A^+)=o(\varepsilon_i),
\end{equation}
which finishes the claim of this step.

\emph{Step 4: construct the currents and identify their projection multiplicities.}
We now write $T_i$ over $A'$ as a sum of currents $S_{i,a}$
with pairwise disjoint supports and prove that they satisfy
\eqref{eq:multiplicity-matching}.

We first substitute the conclusion
\eqref{eq:matching-fine-excess-decay} into the support estimate
\eqref{eq:matching-support-bound} and then bound
$|Y-h_{i,a}(X)|$ for a suitable branch $a$ at each
$(X,Y)\in\supp T_i$ with $X\in A'$.
Let $C$ be the constant in \eqref{eq:matching-support-bound} and set
\begin{equation}\label{eq:matching-distance-tolerance}
 d_i:=Cr_+E_{\sigma,i}(A^+)+\varepsilon_i^2
       =o(\varepsilon_i),
\end{equation}
where the equality follows from
\eqref{eq:matching-fine-excess-decay}.
For each $(X,Y)\in\supp T_i$ with $X\in A'$, the definition of
distance and \eqref{eq:matching-support-bound} allow us to choose
$(X',h_{i,a}(X'))\in\operatorname{graph}(h_i)$ such that
\begin{equation}\label{eq:matching-nearby-graph-point}
 |X-X'|<d_i,\qquad |Y-h_{i,a}(X')|<d_i.
\end{equation}
Since $X,X'\in A^+$ and $A^+$ is a ball, the segment joining them
lies in $A^+$. Integrating $dh_{i,a}$ along this segment and using
\eqref{eq:matching-nearby-graph-point} and
\eqref{eq:matching-height-properties}, we find
\begin{equation}\label{eq:matching-same-base-distance}
 \begin{aligned}
 |Y-h_{i,a}(X)|
 &\leq |Y-h_{i,a}(X')|+|h_{i,a}(X')-h_{i,a}(X)|\\
 &\leq d_i+\|dh_{i,a}\|_{C^0(A^+)}|X-X'|\\
 &\leq (1+O(\varepsilon_i))d_i=o(\varepsilon_i).
 \end{aligned}
\end{equation}
The bound in \eqref{eq:matching-same-base-distance} is independent
of $(X,Y)$ and of the chosen branch $a$. Recalling
\eqref{eq:matching-height-functions}, we conclude that
\begin{equation}\label{eq:matching-vertical-bound}
 \sup_{\substack{(X,Y)\in\supp T_i\\X\in A'}}
       \min_{1\leq a\leq s}|Y-\varepsilon_i\psi_a(X)|
 =o(\varepsilon_i).
\end{equation}

We now construct the currents.
Let
\[
 d_A:=\min_{a<b}\inf_{\overline{A^+}}|\psi_a-\psi_b|>0
\]
when $s>1$, and set $d_A=1$ when $s=1$.
Define the disjoint open sets
\[
 U_{i,a}:=\{(X,Y):X\in A',\
       |Y-\varepsilon_i\psi_a(X)|<d_A\varepsilon_i/4\},
\]
and set $S_{i,a}:=T_i\llcorner U_{i,a}$.
By \eqref{eq:matching-vertical-bound}, every point of
$\supp T_i\cap(A'\times P^\perp)$ lies in $\bigcup_a U_{i,a}$
for all sufficiently large $i$.
The boundaries of the sets $U_{i,a}$ inside $A'\times P^\perp$
do not meet $\supp T_i$. Therefore
\[
 T_i\llcorner(A'\times P^\perp)=\sum_{a=1}^s S_{i,a},
\]
where the $S_{i,a}$ are integral currents with pairwise disjoint
supports and no boundary in $A'\times P^\perp$.
Each $S_{i,a}$ has bounded normal support because it is a restriction
of $T_i$, whose normal support over $A'$ is bounded by
\eqref{eq:magnified-local-bounds}.

  Each  $S_{i,a}$ is locally area minimizing in $A'\times P^\perp$, since for every open $W\Subset A'\times P^\perp$ and every integral
$(k+1)$-current $Q$ with $\supp Q\Subset W$, disjointness of the supports,
local area minimization of $T_i$, and subadditivity of mass give
\[
 \begin{aligned}
 \sum_{b=1}^s\mathbf M(S_{i,b}\llcorner W)
 &=\mathbf M(T_i\llcorner W)\\
 &\leq\mathbf M((T_i+\partial Q)\llcorner W)\\
 &\leq\mathbf M((S_{i,a}+\partial Q)\llcorner W)
       +\sum_{b\ne a}\mathbf M(S_{i,b}\llcorner W).
 \end{aligned}
\]
Subtracting the terms with $b\ne a$ gives the required comparison.
For points in $\supp S_{i,a}$, branch separation forces the minimum
in \eqref{eq:matching-vertical-bound} to be attained at $a$. Hence
\[
 \sup_{\substack{(X,Y)\in\supp S_{i,a}\\X\in A'}}
       |Y-\varepsilon_i\psi_a(X)|=o(\varepsilon_i).
\]

It remains to identify the projection multiplicities.
The constancy theorem \cite[Chapter~6, Theorem~2.41, p.~168]{SimonGMT} gives
\begin{equation}\label{eq:matching-projection-constant}
 (\pi_P)_\#S_{i,a}=d_{i,a}\llbracket A'\rrbracket,
 \qquad d_{i,a}\in\mathbb Z,
\end{equation}
since $A'$ is connected, $S_{i,a}$ has no boundary in the cylinder,
and its normal support is bounded. We show that $d_{i,a}=n_a$ by
counting the points of the current above suitable base points.

Fix positively oriented orthonormal coordinates $x^1,\ldots,x^k$
on $P$, and also write $x^j$ for the function $x^j\circ\pi_P$
on $P\times P^\perp$. For $X\in A'$, write $X^j:=x^j(X)$.
For any $k$-dimensional integral current $S$ in $A'\times P^\perp$,
we apply the slicing construction of
\cite[Chapter~6, Definition~4.4, p.~181]{SimonGMT}
successively to $S$ by $x^1$ at $X^1$,
then by $x^2$ at $X^2$, and so on, through $x^k$ at $X^k$.
By \cite[Chapter~6, Lemmas~4.1 and~4.5, pp.~180--182]{SimonGMT},
each successive slice is an integral current for almost every
value of the next coordinate. Thus the construction gives a
zero-dimensional integral current for almost every $X\in A'$.
We call this current the \defterm{projection slice}
of $S$ at $X$ and denote it by $\langle S,\pi_P,X\rangle$.
Note that this current is supported in
\[
 \{x^1=X^1,\ldots,x^k=X^k\}=\{X\}\times P^\perp.
\]
The slice has locally finite mass by
\cite[Chapter~6, Lemma~4.5(1), p.~182]{SimonGMT},
applied at each successive slicing step.
We now apply this construction to $S=S_{i,a}$.
Since $S_{i,a}$ has bounded normal support, its projection slice
has compact support in the cylinder. Its locally finite mass and
integer coefficients therefore imply that it is a finite sum
\begin{equation}\label{eq:matching-projection-slice}
 \langle S_{i,a},\pi_P,X\rangle
 =\sum_{j=1}^{N(X)}b_j(X)\llbracket(X,Y_j(X))\rrbracket,
 \qquad b_j(X)\in\mathbb Z, \qquad Y_j(X) \in P^\perp,
\end{equation}
where $\llbracket(X,Y_j(X))\rrbracket$ is the current,
which acts on a test function by evaluation at $(X,Y_j(X))$.
The integers $b_j(X)$ may depend on $i$ and $a$. We call $\sum_j b_j(X)$ the \defterm{total signed multiplicity}
of this projection slice.

We claim that, for almost every $X\in A'$, the total signed
multiplicity $\sum_{j=1}^{N(X)}b_j(X)$ of the projection slice
equals the integer $d_{i,a}$ in
\eqref{eq:matching-projection-constant}.
To prove this, let $\omega_P$ be the constant volume
$k$-form of the oriented Euclidean plane $P$.
Applying the slice definition in
\cite[Chapter~6, Definition~4.4, p.~181]{SimonGMT}
and the coarea formula in
\cite[Chapter~6, Lemma~4.1, p.~180]{SimonGMT}
successively to the coordinates $x^1,\ldots,x^k$,
and using the definition of pushforward, gives, for every
$\eta\in C_c^\infty(A')$,
\begin{equation}\label{eq:matching-slice-projection}
 \int_{A'}\eta(X)\sum_{j=1}^{N(X)}b_j(X)\,\dd X
 =((\pi_P)_\#S_{i,a})(\eta\omega_P)
 =d_{i,a}\int_{A'}\eta(X)\,\dd X,
\end{equation}
where the last equality is \eqref{eq:matching-projection-constant}.
Thus
\begin{equation}\label{eq:matching-slice-signed-count}
 \sum_{j=1}^{N(X)}b_j(X)=d_{i,a}
 \quad\text{for almost every }X\in A'.
\end{equation}

We next use the graph representation on $K_i$ to show that $d_{i,a}=n_a$.
For almost every base point $X$, the projection slice of a graph
current consists of its points above $X$, counted with repetition
multiplicity and positive sign, by the graph orientation in
Definition~\ref{def:q-valued}.
Slicing respects restriction to $K_i\times P^\perp$, so
\eqref{eq:tangent-graph-agreement} gives, for almost every
$X\in K_i\cap A'$,
\begin{equation}\label{eq:matching-graph-slice}
 \langle T_i,\pi_P,X\rangle
 =\langle\mathbf G_{u_i},\pi_P,X\rangle
 =\sum_{\ell=1}^q\llbracket(X,u_{i,\ell}(X))\rrbracket,
\end{equation}
where $u_i(X)=\sum_{\ell=1}^q\llbracket u_{i,\ell}(X)\rrbracket$,
counting repetitions.
Since $S_{i,a}=T_i\llcorner U_{i,a}$, the projection slice of $S_{i,a}$ keeps exactly
the summands of \eqref{eq:matching-graph-slice} lying in $U_{i,a}$.
Combining this with \eqref{eq:matching-slice-signed-count}, we obtain
\begin{equation}\label{eq:matching-slice-count}
 d_{i,a}
 =\#\{\ell\in\{1,\ldots,q\}:
          (X,u_{i,\ell}(X))\in U_{i,a}\}
 \quad\text{for almost every }X\in K_i\cap A'.
\end{equation}
To determine this number, define the measurable set
\begin{equation}\label{eq:matching-exceptional-base-set}
 D_i:=\left\{X\in A':
   \mathcal G(\varepsilon_i^{-1}\cdot u_i(X),\widehat\Phi(X))
       \geq\frac{d_A}{8}\right\}.
\end{equation}
By the first line of \eqref{eq:tangent-convergence} and the
definition of strong $L^2$ convergence in
Definition~\ref{def:q-valued},
\begin{equation}\label{eq:matching-exceptional-base-measure}
 |D_i|\leq\frac{64}{d_A^2}
       \int_{A'}\mathcal G(\varepsilon_i^{-1}\cdot u_i,\widehat\Phi)^2
       \,\dd X\longrightarrow0.
\end{equation}
Set $V_i:=(K_i\cap A')\setminus D_i$. Equations
\eqref{eq:tangent-approximation-bounds} and
\eqref{eq:matching-exceptional-base-measure} imply
\begin{equation}\label{eq:matching-counting-base-set}
 |A'\setminus V_i|
 \leq |A'\setminus K_i|+|D_i|\longrightarrow0.
\end{equation}
Since $|A'|>0$, the set $V_i$ has positive measure for all
sufficiently large $i$.

For $X\in V_i$, choose an ordering
$\xi_1(X),\ldots,\xi_q(X)$ of the values of $\widehat\Phi(X)$,
including repetitions, that attains the minimum in the definition
of $\mathcal G(\varepsilon_i^{-1}\cdot u_i(X), \widehat\Phi(X))$. Then \eqref{eq:matching-exceptional-base-set}
gives
\begin{equation}\label{eq:matching-value-pairing}
 |u_{i,\ell}(X)-\varepsilon_i\xi_\ell(X)|
 \leq\varepsilon_i
       \mathcal G(\varepsilon_i^{-1}\cdot u_i(X),\widehat\Phi(X))
 <\frac{d_A\varepsilon_i}{8}
 \qquad(\ell=1,\ldots,q).
\end{equation}
Exactly $n_a$ of the values $\xi_\ell(X)$ equal $\psi_a(X)$.
For each of those indices, \eqref{eq:matching-value-pairing}
implies
\[
 |u_{i,\ell}(X)-\varepsilon_i\psi_a(X)|
 <\frac{d_A\varepsilon_i}{8}
 <\frac{d_A\varepsilon_i}{4},
\]
so $(X,u_{i,\ell}(X))\in U_{i,a}$.
For every other index, $\xi_\ell(X)=\psi_b(X)$ for some $b\ne a$,
and the definition of $d_A$ and
\eqref{eq:matching-value-pairing} give
\[
 \begin{aligned}
 |u_{i,\ell}(X)-\varepsilon_i\psi_a(X)|
 &\geq\varepsilon_i|\psi_b(X)-\psi_a(X)|
       -|u_{i,\ell}(X)-\varepsilon_i\psi_b(X)|\\
 &>\frac{7d_A\varepsilon_i}{8}
  >\frac{d_A\varepsilon_i}{4}.
 \end{aligned}
\]
Such a point does not belong to $U_{i,a}$. Hence exactly $n_a$
indices are counted on the right-hand side of
\eqref{eq:matching-slice-count}.

For all sufficiently large $i$, we can choose $X\in V_i$
outside the null sets where the slice identities fail, because
$V_i$ has positive measure by \eqref{eq:matching-counting-base-set}.
Equation~\eqref{eq:matching-slice-count} then gives $d_{i,a}=n_a$
for every $a$. Substituting into
\eqref{eq:matching-projection-constant} proves the projection
identity in \eqref{eq:multiplicity-matching}. Together with the
support estimate obtained from \eqref{eq:matching-vertical-bound},
this completes the proof.
\end{proof}

\subsection{Exact graphs and convergence of derivatives}
When the projection multiplicity is one, Allard regularity and interior
estimates for the minimal graph equation give an exact sheet-height
function for $T_i$, in the sense of Definition~\ref{def:chart}.
Its height normalized by the cylindrical height excess, as defined in
Section~\ref{sec:normalized-heights}, converges in every fixed $C^b(A)$
norm to the chosen local branch.
We record this standard consequence in Lemma~\ref{lem:extract}, with the details
needed to apply the local curvature obstruction.
\begin{lemma}[Multiplicity-one graph extraction]\label{lem:extract}
In Lemma~\ref{lem:matching}, suppose $n_a=1$.
On the smaller fixed annular sheet chart $A$, for all sufficiently
large $i$ there is an exact sheet-height function $w_i:A\to P^\perp$
for $T_i$ such that
\[
 S_{i,a}\llcorner(A\times P^\perp)=\mathbf G_{w_i},
\]
where $\mathbf G_{w_i}$ is the graph current from
Definition~\ref{def:q-valued}. Its graph is minimal.
For every fixed integer $b\geq0$,
\begin{equation}\label{eq:normalized-sheet-convergence}
 \|w_i\|_{C^b(A)}\leq C_b\varepsilon_i,\qquad
 \varepsilon_i^{-1}w_i\longrightarrow\psi_a
       \quad\hbox{in }C^b(A).
\end{equation}
The constants $C_b$ are independent of $i$. Each graph is an open
subset of $\reg T_i$.
\end{lemma}
\begin{proof}
Fix the index $a$ with $n_a=1$. We first identify the limit of
$S_{i,a}$, then obtain a smooth graph by Allard regularity.
Finally, we estimate its derivatives and prove convergence after
division by $\varepsilon_i$.

\emph{Step 1: identify the limit current.}
Choose a fixed $L>0$ with $A'\subset B_L^P(0)$.
Since $S_{i,a}$ is a restriction of $T_i$,
\eqref{eq:magnified-local-bounds} gives, for every relatively compact
open set $W\Subset A'\times P^\perp$ and all sufficiently large $i$,
\begin{equation}\label{eq:extraction-local-mass}
 \|S_{i,a}\|(W)
 \leq\|T_i\|\bigl(B_L^P(0)\times P^\perp\bigr)
 \leq C_L,\qquad \|\partial S_{i,a}\|(W)=0,
\end{equation}
where the last equality follows from Lemma~\ref{lem:matching}.
Thus the Federer--Fleming compactness theorem
\cite[Chapter~6, Theorem~3.15, p.~175]{SimonGMT}
applies: every subsequence has a further subsequence converging
weakly to an integral current $S_\infty$ in $A'\times P^\perp$.
For the moment, work along such a subsequence and keep the index $i$.
The limit has no boundary, since for every smooth compactly supported
$(k-1)$-form $\zeta$ in the cylinder,
\[
 \partial S_\infty(\zeta)
 =S_\infty(\dd\zeta)
 =\lim_i S_{i,a}(\dd\zeta)
 =\lim_i\partial S_{i,a}(\zeta)=0.
\]

We now show that $S_\infty$ is supported in $A'\times\{0\}$.
The function $\psi_a$ is bounded on $A'$, so the support estimate in
\eqref{eq:multiplicity-matching} gives
\begin{equation}\label{eq:extraction-normal-height}
 \sup_{(X,Y)\in\supp S_{i,a}}|Y|
 \leq\varepsilon_i\|\psi_a\|_{C^0(A')}+o(\varepsilon_i)
 \longrightarrow0.
\end{equation}
Any compact subset of the cylinder disjoint from $A'\times\{0\}$
has positive distance from this plane. By
\eqref{eq:extraction-normal-height}, a test form supported in such
a set evaluates to zero on $S_{i,a}$ for all sufficiently large $i$,
and hence on $S_\infty$. Therefore
\[
 \supp S_\infty\subset A'\times\{0\}.
\]

We next pass the projection identity in
\eqref{eq:multiplicity-matching}, with $n_a=1$, to the limit.
Let $\eta$ be a smooth compactly supported $k$-form on $A'$.
Choose $\chi\in C_c^\infty(P^\perp)$ equal to one on a fixed
neighborhood of $0$. By \eqref{eq:extraction-normal-height}, this
neighborhood contains every normal coordinate of $\supp S_{i,a}$
for all sufficiently large $i$. Since
$\supp S_\infty\subset A'\times\{0\}$, the function $\chi(Y)$
is therefore equal to one on the supports of $S_\infty$ and of
all sufficiently large $S_{i,a}$. The form $\chi(Y)\pi_P^*\eta$
is compactly supported in $A'\times P^\perp$, so it can be used
to test weak convergence. As in \eqref{eq:projected-current-limit},
we obtain
\begin{equation}\label{eq:extraction-projection-limit}
 \begin{aligned}
 ((\pi_P)_\#S_\infty)(\eta)
 &=S_\infty\bigl(\chi(Y)\pi_P^*\eta\bigr)\\
 &=\lim_i S_{i,a}\bigl(\chi(Y)\pi_P^*\eta\bigr)\\
 &=\lim_i((\pi_P)_\#S_{i,a})(\eta)
 =\llbracket A'\rrbracket(\eta).
 \end{aligned}
\end{equation}
Since $A'$ is connected and $S_\infty$ is supported in its plane
and has no boundary, the constancy theorem
\cite[Chapter~6, Theorem~2.41, p.~168]{SimonGMT}
gives $S_\infty=d\llbracket A'\rrbracket$ for an integer $d$.
Here the current on the right is viewed in $A'\times P^\perp$
by the inclusion $X\mapsto(X,0)$.
Equation~\eqref{eq:extraction-projection-limit} forces $d=1$.

As in the proof of Lemma~\ref{lem:localization}, compactness and
uniqueness of the subsequential limit imply that the whole sequence
$S_{i,a}$ converges weakly to $\llbracket A'\rrbracket$ in
$A'\times P^\perp$.
Since the currents are locally area minimizing,
\cite[Chapter~7, Theorem~2.4 and Remark~2.5(1), pp.~214--215]{SimonGMT}
also give local weak convergence of their mass measures and of
their associated varifolds. In particular,
\begin{equation}\label{eq:extraction-plane-convergence}
 S_{i,a}\to\llbracket A'\rrbracket,\qquad
 \|S_{i,a}\|\to
 \mathcal H^k\llcorner(A'\times\{0\})
 \quad\text{weakly locally in }A'\times P^\perp.
\end{equation}
The limiting plane thus has multiplicity one, also as a varifold.

\emph{Step 2: obtain a smooth graph.}
Choose an open ball $A''$ with $A\Subset A''\Subset A'$.
The varifolds associated with $S_{i,a}$ are stationary in
$A'\times P^\perp$, in the sense of
\cite[Chapter~4, Definition~2.4, p.~106]{SimonGMT}, which follows from \cite[Chapter~7, Lemma~1.2, p.~212]{SimonGMT},
since $S_{i,a}$ is locally area minimizing and has no boundary
there by Lemma~\ref{lem:matching}.
Since $S_{i,a}$ are integral, the associated varifolds also have integer multiplicities at least one almost everywhere.
{We next show that the mass ratios in
\eqref{eq:extraction-mass-ratio} converge to one.
Fix an ambient ball $B:=B_r^{\R^{k+m}}((X,0))$, with $X\in A'$
and $\overline B\subset A'\times P^\perp$.
The limiting measure in \eqref{eq:extraction-plane-convergence}
is $\mathcal H^k\llcorner(A'\times\{0\})$. Its mass on $B$ is
$\omega_k r^k$. Its mass on $\partial B$ is
zero, since $\partial B\cap(A'\times\{0\})$ is a
$(k-1)$-dimensional sphere and therefore has zero
$\mathcal H^k$ measure.

We use this boundary condition to pass from weak convergence to
convergence of the masses of $B$. The indicator function
$\mathbf1_B$ is discontinuous on $\partial B$, so it cannot
itself be used as a test function for weak convergence.
Approximate it from below and above by continuous functions
compactly supported in $A'\times P^\perp$. Since $\partial B$
has zero limiting mass, the integrals of these approximations
against the limiting measure tend to the same value
$\omega_k r^k$. Applying
\eqref{eq:extraction-plane-convergence} to the approximations
and using the two resulting bounds gives
\begin{equation}\label{eq:extraction-mass-ratio}
 \frac{\|S_{i,a}\|(B_r^{\R^{k+m}}((X,0)))}{\omega_k r^k}
 \longrightarrow1.
\end{equation}}
The same convergence holds with centers
\[
 W_i\in\supp S_{i,a},\qquad W_i\longrightarrow(X,0),
\]
because comparison with concentric balls centered at $(X,0)$
of slightly larger and smaller radii gives
\begin{equation}\label{eq:extraction-moving-center-mass-ratio}
 \frac{\|S_{i,a}\|(B_r^{\R^{k+m}}(W_i))}{\omega_k r^k}
 \longrightarrow1.
\end{equation}
Such points $W_i\in\supp S_{i,a}$ exist by the positive limiting
mass in every ball around $(X,0)$.

We have now verified the hypotheses of Allard's regularity
theorem in formulation of
\cite[Chapter~5, hypotheses~(5.1) and Theorem~5.2,
pp.~145--146]{SimonGMT}: the integer multiplicities are at least
one almost everywhere, the centers lie in the support, and
stationarity from \cite[Chapter~7, Lemma~1.2, p.~212]{SimonGMT}
makes the generalized mean curvature zero
\cite[Chapter~4, Definition~3.15, p.~111]{SimonGMT}.
The mass ratios centered at $W_i\in\supp S_{i,a}$ satisfy
the required smallness condition for all sufficiently large $i$
by \eqref{eq:extraction-moving-center-mass-ratio}.
Applying this theorem at $W_i$ and decreasing the graph domains
by a fixed factor gives a radius $s>0$ and an exponent
$\alpha\in(0,1)$, both independent of $i$,
orthogonal transformations $Q_i$ of $\R^{k+m}$, and functions
$v_i\in C^{1,\alpha}(B_{2s}^{\R^k}(0);\R^m)$
with $v_i(0)=0$ and $dv_i(0)=0$ such that the maps
\begin{equation}\label{eq:extraction-local-parametrization}
 F_i:B_{2s}^{\R^k}(0)\longrightarrow\R^{k+m},\qquad
 F_i(\xi):=W_i+Q_i(\xi,v_i(\xi)),
\end{equation}
parametrize neighborhoods of $W_i$ in $\supp S_{i,a}$.
Allard's theorem also gives a fixed $b\in(0,s)$ such that
\begin{equation}\label{eq:extraction-all-support-ball}
 \supp S_{i,a}\cap B_b(W_i)\subset F_i(B_s)
 \qquad\text{for all sufficiently large }i.
\end{equation}
The estimate in the theorem gives
\begin{equation}\label{eq:extraction-rotated-graph-bounds}
 \|dv_i\|_{C^0(B_{2s})}
 +\sup_{\substack{\xi,\eta\in B_{2s}\\\xi\ne\eta}}
       \frac{|dv_i(\xi)-dv_i(\eta)|}{|\xi-\eta|^\alpha}
 \leq C,
 \qquad
 \|v_i\|_{C^0(B_{2s})}\leq2s\|dv_i\|_{C^0(B_{2s})},
\end{equation}
where $C$ is independent of $i$ and the second inequality follows
from $v_i(0)=0$ by integration along radial segments.

We next show that the tangent planes of these graphs converge
uniformly to $P$ on the smaller domains $B_s$, in the sense that
\[
 \sup_{\xi\in B_s}
 \|\pi_{T_{i,\xi}}-\pi_P\|\longrightarrow0,
 \qquad T_{i,\xi}:=dF_i|_\xi(\R^k).
\]
Here $T_{i,\xi}$ is the tangent plane to the graph at $F_i(\xi)$,
$\pi_{T_{i,\xi}}$ is its orthogonal projection.
By compactness of the orthogonal group, every subsequence has
a further subsequence on which $Q_i$ converges.
By \eqref{eq:extraction-rotated-graph-bounds} and the
Arzel\`a--Ascoli theorem applied to $(v_i,dv_i)$ on
$\overline{B_s}$, we can pass to a further subsequence on which
$v_i$ converges in $C^1(\overline{B_s})$.
Thus, along the same subsequence,
\[
 Q_i\longrightarrow Q_\infty,\qquad
 v_i\longrightarrow v_\infty\quad\text{in }C^1(\overline{B_s}).
\]
Since $F_i(B_s)\subset\supp S_{i,a}$, the support estimate
\eqref{eq:extraction-normal-height} implies
\[
 F_\infty(\xi):=(X,0)+Q_\infty(\xi,v_\infty(\xi))
 \in P\qquad(\xi\in B_s).
\]
Differentiating at $0$ and using $dv_\infty(0)=0$ gives
\[
 dF_\infty|_0(\zeta)=Q_\infty(\zeta,0)\in P
 \qquad(\zeta\in\R^k),
\]
where the inclusion follows from $F_\infty(B_s)\subset P$.
Since $Q_\infty$ is invertible, $Q_\infty(\R^k\times\{0\})$
has dimension $k$, so we get that $Q_\infty(\R^k\times\{0\})=P$.
Orthogonality then gives
$Q_\infty(\{0\}\times\R^m)=P^\perp$, so
$F_\infty(B_s)\subset P$ forces $v_\infty=0$.
This argument applies to every subsequence $(v_i,Q_i)$ and therefore proves
\begin{equation}\label{eq:extraction-rotated-graph-decay}
 \|\pi_{Q_i(\R^k\times\{0\})}-\pi_P\|\longrightarrow0,
 \qquad
 \|dv_i\|_{C^0(B_s)}\longrightarrow0.
\end{equation}
The tangent plane at $F_i(\xi)$ is
\[
 T_{i,\xi}
 =Q_i\{(\zeta,dv_i(\xi)\zeta):\zeta\in\R^k\}
 =Q_i\operatorname{graph}(dv_i(\xi)).
\]
To estimate its orthogonal projection, write
$E:=\R^k\times\{0\}$. For a linear map
$M:\R^k\to\R^m$, the projection onto its graph is
\[
 \pi_{\operatorname{graph}(M)}
 =\begin{pmatrix}I\\M\end{pmatrix}
   (I+M^TM)^{-1}
   \begin{pmatrix}I&M^T\end{pmatrix}.
\]
With our Euclidean norm of the matrix entries, this gives
\begin{equation}\label{eq:extraction-linear-graph-projection}
 \|\pi_{\operatorname{graph}(M)}-\pi_E\|^2
 =2\operatorname{tr}\bigl(M^TM(I+M^TM)^{-1}\bigr)
 \leq2|M|^2,
\end{equation}
where $\operatorname{tr}$ denotes the trace. The last inequality follows
by diagonalizing $M^TM$.
Since $\pi_{Q_i(E)}=Q_i\pi_EQ_i^T$ and the matrix norm is
invariant under orthogonal transformations, the triangle
inequality and \eqref{eq:extraction-linear-graph-projection} give
\[
 \begin{aligned}
 \|\pi_{T_{i,\xi}}-\pi_P\|
 &\leq
 \|\pi_{\operatorname{graph}(dv_i(\xi))}-\pi_E\|
 +\|\pi_{Q_i(E)}-\pi_P\|\\
 &\leq\sqrt2\,|dv_i(\xi)|
 +\|\pi_{Q_i(E)}-\pi_P\|.
 \end{aligned}
\]
Taking the supremum over $\xi$ therefore gives
\[
 \sup_{\xi\in B_s}\|\pi_{T_{i,\xi}}-\pi_P\|
 \leq\sqrt2\,\|dv_i\|_{C^0(B_s)}
      +\|\pi_{Q_i(E)}-\pi_P\|
 \longrightarrow0
\]
by \eqref{eq:extraction-rotated-graph-decay}, which proves the
claimed uniform convergence of tangent planes.

Now, we claim that, for all sufficiently large $i$, these local graphs
can be expressed over a fixed open ball $B\subset P$ centered at
$X$. To prove this, set
\begin{equation}\label{eq:extraction-projected-maps}
 H_i:=\pi_P\circ F_i,\qquad
 N_i:=\pi_{P^\perp}\circ F_i,\qquad
 L_i\xi:=\pi_P Q_i(\xi,0).
\end{equation}
By \eqref{eq:extraction-rotated-graph-decay}, the smallest
singular value of $L_i:\R^k\to P$, which is $\sigma_{\min}(L_i)  :=\min_{\substack{\zeta\in\R^k\\|\zeta|=1}}|L_i\zeta|$, approaches $1$.

By definition, for $\xi\in B_{2s}$, we have
\begin{equation}\label{eq:extraction-base-map-expansion}
 H_i(\xi) =\pi_P(W_i)+L_i\xi+\pi_PQ_i(0,v_i(\xi)).
\end{equation}
For $\xi,\eta\in B_s$, we have
\begin{equation}\label{eq:extraction-base-map-difference}
 \begin{aligned}
 |H_i(\xi)-H_i(\eta)|
 &=|L_i(\xi-\eta)+\pi_PQ_i(0,v_i(\xi)-v_i(\eta))|\\
 &\geq |L_i(\xi-\eta)|-|v_i(\xi)-v_i(\eta)|.
 \end{aligned}
\end{equation}
Since $\sigma_{\min}(L_i)\to1$, for all sufficiently large
$i$ we have
\begin{equation}\label{eq:extraction-linear-map-lower-bound}
 |L_i(\xi-\eta)|\geq\tfrac12|\xi-\eta|.
\end{equation}
The segment from $\eta$ to $\xi$ lies in $B_s$, since $B_s$ is
convex. Thus
\[
 v_i(\xi)-v_i(\eta)
 =\int_0^1dv_i\bigl(\eta+t(\xi-\eta)\bigr)[\xi-\eta]\,\dd t,
\]
and hence
\begin{equation}\label{eq:extraction-height-difference}
 |v_i(\xi)-v_i(\eta)|
 \leq\|dv_i\|_{C^0(B_s)}|\xi-\eta|.
\end{equation}
Combining \eqref{eq:extraction-base-map-difference},
\eqref{eq:extraction-linear-map-lower-bound}, and
\eqref{eq:extraction-height-difference} gives, for all sufficiently
large $i$,
\begin{equation}\label{eq:extraction-base-map-lower-bound}
 |H_i(\xi)-H_i(\eta)|
 \geq\bigl(\tfrac12-\|dv_i\|_{C^0(B_s)}\bigr)|\xi-\eta|
 \geq\tfrac14|\xi-\eta|,
\end{equation}
where the last inequality uses
$\|dv_i\|_{C^0(B_s)}\to0$ from
\eqref{eq:extraction-rotated-graph-decay}.
Thus, for all sufficiently large $i$, each $H_i$ is one-to-one
and has invertible differential,
with
\begin{equation}\label{eq:extraction-inverse-differential-bound}
 \bigl|(dH_i|_\xi)^{-1}\zeta\bigr|
 \leq4|\zeta|\qquad(\xi\in B_s,\ \zeta\in P).
\end{equation}
By continuity of $H_i$, inequality~\eqref{eq:extraction-base-map-lower-bound}
also holds for $\xi,\eta\in\overline{B_s}$. Thus every point of $H_i(\partial B_s)$ is at least $s/4$ from $H_i(0)$,
for all sufficiently large $i$.
The inverse function theorem shows that $H_i(B_s)$ is open.
Together with continuity of $H_i$ and compactness of $\overline{B_s}$, this gives
\[
 \partial\bigl(H_i(B_s)\bigr)\subset H_i(\partial B_s).
\]
It follows that
$B_{s/4}^P(H_i(0))\subset H_i(B_s)$.
Since $H_i(0)=\pi_P(W_i)\to X$ as $i\to\infty$, we therefore have
\[
 B_{s/8}^P(X)\subset H_i(B_s)
 \qquad\text{for all sufficiently large }i.
\]
Let $B:=B_\tau^P(X)$, where
$0<\tau<\min\{s/8,b/4\}$ and $b$ is from
\eqref{eq:extraction-all-support-ball}. Then
\[
 B\subset H_i(B_s)
 \qquad\text{for all sufficiently large }i.
\]
We define
\[
 w_i^{\mathrm{loc}}:B\longrightarrow P^\perp,\qquad
 w_i^{\mathrm{loc}}:=N_i\circ H_i^{-1}.
\]
Writing $\delta_i:=\|\pi_{Q_i(E)}-\pi_P\|$, the formula
\eqref{eq:extraction-projected-maps} gives
\[
 |dN_i(\xi)[\zeta]|
 \leq\bigl(\delta_i+\|dv_i\|_{C^0(B_s)}\bigr)|\zeta|.
\]
The chain rule and \eqref{eq:extraction-inverse-differential-bound} give
\begin{equation}\label{eq:extraction-local-derivative-bound}
 \left|dw_i^{\mathrm{loc}}\big|_{H_i(\xi)}(z)\right|
 \leq4\bigl(\delta_i+\|dv_i\|_{C^0(B_s)}\bigr)|z|
 \qquad(\xi\in H_i^{-1}(B),\ z\in P).
\end{equation}
Since $\delta_i+\|dv_i\|_{C^0(B_s)}\to0$ by
\eqref{eq:extraction-rotated-graph-decay}, this estimate implies
that $dw_i^{\mathrm{loc}}\to0$ uniformly on $B$.
Moreover, the uniform H\"older bounds in
\eqref{eq:extraction-rotated-graph-bounds} and the uniform
Lipschitz bounds for $H_i^{-1}$ give uniform
$C^{1,\alpha}$ bounds for $w_i^{\mathrm{loc}}$.

To see that the graph of $w_i^{\mathrm{loc}}$ contains all support of $S_{i,a}$ above $B$, note that for large $i$, $|W_i-(X,0)|<b/4$ and the bound \eqref{eq:extraction-normal-height} is less than $b/4$.
Every $(Y,Z)\in\supp S_{i,a}$ with $Y\in B$ then satisfies
$|(Y,Z)-W_i|<3b/4$.
By \eqref{eq:extraction-all-support-ball}, $(Y,Z)$ belongs to $F_i(B_s)$.
Writing $(Y,Z)=F_i(\xi)$ with $\xi\in B_s$, we have
$H_i(\xi)=Y$, so injectivity of $H_i$ gives
\[
 Z=N_i(\xi)=N_i(H_i^{-1}(Y))=w_i^{\mathrm{loc}}(Y).
\]
Thus, for all sufficiently large $i$,
\[
 \supp S_{i,a}\cap(B\times P^\perp)
 =\operatorname{graph}(w_i^{\mathrm{loc}}).
\]

We next combine these local graphs into a single graph over $A''\Subset A'$.
To do this, choose finitely many such base balls covering
$\overline{A''}$. For all sufficiently large $i$, the corresponding
local graphs cover $\supp S_{i,a}\cap(A''\times P^\perp)$,
because each graph describes all the support above its corresponding
ball, as shown in the preceding paragraph.
On overlaps the graph functions agree, since each describes the entire support above the overlap.
They therefore define one function $w_i:A''\to P^\perp$.
The projection identity in \eqref{eq:multiplicity-matching}
forces its current multiplicity to be one,  and orientation to be the graph orientation. Hence, for some fixed constant $C>0$ independent of $i$,
\begin{equation}\label{eq:extraction-graph-representation}
 S_{i,a}\llcorner(A''\times P^\perp)=\mathbf G_{w_i},
 \qquad
 \|w_i\|_{C^{1,\alpha}(A'')}\leq C,\qquad
 \|dw_i\|_{C^0(A'')}\longrightarrow0.
\end{equation}
The graph is stationary by
\cite[Chapter~7, Lemma~1.2, p.~212]{SimonGMT}.
It is therefore smooth and minimal by
\cite[Chapter~5, Theorem~6.1, pp.~150--151]{SimonGMT}.
Since the graph lies in $U_{i,a}$, which is open and disjoint
from the sets $U_{i,b}$ containing the other currents $S_{i,b}$
by Lemma~\ref{lem:matching}, each graph point has a neighborhood
in which $T_i=S_{i,a}$. Thus near each point of the graph, $T_i$ is exactly
this smooth multiplicity-one graph $\mathbf G_{w_i}$.
Hence, $\operatorname{graph}(w_i)$ is an open subset of
$\reg T_i$, so $w_i|_A$ is an exact sheet-height function
for $T_i$ in the sense of Definition~\ref{def:chart}.

\emph{Step 3: estimate the derivatives of $w_i$.}
We use interior Schauder estimates to show that the derivatives
of $w_i$ of every fixed order are $O(\varepsilon_i)$ on each
open ball $B\Subset A''$, with constants independent of $i$.
Equations \eqref{eq:multiplicity-matching} and
\eqref{eq:extraction-graph-representation}, together with the
boundedness of $\psi_a$ on $A''$, give
\begin{equation}\label{eq:extraction-height-estimate}
 \|w_i-\varepsilon_i\psi_a\|_{C^0(A'')}=o(\varepsilon_i),
 \qquad \|w_i\|_{C^0(A'')}\leq C\varepsilon_i.
\end{equation}
By \eqref{eq:induced-graph-metric}, the induced metric is
$G_i=\delta+(dw_i)^Tdw_i$. The minimal graph equation is
\begin{equation}\label{eq:minimal-graph-equation}
 \sum_{r,s=1}^k(G_i)^{rs}\partial_r\partial_sw_i^\mu=0,
 \qquad \mu=1,\ldots,m,
\end{equation}
where $((G_i)^{rs})$ is the inverse matrix of $G_i$ and
$\partial_r=\partial/\partial X^r$ in orthonormal coordinates on $P$.
Each component therefore satisfies a scalar linear equation
with the same principal coefficient matrix, once $w_i$ is fixed.

For large $i$, \eqref{eq:extraction-graph-representation} gives
$|dw_i|\leq1$. Thus, for every $\zeta\in\R^k$,
\[
 |\zeta|^2\leq\zeta^TG_i\zeta
 =|\zeta|^2+|dw_i\zeta|^2\leq2|\zeta|^2.
\]
The eigenvalues of $G_i$ therefore lie in $[1,2]$, so those
of $G_i^{-1}$ lie in $[1/2,1]$. Consequently,
\begin{equation}\label{eq:extraction-ellipticity}
 \frac{1}{2}|\zeta|^2
 \leq\sum_{r,s=1}^k(G_i)^{rs}\zeta_r\zeta_s
 \leq|\zeta|^2.
\end{equation}
The uniform bound $\|w_i\|_{C^{1,\alpha}(A'')}\leq C$ in
\eqref{eq:extraction-graph-representation} also gives a uniform
$C^\alpha(A'')$ bound for the coefficients $(G_i)^{rs}$, since matrix inversion
is smooth on these uniformly positive definite matrices.
Thus the interior Schauder estimate
\cite[Corollary~6.3]{GT}, applied componentwise, yields for each
open ball $B\Subset A''$
\begin{equation}\label{eq:extraction-schauder}
 \|w_i\|_{C^{2,\alpha}(B)}
 \leq C_B\|w_i\|_{C^0(A'')}
 \leq C_B\varepsilon_i,
\end{equation}
where $C_B$ may depend on $\dist(B,\partial A'')>0$ but is
independent of $i$.

Applying the interior Schauder estimate repeatedly to
\eqref{eq:minimal-graph-equation} and its differentiated equations
on nested open balls, as in the proof of
\cite[Theorem~6.17]{GT}, gives
\begin{equation}\label{eq:extraction-derivative-bounds}
 \|w_i\|_{C^b(B)}\leq C_{b,B}\varepsilon_i
 \quad\text{for every fixed integer }b\geq0
 \text{ and open ball }B\Subset A'',
\end{equation}
where $C_{b,B}$ is independent of $i$.

\emph{Step 4: prove $w_i/\varepsilon_i\to\psi_a$ in $C^b(A)$.}
Set $f_i:=w_i/\varepsilon_i$. Dividing
\eqref{eq:extraction-height-estimate} by $\varepsilon_i$ gives
\begin{equation}\label{eq:extraction-normalized-c0}
 \|f_i-\psi_a\|_{C^0(A'')}\longrightarrow0.
\end{equation}
Fix $b\geq0$ and choose an open ball $B$ with
$A\Subset B\Subset A''$. By
\eqref{eq:extraction-derivative-bounds},
\[
 \|f_i\|_{C^{b+1}(B)}\leq C_{b+1,B}.
\]
The Arzel\`a--Ascoli theorem, applied to the coordinate derivatives
through order $b$ on $\overline A$, gives a subsequence
converging in $C^b(A)$. Every such limit equals $\psi_a$
by \eqref{eq:extraction-normalized-c0}.
Since this argument applies to every subsequence,
the whole sequence converges to $\psi_a$ in $C^b(A)$.
Together with \eqref{eq:extraction-derivative-bounds}, this proves
\eqref{eq:normalized-sheet-convergence}.
\end{proof}

\subsection{The pointwise obstruction}
Under the hypotheses of Theorem~\ref{thm:intro-branch}(i), we now
prove that a nonzero cyclic component cannot have repetition
multiplicity one. Lemmas~\ref{lem:matching} and~\ref{lem:extract}
would turn such a component into exact sheet-height functions to
which Theorem~\ref{thm:local} applies.

\begin{proof}[Proof of Theorem~\ref{thm:intro-branch}(i)]
Fix $Z\in\mathcal R_q(T)$ and a cylindrical normalized tangent
function $\widehat\Phi$ at $Z$. Put
$\nu:=N_{T,\mathrm{Pl}}(Z)>1$, where the inequality follows from
\eqref{eq:planar-frequency-lower-bound}. Assume either that
$k=3$ and $\nu\ne2$, or that $k\geq4$.

\emph{Step 1: obtain exact sheet-height functions from a component
of repetition multiplicity one.}
Suppose that a nonzero cyclic component $F_j$ of
$F:=\widehat\Phi$, in the notation of
Section~\ref{sec:cyclic-components}, has $m_j=1$.
Choose one of its local branches $\varphi_{j,l}$ from
\eqref{eq:branches}. As in the setup of Lemma~\ref{lem:matching},
choose nested profile-separated annular sheet charts
$A\Subset A'\Subset A^+$ and an index $a$ such that
$\psi_a=\varphi_{j,l}$. The multiplicity of this local branch
in $\widehat\Phi|_{A^+}$ is $n_a=m_j=1$.

Choose a sequence of radii realizing $\widehat\Phi$, and let
$T_i$ and $\varepsilon_i$ be as in \eqref{eq:tangent-scales}.
Here $\varepsilon_i>0$ for all sufficiently large $i$, and
$\varepsilon_i\to0$ by \eqref{eq:cylindrical-height-decay}.
Lemma~\ref{lem:matching} gives the currents $S_{i,a}$ with
projection multiplicity $n_a=1$. Lemma~\ref{lem:extract} therefore
gives exact sheet-height functions $w_i:A\to P_Z^\perp$ for $T_i$.
Define $f_i:=\varepsilon_i^{-1}w_i$. By
\eqref{eq:normalized-sheet-convergence},
\begin{equation}\label{eq:pointwise-normalized-branch}
 f_i\longrightarrow\varphi_{j,l}\quad\text{in }C^b(A)
 \qquad\text{for every fixed integer }b\geq0.
\end{equation}

\emph{Step 2: apply the local graph obstruction.}
By \eqref{eq:branches}, using the angle in
Definition~\ref{def:chart}, we have
\[
 \varphi_{j,l}(z,t)=\operatorname{Re}(\widetilde c z^\nu),
 \qquad \widetilde c:=e^{2\pi i\nu l}c_j\ne0.
\]
Since $w_i=\varepsilon_i f_i$, the induced metric of
$\operatorname{graph}(w_i)$, expressed on $A$, by \eqref{eq:induced-graph-metric}, is
\begin{equation}\label{eq:pointwise-graph-metric}
 G_i=\delta+(dw_i)^Tdw_i
     =\delta+\varepsilon_i^2(df_i)^Tdf_i.
\end{equation}
We next show that these metrics are locally conformally flat.
The graph $\operatorname{graph}(w_i)$ is an open subset of $\reg T_i$ by
Lemma~\ref{lem:extract}. Under the map $W\mapsto Z+\rho_iW$, where $W\in\reg T_i\subset\R^{k+m}$,
the induced metric on $\reg T_i$ is $\rho_i^{-2}$ times the
pullback of the induced metric on $\reg T$. The latter is locally
conformally flat by hypothesis, and this property is preserved by
pullback under a diffeomorphism and multiplication by a positive
constant, as recalled in Section~\ref{sec:normalized-heights}.

On the other hand, \eqref{eq:pointwise-normalized-branch} gives the
$C^3(A)$ convergence required by Theorem~\ref{thm:local} when
$k=3$, and the $C^2(A)$ convergence required when $k\geq4$.
That theorem says that the metrics
\eqref{eq:pointwise-graph-metric} cannot be locally conformally
flat for all sufficiently large $i$ under our assumptions on
$k$ and $\nu$, which gives the desired contradiction.

\emph{Step 3: count the values of the normalized tangent function.}
Every repetition multiplicity is a positive integer. Excluding
$m_j=1$ for each nonzero cyclic component therefore gives
$m_j\geq2$. The decomposition in
Section~\ref{sec:cyclic-components} now yields
\begin{equation}\label{eq:pointwise-multiplicity-count}
 q=m_0+q_0\sum_{j=1}^Jm_j
       \geq m_0+2q_0J,
 \qquad\text{hence}\qquad 2q_0J\leq q-m_0.
\end{equation}
The choice of $\widehat\Phi$ was arbitrary, so this holds for every
cylindrical normalized tangent function at $Z$.
\end{proof}

\subsection{Stratification and the dimension conclusion}\label{sec:stratification}

We next use stratification to estimate the Hausdorff dimension of
the points where the counting inequality
\eqref{eq:pointwise-multiplicity-count} rules out every cylindrical
normalized tangent function. We then use the almost-everywhere
existence and uniqueness results to obtain the conclusion on the
full branch set.

\paragraph{The strata and the dimension estimate.}
We first state the Hausdorff dimension bound for the strata in
\eqref{eq:spine-strata} under assumptions~\ref{ass:stratification-cylinder}
and~\ref{ass:stratification-branch-decay} below, which are
\cite[Hypotheses~1.4--1.5, p.~14]{KWIII}, respectively,
expressed in our notation. In the paragraph
``Localization and a countable cover'' below, we verify them
for suitable localized portions of the original current after
translation, rotation, and rescaling.
Fix an integer $q\geq2$ and $\mu\in(0,1)$, write
$P_0:=\R^k\times\{0\}$, and fix
$\epsilon_0=\epsilon_0(k,m,q,\mu)\in(0,1)$ supplied by
\cite[Theorem~3.7]{KWIII}. Assume the following:
\begin{enumerate}[label=(\roman*),ref=(\roman*),leftmargin=2em]
\item\label{ass:stratification-cylinder} $T$ is a locally area-minimizing integral $k$-current in
$\mathbf C_1(0,P_0)=B_1^{P_0}(0)\times P_0^\perp$ and satisfies
\begin{equation}\label{eq:stratification-cylinder-hypotheses}
 \begin{gathered}
 (\partial T)\llcorner\mathbf C_1(0,P_0)=0,\qquad
 (\pi_{P_0})_\#T=q\llbracket B_1^{P_0}(0)\rrbracket,\\
 \sup_{W\in\supp T\cap\mathbf C_1(0,P_0)}
       \dist(W,P_0)<\epsilon_0.
 \end{gathered}
\end{equation}
\item\label{ass:stratification-branch-decay} For $\HH^{k-2}$-almost every branch point
$Z\in\sing T\cap\mathbf C_{1/2}(0,P_0)$ with
$\Theta_k(\|T\|,Z)=q$, the point $Z$ is a flat singular point
in the sense of Definition~\ref{def:density-branch}, with
tangent plane $P_Z$, and satisfies
\begin{equation}\label{eq:stratification-branch-decay-hypothesis}
 \sup_{0<\rho\leq1/8}\rho^{-k-2-2\mu}
 \int_{B_{2\rho}(Z)}\dist^2(W,Z+P_Z)\,\dd\|T\|(W)
 <\infty.
\end{equation}
\end{enumerate}
Define $\mathcal B_{q,\mu}$ to be the set of flat singular points
of density $q$ in $\mathbf C_{1/2}(0,P_0)$ satisfying
\eqref{eq:rapid} with the fixed exponent $\mu$.
For comparison, \cite[Definition~1.6, p.~14]{KWIII} defines
the set of such flat singular points satisfying
\eqref{eq:stratification-branch-decay-hypothesis}.
Since replacing $r$ by $2\rho$ in \eqref{eq:rapid} changes only
the constant $C$, and local finiteness of mass extends the bound
to all $0<\rho\leq1/8$, \eqref{eq:rapid} implies
\eqref{eq:stratification-branch-decay-hypothesis}.
Conversely, substituting $\rho=r/2$ in
\eqref{eq:stratification-branch-decay-hypothesis} gives
\eqref{eq:rapid} for sufficiently small $r>0$.
Thus $\mathcal B_{q,\mu}$ is precisely the set in their definition.

For $0\leq j\leq k-2$, define the \defterm{$j$th stratum} by
\begin{equation}\label{eq:spine-stratum-definition}
 \begin{split}
 S_{q,\mu}^{(j)}
 :=\{Z\in\mathcal B_{q,\mu}:{}&
       N_{T,\mathrm{Pl}}(Z)\ne2,\
       \dim\operatorname{spine}\widehat\Phi\leq j\\
     &\text{for every normalized tangent function }
       \widehat\Phi\text{ at }Z\}.
 \end{split}
\end{equation}
The equivalence between their definition of the spine using frequency
and our definition by translation invariance is established in
Section~\ref{sec:planar-frequency}.
By \cite[\S3.3, paragraph preceding Theorem~3.7, p.~31]{KWIII},
every normalized tangent function at a point of
$\mathcal B_{q,\mu}\cap\{N_{T,\mathrm{Pl}}\ne2\}$ has spine
dimension at most $k-2$.
Hence, combining this with the dimension estimate
\cite[Theorem~3.7(i), pp.~31--32]{KWIII}, we get:
\begin{equation}\label{eq:spine-strata}
 \mathcal B_{q,\mu}\cap\{N_{T,\mathrm{Pl}}\ne2\}
       =S_{q,\mu}^{(k-2)},\qquad
 \dim_{\HH}S_{q,\mu}^{(j)}\leq j.
\end{equation}

In particular, if $Z\in\mathcal B_{q,\mu}$ has planar frequency
different from $2$ and does not belong to $S_{q,\mu}^{(k-3)}$,
then \eqref{eq:spine-stratum-definition} gives at least one
normalized tangent function $\widehat\Phi$ at $Z$ with
$\dim\operatorname{spine}\widehat\Phi>k-3$.
The first equality in \eqref{eq:spine-strata} bounds that dimension
by $k-2$. Since it is an integer,
\begin{equation}\label{eq:stratum-cylindrical-tangent}
 \dim\operatorname{spine}\widehat\Phi=k-2.
\end{equation}
This is precisely the definition of a cylindrical normalized
tangent function in Section~\ref{sec:planar-frequency}.

\paragraph{Localization and a countable cover.}
We justify applying \eqref{eq:spine-strata} to all of
$\mathcal R_q(T)$. We first verify
assumptions~\ref{ass:stratification-cylinder}
and~\ref{ass:stratification-branch-decay} for localized currents
after translation, rotation, and rescaling.
Using these currents, we construct a countable cover of
$\mathcal R_q(T)$. The image of each set of this cover under the translation, rotation,
and rescaling used for the corresponding localized current
is contained in $\mathcal B_{q,\mu}$, so we can apply~\eqref{eq:spine-strata}
to it.

Choose a compact exhaustion $(K_\ell)_{\ell\geq1}$ of $\mathcal U$,
where each $K_\ell$ is compact,
$K_\ell\subset\operatorname{int}K_{\ell+1}$, and
$\bigcup_{\ell\geq1}K_\ell=\mathcal U$.
For each $Z\in\mathcal R_q(T)$, choose a member $K:=K_\ell$
of this exhaustion with $Z\in\operatorname{int}K$, and an open
ambient neighborhood $V_Z$ of $Z$ satisfying
$\overline{V_Z}\subset\operatorname{int}K$. Theorem~X(b)(ii)--(iv) in
\cite[pp.~7--8]{KWIII} supplies an exponent $\mu_K>0$ for which
$\HH^{k-2}$-almost every branch point in $K$ has a unique tangent
plane and satisfies \eqref{eq:rapid} with exponent $\mu_K$.

For the countable cover, we claim that every point
$Z\in\mathcal R_q(T)$ satisfies \eqref{eq:rapid} with exponent
$\mu=1/N$ for some integer $N\geq2$ with $\mu\leq\mu_K$.
By Definition~\ref{def:rapid}, choose some exponent $\mu_0>0$
for which \eqref{eq:rapid} holds at $Z$.
Choose $N$ with $1/N\leq\min\{\mu_0,\mu_K\}$ and decrease
$r_0$ so that $r_0\leq1$. Then
\begin{equation}\label{eq:rapid-exponent-comparison}
 r^{k+2+2\mu_0}\leq r^{k+2+2/N}\qquad(0<r<r_0),
\end{equation}
so \eqref{eq:rapid} holds at $Z$ with exponent $\mu$.
Applying \eqref{eq:rapid-exponent-comparison} with $\mu_0=\mu_K$
also shows that the preceding almost-everywhere statement holds
with the same exponent $\mu$.
The constant $C$ and the initial radius $r_0$ in \eqref{eq:rapid}
may depend on the point $Z$.

For each $Z\in\mathcal R_q(T)$ satisfying \eqref{eq:rapid}
with this fixed exponent $\mu=1/N$, apply
Lemma~\ref{lem:localization} at $Z$ to obtain a localized current.
Choose the radius $R>0$ in the proof of
Lemma~\ref{lem:localization} sufficiently small that
$\overline{B_{3R}(Z)}\subset V_Z$.
The lemma gives local area minimization and absence of boundary
in $\mathbf C_R(Z,P_Z)$.
For sufficiently small $0<r<R$, restrict the retained current
to $\mathbf C_r(Z,P_Z)$.
The support bound \eqref{eq:vanishing-relative-normal-height},
together with \eqref{eq:localized-cylinder-ball}, makes the
normal support arbitrarily small after rescaling
$\mathbf C_r(Z,P_Z)$ to the unit cylinder $\mathbf C_1(0,P_Z)$
using
\[
 W\longmapsto\frac{W-Z}{r}.
\]
Equation~\eqref{eq:localized-projection-identity} gives the
projection multiplicity $q$.
Finally, rotate $P_Z$ to $P_0$ and denote the transformed
localized current again by $T$.
The local area minimization and the boundary, projection,
and support bounds just established verify assumption~\ref{ass:stratification-cylinder}.

We next verify assumption~\ref{ass:stratification-branch-decay} for this
transformed current.
Before translation, rotation, and rescaling, the localized current
has a unique tangent plane and satisfies \eqref{eq:rapid}
with exponent $\mu$ at $\HH^{k-2}$-almost every branch point,
because its retained support lies in $K$ and it agrees locally
with the original current at these points.
Translation, rotation, and dilation preserve tangent-cone
uniqueness, the exponent $\mu$, and Hausdorff-null sets.
By the same argument used above to establish the equivalence of
\eqref{eq:rapid} and \eqref{eq:stratification-branch-decay-hypothesis},
local finiteness of mass bounds the expression under the supremum in
\eqref{eq:stratification-branch-decay-hypothesis} at the remaining
radii up to $\rho=1/8$.
Together with the decay estimate at small radii, this shows that
the supremum is finite, verifying
assumption~\ref{ass:stratification-branch-decay}.

We now construct a countable cover of $\mathcal R_q(T)$.
By Lemma~\ref{lem:localization}, the retained current agrees
with the original current on some open neighborhood
$U_Z\subset V_Z$ of $Z$ in the ambient space.
Since $Z$ is mapped to the origin, choose $U_Z$ small enough
that its image under the translation, rotation, and rescaling
is contained in the inner half-cylinder $\mathbf C_{1/2}(0,P_0)$.
For each member $K_\ell$ of the exhaustion and each integer
$N\geq2$ with $1/N\leq\mu_{K_\ell}$, take a countable subcover
of these neighborhoods covering the points of
$\mathcal R_q(T)\cap\operatorname{int}K_\ell$ that satisfy
\eqref{eq:rapid} with exponent $1/N$.
Intersect the neighborhoods with that set of points.
Since every point of $\mathcal R_q(T)$ lies in some
$\operatorname{int}K_\ell$ and satisfies \eqref{eq:rapid} with
exponent $1/N$ for all sufficiently large $N$, this gives a
countable collection of subsets $\mathcal F_h$ with
\begin{equation}\label{eq:rapid-countable-cover}
 \mathcal R_q(T)=\bigcup_{h=1}^{\infty}\mathcal F_h,
\end{equation}
each of which is mapped into a set $\mathcal B_{q,\mu}$ as above.
Translations, rotations, and dilations preserve planar frequency,
spine dimension, repetition multiplicities, and Hausdorff dimension.
Thus \eqref{eq:spine-strata} can be applied on every
$\mathcal F_h$.

\begin{proof}[Proof of Theorem~\ref{thm:intro-branch}(ii)]
Recall that
\[
 \mathcal E_q=
 \{Z\in\mathcal R_q(T):\nu:=N_{T,\mathrm{Pl}}(Z)\notin\mathbb Q
 \ \text{or}\;
 \nu=p_0/q_0\notin\mathbb Z,\ \gcd(p_0,q_0)=1,\ 2q_0>q\}.
\]
Fix a set $\mathcal F_h$ from \eqref{eq:rapid-countable-cover},
and let $\Psi_h$ be the composition of the translation, rotation,
and rescaling associated with this set, so that
$\Psi_h(\mathcal F_h)\subset\mathcal B_{q,\mu}$.
Let $Z_0\in\mathcal E_q\cap\mathcal F_h$ and set $Z:=\Psi_h(Z_0)$.
Since $\Psi_h$ preserves planar frequency,
the planar frequency at $Z$ equals $\nu:=N_{T,\mathrm{Pl}}(Z_0) \neq 2$.

Suppose for contradiction that $Z\notin S_{q,\mu}^{(k-3)}$, as in definition~\eqref{eq:spine-stratum-definition}.
By \eqref{eq:stratum-cylindrical-tangent}, there is a cylindrical
normalized tangent function at $Z$. The classification in
Section~\ref{sec:cyclic-components} makes $\nu$ a rational
number $p_0/q_0$ in lowest terms. Since $Z_0 \in \mathcal E_q$, we have $2q_0>q$.
Since $m_0\geq0$ and $J\geq1$, Theorem~\ref{thm:intro-branch}(i),
specifically \eqref{eq:pointwise-multiplicity-count}, gives
\[
 q\geq m_0+2q_0J\geq2q_0,
\]
which is a contradiction. So, we have proved
\begin{equation}\label{eq:exceptional-set-stratum}
 \Psi_h(\mathcal E_q\cap\mathcal F_h)
       \subset S_{q,\mu}^{(k-3)}.
\end{equation}

By \eqref{eq:exceptional-set-stratum},
\eqref{eq:spine-strata}, and preservation of Hausdorff dimension
under $\Psi_h$, each set $\mathcal E_q\cap\mathcal F_h$
has Hausdorff dimension at most $k-3$.
Hence,~\eqref{eq:rapid-countable-cover} gives
\[
 \dim_{\HH}\mathcal E_q
 =\dim_{\HH}\bigcup_{h=1}^{\infty}
       (\mathcal E_q\cap\mathcal F_h)
 \leq k-3.
\]
Finally, for $q=2$ or $q=3$, note that any nonintegral rational $p_0/q_0$ in
lowest terms has $q_0 \geq 2$ and hence $2q_0\geq4>q$. Thus, every
point of $\mathcal R_2(T)\cup\mathcal R_3(T)$ with $\nu \not\in \mathbb{Z}$ belongs to
$\mathcal E_2\cup\mathcal E_3$, which gives the required
\[
 \dim_{\HH}\{Z\in\mathcal R_2(T)\cup\mathcal R_3(T):
                    N_{T,\mathrm{Pl}}(Z)\notin\mathbb Z\}\leq k-3.
\]
\end{proof}

\begin{proof}[Proof of Theorem~\ref{thm:intro-branch}(iii)]
By \cite[Theorem~X(b)(ii)--(iv), pp.~7--8]{KWIII},
the branch points of density $q$ outside $\mathcal R_q(T)$
have zero $\HH^{k-2}$ measure.
The application of \cite[Theorem~1.8(ii)--(iii)]{KWIII}
to each localized current is justified in
Section~\ref{sec:planar-frequency}.
Since the cover \eqref{eq:rapid-countable-cover} is countable,
we conclude that at $\HH^{k-2}$-almost every point of
$\mathcal R_q(T)$ with planar frequency different from $2$,
there is a unique cylindrical normalized tangent function.

At each such point, Theorem~\ref{thm:intro-branch}(i) applies to
this normalized tangent function. Every nonzero cyclic component
therefore has repetition multiplicity $m_j\geq2$, and
\eqref{eq:pointwise-multiplicity-count} gives
$2q_0J\leq q-m_0$, as claimed.

For $q=2,3$, the nonintegral-frequency branch points split into
those in $\mathcal R_q(T)$ and those outside it.
Part~(ii) bounds the Hausdorff dimension of the former set by
$k-3$. Since $k-3<k-2$, this set has zero $\HH^{k-2}$ measure.
The latter set has zero $\HH^{k-2}$ measure by above. Hence, the nonintegral-frequency
points of densities $2$ and $3$ have zero $\HH^{k-2}$ measure.
\end{proof}

\begin{remark}
When $k=3$, the full branch set has Hausdorff dimension at most $1$
by Almgren's theorem \cite[Chapter~7, \S4, p.~223]{SimonGMT}.
Under the hypotheses of Theorem~\ref{thm:intro-branch}, part~(iii)
shows that the branch points of densities $2$ and $3$ with nonintegral
planar frequency have zero $\HH^1$ measure. Part~(ii) shows that the subset of these points in
$\mathcal R_2(T)\cup\mathcal R_3(T)$ is either empty or has
Hausdorff dimension zero. Note that neither a set with Hausdorff dimension zero nor an
$\HH^{k-2}$-null set need be empty or even discrete.
\end{remark}

\section{Flat critical points and unique continuation}\label{sec:flatness}
We first show that weak conformality rules out nonzero homogeneous polynomial maps
of degree greater than one when the domain dimension is at least three.
We use this result to obtain Taylor flatness at critical points. We then prove the continuation
criterion under the curvature bound and examine its hypotheses
using the Bor\r{u}vka cone example.

\subsection{Homogeneous-polynomial rigidity}

As shown in Proposition~\ref{prop:polynomial-equivalence} below,
the following Theorem~\ref{thm:homogeneous-polynomial} is equivalent to the result of
P\'erez Fern\'andez--P\'erez \cite[Theorem~4.1]{PFP}.
Their statement does not assume homogeneity.
Their proof uses restrictions to affine two-dimensional planes and a Hessian argument.
We give a proof using spherical harmonics, as in the conventions of Section~\ref{sec:spherical-definitions}.

\begin{theorem}[Homogeneous weakly-conformal polynomials]
\label{thm:homogeneous-polynomial}
Let $n\geq k\geq3$, and let $P:\R^k\to\R^n$ be a nonzero homogeneous
polynomial of degree $s\geq1$.  If $P$ is weakly-conformal, meaning that
for some function $\Lambda\geq0$ we have
\[
 (dP_x)^TdP_x=\Lambda(x)^2I_k
 \qquad\text{for every }x\in\R^k,
\]
then $s=1$.
\end{theorem}

\begin{proof}
Set $d_{\rm sph}:=k-1$, write $x=r\theta$ for
$r:=|x|>0$ and $\theta:=x/|x|\in S^{d_{\rm sph}}$, and define
\[
 \psi:=P|_{S^{d_{\rm sph}}}:S^{d_{\rm sph}}\to\R^n.
\]
Write $d_S\psi$ for its differential along the sphere and
$|d_S\psi|$ for its Hilbert--Schmidt norm. All integrals in this proof use the volume measure
induced by the round metric on $S^{d_{\rm sph}}$.
We compute $\int_{S^{d_{\rm sph}}}|d_S\psi|^2$ using weak conformality
and bound the same integral using spherical harmonics.

By homogeneity, $P(r\theta)=r^s\psi(\theta)$, so differentiation gives,
for $X\in T_\theta S^{d_{\rm sph}}$ viewed as a Euclidean vector at $r\theta$,
\[
 dP_{r\theta}(\theta)=sr^{s-1}\psi(\theta),
 \qquad
 dP_{r\theta}(X)=r^{s-1}d_S\psi_\theta(X).
\]
Weak conformality makes the images under $dP_x$ of all unit vectors
have the same length at each $x$. At $x=\theta$, this common length
is $s|\psi(\theta)|$, since $|\theta|=1$ and
$dP_\theta(\theta)=s\psi(\theta)$.
Thus, $|d_S\psi|^2=d_{\rm sph}s^2|\psi|^2$,
and integrating gives
\begin{equation}\label{eq:spherical-exact}
 \int_{S^{d_{\rm sph}}}|d_S\psi|^2
 =d_{\rm sph}s^2\int_{S^{d_{\rm sph}}}|\psi|^2.
\end{equation}

By the Fischer decomposition \cite[Theorem~1.3, p.~24]{Gallier}, each
homogeneous scalar polynomial $Q$ of degree $s$ on
$\R^{d_{\rm sph}+1}$ has a unique decomposition
\[
 Q=H_s+|x|^2H_{s-2}+|x|^4H_{s-4}+\cdots
\]
into homogeneous harmonic polynomials $H_h$, where each $H_h$ has degree $h$.
Consequently, each component $\psi^\mu$, $1\leq\mu\leq n$,
is a sum of spherical harmonics of degrees $s,s-2,\ldots$.

For a spherical harmonic $f_h$ of degree $h$, write
$\lambda_h:=h(h+d_{\rm sph}-1)$.
Equation~\eqref{eq:spherical-harmonic-eigenvalue} gives
$-\Delta_{S^{d_{\rm sph}}}f_h=\lambda_h f_h$.
For spherical harmonics $f_h,f_j$ of distinct degrees, integration
by parts shows that $f_h$ and $f_j$ are $L^2$-orthogonal, as are
their differentials $d_Sf_h$ and $d_Sf_j$.
Integration by parts also gives
\[
 \int_{S^{d_{\rm sph}}}|d_Sf_h|^2
 =\lambda_h\int_{S^{d_{\rm sph}}}|f_h|^2.
\]
Every spherical harmonic in the decomposition of $\psi^\mu$ has degree $h\leq s$, so
$\lambda_h\leq s(s+d_{\rm sph}-1)$.
Using the $L^2$-orthogonality of the spherical harmonics in the decomposition
of $\psi^\mu$ and of their differentials, we obtain
\[
 \int_{S^{d_{\rm sph}}}|d_S\psi^\mu|^2
 \leq s(s+d_{\rm sph}-1)
       \int_{S^{d_{\rm sph}}}|\psi^\mu|^2.
\]
Summing over $\mu$ gives
\begin{equation}\label{eq:spherical-upper}
 \int_{S^{d_{\rm sph}}}|d_S\psi|^2
 \leq s(s+d_{\rm sph}-1)\int_{S^{d_{\rm sph}}}|\psi|^2.
\end{equation}
Since $P\ne0$ and is homogeneous, its restriction $\psi$ is not
identically zero, so $\int_{S^{d_{\rm sph}}}|\psi|^2>0$.
Combining \eqref{eq:spherical-exact} and \eqref{eq:spherical-upper}
and using $\int_{S^{d_{\rm sph}}}|\psi|^2>0$ and $s>0$, we obtain
\[
 d_{\rm sph}s^2\leq s(s+d_{\rm sph}-1).
\]
Therefore,
\[
 (d_{\rm sph}-1)(s-1)\leq0.
\]
Because $d_{\rm sph}=k-1\geq2$ and $s\geq1$, this forces $s=1$.
\end{proof}
\begin{proposition}
\label{prop:polynomial-equivalence}
For $n\geq k\geq3$, Theorem~\ref{thm:homogeneous-polynomial} is equivalent
to \cite[Theorem~4.1]{PFP}, which states that every weakly-conformal
polynomial $u:\R^k\to\R^n$ is affine.
\end{proposition}

\begin{proof}
Since \cite[Theorem~4.1]{PFP} does not assume homogeneity, it implies
Theorem~\ref{thm:homogeneous-polynomial}.

Conversely, suppose that a weakly-conformal polynomial
$u:\R^k\to\R^n$ has degree $s\geq2$.
For $0\leq j\leq s$, let $P_j$ be its homogeneous part of degree $j$,
so that $u=P_s+\cdots+P_0$ and $P_s\ne0$.
Weak conformality means that
\begin{equation}\label{eq:polynomial-conformality}
 \langle\partial_i u,\partial_j u\rangle
 =\frac{1}{k}|du|^2\delta_{ij}.
\end{equation}
Since $\partial_iP_s$ has degree $s-1$ and all derivatives of
$P_{s-1}+\cdots+P_0$ have degree at most $s-2$, every product
involving a derivative of a lower-degree part has degree at most
$2s-3$. Equating the homogeneous terms of degree $2s-2$ in
\eqref{eq:polynomial-conformality} therefore gives
\[
 \langle\partial_iP_s,\partial_jP_s\rangle
 =\frac{1}{k}|dP_s|^2\delta_{ij}.
\]
Thus $P_s$ is a nonzero homogeneous weakly-conformal polynomial
of degree $s\geq2$, contradicting Theorem~\ref{thm:homogeneous-polynomial}.
Thus, every weakly-conformal polynomial has degree at most one.
\end{proof}

\begin{remark}
The restriction $k\geq3$ is sharp for this statement.  When $k=2$, the
factor $k-2$ in the final inequality is zero, and the polynomials
$z\mapsto z^s$, with $z=x+iy$ identifying $\R^2$ with $\C$, give nonzero
homogeneous conformal polynomials of every degree $s\geq1$.
\end{remark}

\subsection{Proof of critical-point flatness}
We apply Theorem~\ref{thm:homogeneous-polynomial} to the first nonzero
homogeneous term in the Taylor expansion at a critical point of a weakly-conformal map $u$.
We then bound the derivatives of $u$ in terms of the distance to the
critical set and use these bounds to exclude critical points of
nonconstant smooth maps under the curvature assumption of
Theorem~\ref{thm:intro-continuation}.

\begin{flatnessrestatement}[Finite-order flatness]
Let $u:(L^k,g)\to(M^n,h)$ be $C^j$, where $j\geq2$, $n\geq k\geq3$,
and the metrics are smooth.  If $u$ is {\color{black}weakly-conformal}, then
for every $p\in L$ with $du(p)=0$, we have
\[
 D^bu(p)=0\quad(2\leq b\leq j)
\]
in any choice of local coordinates on the domain and target.
In local domain coordinates centered at $p$, as $x\to0$,
\[
 \begin{aligned}
 u(x)-u(p)&=o(|x|^j),\\
 D^bu(x)&=o(|x|^{j-b})\qquad(1\leq b\leq j).
 \end{aligned}
\]
If, in addition, $u\in C^{j,\alpha}$, $0<\alpha\leq1$, then
\[
 \begin{aligned}
 u(x)-u(p)&=O(|x|^{j+\alpha}),\\
 D^bu(x)&=O(|x|^{j+\alpha-b})\qquad(1\leq b\leq j).
 \end{aligned}
\]
In particular, every critical point of a smooth {\color{black}weakly-conformal} map is Taylor-flat.
\end{flatnessrestatement}

\begin{proof}
Choose normal coordinates on the domain and target centered at $p$
and $u(p)$, respectively.
Write the coordinate representation again as $u$, with $u(0)=0$.
Let $\partial_i$ denote differentiation in the $i$th domain
coordinate, let $g_{ij}$ be the domain metric coefficients, and use
$|du|_{g,h}$ for the Hilbert--Schmidt norm.

Let $s\in\{2,\ldots,j\}$ be the first order at which a derivative of $u$
at $0$ is nonzero, if such an order exists.
Define the nonzero homogeneous polynomial $P_s:\R^k\to\R^n$
componentwise by
\[
 P_s^\mu(x):=\frac{1}{s!}
 \sum_{i_1,\ldots,i_s=1}^k
 (\partial_{i_1}\cdots\partial_{i_s}u^\mu)(0)
 x^{i_1}\cdots x^{i_s},
 \qquad 1\leq\mu\leq n.
\]
For each fixed $x\in\R^k$, Taylor's theorem applied to $u$ and $du$ gives
\[
 u(tx)=t^sP_s(x)+o(t^s),\qquad
 du(tx)=t^{s-1}dP_s(x)+o(t^{s-1})
 \quad(t\downarrow0).
\]
Weak conformality of $u$ means
\begin{equation}\label{eq:flatness-rescaled-conformality}
 h_{u(tx)}\bigl(\partial_i u(tx),\partial_j u(tx)\bigr)
 =\frac{1}{k}|du(tx)|_{g,h}^2g_{ij}(tx).
\end{equation}
As $t\downarrow0$,
\[
 g_{ij}(tx)\longrightarrow\delta_{ij},\qquad
 h_{\alpha\beta}(u(tx))\longrightarrow\delta_{\alpha\beta}.
\]
Dividing \eqref{eq:flatness-rescaled-conformality} by $t^{2s-2}$
and using these convergences and the Taylor expansions gives
\[
 \langle\partial_iP_s(x),\partial_jP_s(x)\rangle
 =\frac{1}{k}|dP_s(x)|^2\delta_{ij}.
\]
The polynomial $P_s$ is therefore weakly-conformal, so
Theorem~\ref{thm:homogeneous-polynomial} forces $s=1$,
contradicting $s\geq2$. Thus, every positive-order coordinate derivative
through order $j$ vanishes at $p$. By the chain rule, this vanishing holds in every choice of local coordinates.
Taylor's theorem applied to $u$ and its coordinate derivatives gives
the stated small-$o$ estimates. If, in addition, $u\in C^{j,\alpha}$,
the H\"older continuity of the derivatives of order $j$ gives the
stated big-$O$ estimates.
\end{proof}
\begin{corollary}[Uniform flatness near the zero set]\label{cor:uniform}
Let $u:(L^k,g)\to(M^n,h)$ be a $C^j$ weakly-conformal map, where
$j\geq2$, $n\geq k\geq3$, and $g,h$ are smooth Riemannian metrics.
Assume $Z_u:=\{p\in L:du(p)=0\}\ne\varnothing$, and set
$r_u(x):=\dist_g(x,Z_u)$. In fixed domain and target local coordinates,
\[
 |D^bu(x)|=o\bigl(r_u(x)^{j-b}\bigr)
 \qquad(1\leq b\leq j)
\]
as $r_u(x)\downarrow0$. If, in addition, $u\in C^{j,\alpha}$,
$0<\alpha\leq1$, then
\[
 |D^bu(x)|=O\bigl(r_u(x)^{j+\alpha-b}\bigr)
 \qquad(1\leq b\leq j)
\]
as $r_u(x)\downarrow0$.
For smooth $u$, for every integer $b\geq1$,
\begin{equation}\label{eq:uniform-smooth-flatness}
 |du(x)|\leq C_b\,r_u(x)^b.
\end{equation}
All estimates are uniform on each fixed compact subset of the domain.
\end{corollary}
\begin{proof}
Fix $p\in L$, and let $U$ and $W$ be local coordinate neighborhoods
of $p$ and $u(p)$ in the domain and target, respectively.
Identify points of $U$ with their coordinate vectors. After shrinking
$U$, we may assume that its coordinate image is convex and, by continuity
and positive definiteness of $g$, that there exists a constant $A\geq1$
such that
\begin{equation}\label{eq:uniform-coordinate-distance}
 |x-z|\leq A\,\dist_g(x,z)\qquad(x,z\in U).
\end{equation}
Choose $K$ compact and $V$ open so that
\[
 K\subset V\Subset U\subset L,
 \qquad u(\overline V)\subset W\subset M.
\]
For $z,x\in U$, write
\[
 [z,x]:=\{z+t(x-z):0\leq t\leq1\}\subset U.
\]
By the definition of $r_u$, compactness of $K\subset V$, and
\eqref{eq:uniform-coordinate-distance}, there exists $\delta>0$
such that, for every $x\in K$ with $0<r_u(x)<\delta$, there is
$z\in Z_u\cap V$ satisfying
\begin{equation}\label{eq:uniform-critical-point-choice}
 [z,x]\subset V,\qquad |x-z|\leq 2A r_u(x).
\end{equation}
Using $[z,x]\subset V$ from \eqref{eq:uniform-critical-point-choice}
and uniform continuity of $D^ju$ on $\overline V$, we obtain from
Theorem~\ref{thm:intro-flatness}
\[
 |D^bu(x)|=o\bigl(|x-z|^{j-b}\bigr)
 \qquad(1\leq b\leq j),
\]
uniformly in $x,z$ as $|x-z|\to0$. Combining this with the distance bound in
\eqref{eq:uniform-critical-point-choice} proves the uniform small-$o$ estimates.
If $u\in C^{j,\alpha}$, H\"older continuity of $D^ju$ on $\overline V$
gives
\[
 |D^bu(x)|=O\bigl(|x-z|^{j+\alpha-b}\bigr)
 \qquad(1\leq b\leq j),
\]
uniformly in $x,z$, and the same distance bound gives the stated
big-$O$ estimates.
For smooth $u$, for every integer $N\geq1$, taking $j=N+1$ gives
 a constant $C_N>0$ such that
\[
 |du(x)|\leq C_Nr_u(x)^N.
\]
A finite cover with local coordinate neighborhoods gives uniform estimates on every fixed
compact subset.
\end{proof}
\subsection{Continuation under an inverse-distance curvature bound}
We use the Yamabe equation and Theorem 1.1 to exclude critical points of a nonconstant smooth weakly-conformal map under the stated curvature bound. Recall that $\lambda$ is the dilation of $u$ and $\gamma =u^{*}h = \lambda^2 g$ is the pullback metric on $L\setminus Z_u$, as defined in Section~\ref{sec:energy}.

\begin{continuationrestatement}[Continuation under an inverse-distance curvature bound]
Let $L$ be connected and let $u:(L^k,g)\to(M^n,h)$ be smooth and weakly-conformal, $k\geq3$. Assume $Z_u\ne\varnothing$ and that for every
compact $K\subset L$ there is a constant $C_K\geq0$ such that
\[
 \lambda(x)^2\operatorname{Scal}_\gamma(x)
       \geq-\frac{C_K}{r_u(x)^2}
       \qquad(x\in K\setminus Z_u).
\]
Then $u$ is constant.
  It suffices that
$r_u\lambda|\mathrm{II}|_\gamma$  is bounded on $K\setminus Z_u$
for every compact $K\subset L$.
\end{continuationrestatement}

\begin{proof}
Set
\[
 c_k:=\frac{4(k-1)}{k-2},\qquad
 w:=\lambda^{(k-2)/2}\quad\text{on }L,\qquad
 V:=\frac{\operatorname{Scal}_g-
                \lambda^2\operatorname{Scal}_\gamma}{c_k}
       \quad\text{on }L\setminus Z_u.
\]
The function $w$ is continuous on $L$ and smooth and positive on
$L\setminus Z_u$.
The Yamabe equation
\cite[equation~(11.4.3), p.~331]{BW}, applied to
$\gamma=\lambda^2g$ on $L\setminus Z_u$, gives
\begin{equation}\label{eq:yamabe}
 \Delta_g w=Vw.
\end{equation}

Suppose, for a contradiction, that $u$ is nonconstant. Then
$L\setminus Z_u\ne\varnothing$. Since $Z_u$ is nonempty and closed
and $L$ is connected, choose
$p\in Z_u\cap\overline{L\setminus Z_u}$.
Choose a compact set $K\subset L$ with $p\in\operatorname{int}K$.
The functions $\operatorname{Scal}_g$ and $r_u$ are bounded on $K$, so the curvature
assumption gives a constant
$A>0$ such that
\begin{equation}\label{eq:continuation-potential-bound}
 V(x)\leq A/r_u(x)^2\qquad(x\in K\setminus Z_u).
\end{equation}

  We next choose a ball in $L\setminus Z_u$ whose boundary intersects $Z_u$.
For a point $a \in L$, denote by $\operatorname{inj}_g(a)$ its injectivity radius. Choose a relatively compact neighborhood $U$ of $p$.
Then there exists $\iota>0$ such that
\[
 \operatorname{inj}_g(y)\geq\iota\qquad(y\in\overline U).
\]
Since $p\in U\cap\operatorname{int}K$, choose $\varepsilon\in(0,\iota)$
so small that
\[
 \overline{B_{3\varepsilon}^g(p)}
 \subset U\cap\operatorname{int}K.
\]
The set $\overline{B_{3\varepsilon}^g(p)}$ is compact because it is closed and contained in the compact set $K$.
As $p\in\overline{L\setminus Z_u}$, choose
$a\in B_\varepsilon^g(p)\setminus Z_u$ and set $R:=r_u(a)$.
Closedness of $Z_u$ and $p\in Z_u$ give
\[
 0<R\leq\dist_g(a,p)<\varepsilon<\iota
 \leq\operatorname{inj}_g(a).
\]
For $x\in\overline{B_R^g(a)}$, the triangle inequality gives
\[
 \dist_g(x,p)\leq\dist_g(x,a)+\dist_g(a,p)
 \leq R+\dist_g(a,p)<2\varepsilon.
\]
Thus,
\[
 \overline{B_R^g(a)}\subset B_{2\varepsilon}^g(p)
 \subset\operatorname{int}K.
\]
By the definition of $R$, for each integer $j\geq1$ choose $q_j\in Z_u$
such that
\[
 R\leq\dist_g(a,q_j)<R+\frac{\varepsilon}{j}.
\]
Then
\[
 \dist_g(p,q_j)
 \leq\dist_g(p,a)+\dist_g(a,q_j)
 <\varepsilon+R+\frac{\varepsilon}{j}<3\varepsilon.
\]
Thus $(q_j)$ lies in the compact set $\overline{B_{3\varepsilon}^g(p)}$.
So, it has a subsequence which converges to some point $q \in \overline{B_{3\varepsilon}^g(p)}$. Since $Z_u$ is closed
and the distance function is continuous,
\[
 q\in Z_u,\qquad
 \dist_g(a,q)=\lim_{\ell\to\infty}\dist_g(a,q_{j_\ell})=R.
\]
Since $R<\operatorname{inj}_g(a)$, this gives
$q\in Z_u\cap\partial B_R^g(a)$.

By the definition of $R$, the open ball $B_R^g(a)$ contains no
critical points. Define
\[
 \delta:\overline{B_R^g(a)}\to[0,R],\qquad
 \delta(x):=R-\dist_g(a,x).
\]
For $z\in Z_u$, the triangle inequality gives
$\dist_g(x,z)\geq R-\dist_g(a,x)=\delta(x)$, and taking the infimum
over $z$ yields $r_u(x)\geq\delta(x)$.
Since $R<\operatorname{inj}_g(a)$, the distance from $a$ is smooth
on a neighborhood of $\partial B_R^g(a)$, so
$|\nabla\delta|_g=1$ and $\Delta_g\delta$ is bounded there.

Choose $\beta>1$ with $\beta(\beta-1)>A$, and define the operator
$\mathcal L:=\Delta_g-A\delta^{-2}$ on $B_R^g(a)\setminus\{a\}$.
By \eqref{eq:yamabe}, \eqref{eq:continuation-potential-bound},
$r_u\geq\delta$, and the preceding bounds on $\delta$,
there is $\eta\in(0,R)$ such that
\begin{equation}\label{eq:continuation-barrier}
 \mathcal Lw\leq0,\qquad
 \mathcal L(\delta^\beta)
 =\delta^{\beta-2}
   \bigl(\beta(\beta-1)-A+\beta\delta\Delta_g\delta\bigr)>0
 \quad(0<\delta<\eta).
\end{equation}
Since $w>0$ on the compact set
$\{x\in B_R^g(a):\delta(x)=\eta\}$,
choose $b>0$ so that $w\geq b\delta^\beta$ there.
On $\partial B_R^g(a)$, we have $w\geq0=b\delta^\beta$.
We claim that
\begin{equation}\label{eq:continuation-lower-bound}
 w\geq b\delta^\beta\qquad(0<\delta<\eta).
\end{equation}
Set
\[
 E:=\{x\in L:R-\eta\leq\dist_g(a,x)\leq R\}.
\]
The set $E$ is compact because it is closed and contained in
$\overline{B_R^g(a)}\subset K$.
Define $F:=w-b\delta^\beta$ on $\overline{B_R^g(a)}$.
If \eqref{eq:continuation-lower-bound} failed, continuity of $F$ and
compactness of $E$ would give $x_*\in E$ such that
\[
 F(x_*)=\min_E F<0.
\]
For $x\in E$ with $\delta(x)=0$ or $\delta(x)=\eta$, we have $F(x)\geq0$.
Thus $0<\delta(x_*)<\eta$. Hence $F$ is smooth near $x_*$, so
$\Delta_gF(x_*)\geq0$ and, applying $\mathcal L$ to the restriction of $F$
to $B_R^g(a)\setminus\{a\}$, we get
\[
 (\mathcal LF)(x_*)
 =\Delta_gF(x_*)-A\delta(x_*)^{-2}F(x_*)>0.
\]
But \eqref{eq:continuation-barrier} gives
$\mathcal LF=\mathcal Lw-b\mathcal L(\delta^\beta)<0$
where $0<\delta<\eta$, a contradiction.

Recall that $q\in Z_u\cap\partial B_R^g(a)$ and $\dist_g(a,q)=R$.
Since $R<\operatorname{inj}_g(a)$, choose the minimizing geodesic joining
$a$ to $q$ and parametrize it by arc length starting at $q$:
\[
 \sigma:[0,R]\to L,\qquad
 \sigma(0)=q,\quad \sigma(R)=a,\quad
 |\dot\sigma(t)|_g=1\qquad(0\leq t\leq R).
\]
For $0\leq s<t\leq R$, the restricted curve $\sigma|_{[s,t]}$ is minimizing.
Thus
\[
 \dist_g(q,\sigma(t))=t,\qquad
 \dist_g(a,\sigma(t))=R-t\qquad(0\leq t\leq R).
\]
Hence $\delta(\sigma(t))=t$. For $0<t<\eta$,
\eqref{eq:continuation-lower-bound} therefore gives
\[
 w(\sigma(t))\geq bt^\beta.
\]
Choose an integer $N\geq1$ such that $N(k-2)/2>\beta$.
Since $q\in Z_u$,
\[
 r_u(\sigma(t))\leq\dist_g(\sigma(t),q)=t
 \qquad(0\leq t\leq R).
\]
The estimate \eqref{eq:uniform-smooth-flatness} in
Corollary~\ref{cor:uniform}, together with
$\lambda=|du|_{g,h}/\sqrt{k}$, gives constants $C_N>0$ and
$t_N\in(0,\eta)$ such that
\[
 \lambda(\sigma(t))\leq C_Nr_u(\sigma(t))^N
 \leq C_Nt^N\qquad(0<t<t_N).
\]
Using $w=\lambda^{(k-2)/2}$, we obtain
\[
 0<b\leq\frac{w(\sigma(t))}{t^\beta}
 \leq C_N^{(k-2)/2}t^{N(k-2)/2-\beta}
 \qquad(0<t<t_N).
\]
Letting $t\downarrow0$ contradicts $b>0$, proving that $u$ is constant.

We now prove that boundedness of $r_u\lambda|\mathrm{II}|_\gamma$
on $K\setminus Z_u$ for every compact $K\subset L$ implies the
scalar-curvature bound in the hypothesis.

At a point $p\in L\setminus Z_u$, choose a $\gamma$-orthonormal basis
$e_1,\ldots,e_k$ of $T_pL$ and set $E_i:=du_p(e_i)$ for $1\leq i\leq k$.
Since $\gamma=u^*h$, the vectors $E_1,\ldots,E_k$ form an
$h$-orthonormal basis of $du_p(T_pL)\subset T_{u(p)}M$. Define
\[
 \begin{aligned}
 \operatorname{Scal}_{M,\mathrm{tan}}
 &:=\sum_{i,j=1}^k h\bigl(R^h(E_i,E_j)E_j,E_i\bigr),\\
 \tr_\gamma\mathrm{II}
 &:=\sum_{i=1}^k\mathrm{II}(E_i,E_i).
 \end{aligned}
\]
Both sums are independent of the chosen orthonormal basis.
 Contracting the Gauss equation \eqref{eq:gauss-equation}  gives
 \begin{equation}\label{eq:contracted-gauss}
 \begin{aligned}
 \operatorname{Scal}_\gamma
 &=\operatorname{Scal}_{M,\mathrm{tan}}
   +|\tr_\gamma\mathrm{II}|_h^2-|\mathrm{II}|_\gamma^2\\
 &\geq\operatorname{Scal}_{M,\mathrm{tan}}
                         -|\mathrm{II}|_\gamma^2.
 \end{aligned}
\end{equation}
As an aside, for a minimal immersion, including a smooth weakly-conformal
$k$-harmonic immersion by \cite[Corollary~4, pp.~221--222]{Takeuchi},
$\tr_\gamma\mathrm{II}=0$, giving
 \begin{equation}\label{eq:minimal-contracted-gauss}
 \operatorname{Scal}_\gamma
 =\operatorname{Scal}_{M,\mathrm{tan}}-|\mathrm{II}|_\gamma^2.
\end{equation}
For any compact $K\subset L$, $\lambda$ and $r_u$ are bounded on $K$ and $|R^h|_h$ is bounded on $u(K)$. If for some $M_K\geq0$,
$r_u\lambda|\mathrm{II}|_\gamma\leq M_K$ on $K\setminus Z_u$,
choose $C_K'$ with
$|\lambda^2\operatorname{Scal}_{M,\mathrm{tan}}|\leq C_K'$ on $K\setminus Z_u$.
Equation~\eqref{eq:contracted-gauss}  then gives
\[
 \lambda^2\operatorname{Scal}_\gamma
 \geq-C_K'-\frac{M_K^2}{r_u^2}
 \geq-\frac{C_K'\sup_K r_u^2+M_K^2}{r_u^2},
\]
which is the required curvature bound.
\end{proof}
\begin{corollary}[Necessary curvature growth at a critical boundary point]
\label{cor:curvature-growth}
If a smooth weakly-conformal map, $k\geq3$, has a point
$p\in Z_u\cap\overline{L\setminus Z_u}$, then
\[
 \limsup_{\substack{x\to p\\x\notin Z_u}}
 \bigl(-r_u(x)^2\lambda(x)^2\operatorname{Scal}_\gamma(x)\bigr)
 =+\infty.
\]
In particular,
\[
 \limsup_{\substack{x\to p\\x\notin Z_u}}
       r_u(x)\lambda(x)|\mathrm{II}|_\gamma(x)=+\infty.
\]
\end{corollary}
\begin{proof}
If the first limsup were not $+\infty$, there would be a neighborhood
$U$ of $p$ and a constant $C\geq0$ such that
\[
 \lambda(x)^2\operatorname{Scal}_\gamma(x)
 \geq-\frac{C}{r_u(x)^2}
 \qquad(x\in U\setminus Z_u).
\]
On a smaller compact neighborhood of $p$, chosen to be a sufficiently
small closed geodesic ball, the local argument in the proof of
Theorem~\ref{thm:intro-continuation} then gives a contradiction with Corollary~\ref{cor:uniform}.
This proves the first claim.

For the second claim, \eqref{eq:contracted-gauss} gives
\[
 -r_u^2\lambda^2\operatorname{Scal}_\gamma
 \leq (r_u\lambda|\mathrm{II}|_\gamma)^2
       -r_u^2\lambda^2\operatorname{Scal}_{M,\mathrm{tan}}.
\]
Since $|R^h|_h$ is bounded near $u(p)$,
$\operatorname{Scal}_{M,\mathrm{tan}}(x)$ is bounded for
$x\in L\setminus Z_u$ near $p$. Moreover, $r_u(x)\to0$
and $\lambda(x)\to0$ as $x\to p$. Hence
$r_u^2\lambda^2\operatorname{Scal}_{M,\mathrm{tan}}\to0$,
so the first claim implies the second.
\end{proof}

\begin{remark}
Corollary~\ref{cor:curvature-growth} says that, for every such $p$,
there exists a sequence $x_i\in L\setminus Z_u$ converging to $p$
along which both
$-r_u(x_i)^2\lambda(x_i)^2\operatorname{Scal}_\gamma(x_i)$ and
$r_u(x_i)\lambda(x_i)|\mathrm{II}|_\gamma(x_i)$ tend to $+\infty$.
It does not assert this convergence along every such sequence.
\end{remark}

\begin{remark}
Theorem~\ref{thm:intro-flatness} implies that a real-analytic
weakly-conformal map between real-analytic manifolds, with connected
domain of dimension $k\geq3$, is constant or has no critical points.

For nonconstant smooth maps in positive codimension, we cannot exclude critical points without the curvature hypothesis
of Theorem~\ref{thm:intro-continuation}. To our knowledge, this problem remains open even under the additional
assumption that the map is $k$-harmonic or a conformal $\phi$-curve.
In equal dimensions, nonconstant smooth weakly-conformal maps have no
critical points \cite[Theorem~11.4.6]{BW}.
\end{remark}

\subsection{The Bor\r{u}vka cone: a critical map with a nonplanar tangent cone}
\label{sec:boruvka}
We use the Bor\r{u}vka cone to illustrate Theorem~\ref{thm:intro-flatness}
and the curvature bound in Theorem~\ref{thm:intro-continuation}, and to
show that the image of a critical point of a conformal $\phi$-curve
need not be a branch point of its image current.
We use the sphere $\beta$, cone $\mathcal C_\beta$, and calibrated
orientation fixed in Section~\ref{sec:spherical-definitions}.
  Example~\ref{ex:boruvka-quadratic} concerns $\beta(S^2)$ and
$\mathcal C_\beta\setminus\{0\}$, whereas here we focus on what happens
at the origin.

\paragraph{The weakly-conformal parametrization.}
Write $\delta_3$ and $\delta_7$ for the Euclidean metrics on $\R^3$
and $\R^7$, and set $a:=\sqrt6$, $\alpha_0:=a-2\in(0,1)$.
For $x\ne0$, write $r:=|x|$ and $\theta:=x/|x|\in S^2$.
Define
\begin{equation}\label{eq:boruvka-parametrization}
 u:\R^3\longrightarrow\R^7,\qquad
 u(0):=0,\qquad u(r\theta):=r^a\beta(\theta).
\end{equation}
Set $R:=|u(r\theta)|=r^a$.
Since $|\beta|=1$ and $d\beta$ is bounded on $S^2$, we have
$|u(x)|=r^a$ and $|du(x)|=O(r^{a-1})$.
As $a>1$, $u$ extends to a $C^1$ map at $0$ with $du(0)=0$.

Substituting $R=r^a$ into the cone metric
$\dd R^2+6R^2g_{S^2}$ from Example~\ref{ex:boruvka-quadratic} gives
\begin{equation}\label{eq:boruvka-pullback-metric}
 u^*\delta_7
 =a^2r^{2a-2}\dd r^2+6r^{2a}g_{S^2}
 =6r^{2a-2}(\dd r^2+r^2g_{S^2}).
\end{equation}
Since $\delta_3=\dd r^2+r^2g_{S^2}$ in the coordinates $(r,\theta)$,
$u$ is weakly-conformal with dilation
$\lambda=\sqrt6\,r^{a-1}$ and critical set $Z_u=\{0\}$.
With the calibrated orientation from
Section~\ref{sec:spherical-definitions},   the submanifold $\mathcal C_\beta\setminus\{0\}$ is
associative, so $u^*\varphi_0=\lambda^3\vol_{\delta_3}$ away
from $0$. Both sides vanish at $0$, which makes $u$ a conformal
$\varphi_0$-curve in the sense of Section~\ref{sec:smith-scope}.

\paragraph{Relation to Theorem~\ref{thm:intro-flatness} and Corollary~\ref{cor:uniform}.}
Since $u(r\theta)=r^{2+\alpha_0}\beta(\theta)$, with $\beta$ smooth
and $0<\alpha_0<1$, we have $u\in C^{2,\alpha_0}$ near $0$, with
\begin{equation}\label{eq:boruvka-derivative-bounds}
 \begin{gathered}
 du(0)=D^2u(0)=0,\\
 |D^bu(x)|=O\!\left(|x|^{2+\alpha_0-b}\right)
 \qquad(b=0,1,2).
 \end{gathered}
\end{equation}
Since $r_u(x)=|x|$, these are the estimates in
Theorem~\ref{thm:intro-flatness} and Corollary~\ref{cor:uniform}
with $j=2$ and $\alpha=\alpha_0$. However, for each fixed $\theta\in S^2$,
\[
 \left|\frac{\dd^3}{\dd r^3}u(r\theta)\right|
 =a(a-1)(a-2)r^{a-3}\longrightarrow\infty
 \qquad(r\downarrow0),
\]
so $u\notin C^3$ near $0$.

\paragraph{Relation to Theorem~\ref{thm:intro-continuation} and Corollary~\ref{cor:curvature-growth}.}
We compute $\lambda^2\operatorname{Scal}_\gamma$ and
$\lambda|\mathrm{II}|_\gamma$ to show that the conclusion of
Theorem~\ref{thm:intro-continuation} can fail for $C^{2,\alpha_0}$ maps
satisfying its curvature hypothesis.

Let $\iota:\mathcal C_\beta\setminus\{0\}\hookrightarrow\R^7$ be the
inclusion and set $h_{\mathcal C}:=\iota^*\delta_7$.
Since $\gamma=u^*h_{\mathcal C}$ when $u$ is regarded as a map into
$\mathcal C_\beta\setminus\{0\}$, the map
\[
 u:(\R^3\setminus\{0\},\gamma)\longrightarrow
 (\mathcal C_\beta\setminus\{0\},h_{\mathcal C})
\]
is a local isometry. Hence, the second fundamental form $\mathrm{II}$
of $u$, viewed as an immersion into $\R^7$, is the pullback of that of
$\iota$.
The sectional curvatures of $h_{\mathcal C}$ computed at $R=1$ in
Example~\ref{ex:boruvka-quadratic} scale by $R^{-2}$ under
$Z\mapsto RZ$. Together with minimality and
\eqref{eq:minimal-contracted-gauss}, this gives
\[
 \operatorname{Scal}_\gamma=-\frac{5}{3R^2},
 \qquad |\mathrm{II}|_\gamma^2=\frac{5}{3R^2}.
\]
Substituting $R=r^a$ and $\lambda=\sqrt6\,r^{a-1}$ yields
\begin{equation}\label{eq:boruvka-curvature}
 \lambda^2|\operatorname{Scal}_\gamma|=\frac{10}{r^2},
 \qquad
 \lambda|\mathrm{II}|_\gamma=\frac{\sqrt{10}}{r}.
\end{equation}
Since $r_u=r$, these identities give
\[
 r_u^2\lambda^2\operatorname{Scal}_\gamma=-10,\qquad
 r_u\lambda|\mathrm{II}|_\gamma=\sqrt{10},
\]
which verify the curvature hypothesis of Theorem~\ref{thm:intro-continuation}.
The theorem nevertheless does not apply,
because $u$ is $C^{2,\alpha_0}$ but not even $C^3$ at the origin.
Thus a nonconstant $C^{2,\alpha_0}$ map can satisfy the curvature
bound and still have a critical point.
Moreover, $0\in Z_u\cap\overline{\R^3\setminus Z_u}$, but the two
limsups in Corollary~\ref{cor:curvature-growth} equal $10$ and
$\sqrt{10}$, respectively, rather than $+\infty$.
Hence replacing smoothness by $C^{2,\alpha_0}$ regularity would make
both Theorem~\ref{thm:intro-continuation} and
Corollary~\ref{cor:curvature-growth} false.
However, this example does not show whether the curvature bound is necessary
for smooth maps.

\paragraph{Relation to Theorem~\ref{thm:intro-branch}.}
We show that the vertex of the image current is not a branch point
in the sense of Definition~\ref{def:density-branch}.
Let $T_\beta$ be the current of integration over
$\mathcal C_\beta$, with multiplicity one and the calibrated
orientation.
By \eqref{eq:boruvka-sphere}, the sphere $\beta(S^2)$ has area
$6\cdot4\pi=24\pi$, so
\begin{equation}\label{eq:boruvka-mass}
 \|T_\beta\|(B_R(0))
 =\int_0^R\operatorname{Area}\bigl(t\beta(S^2)\bigr)\,dt
 =\int_0^R24\pi t^2\,dt=8\pi R^3.
\end{equation}
This gives locally finite mass. We now show that $\partial T_\beta=0$.
For an arbitrary smooth compactly supported $2$-form $\eta$ on $\R^7$,
we need to check that
\[
 \partial T_\beta(\eta)=T_\beta(\dd\eta)=0.
\]
Stokes' theorem on $\mathcal C_\beta\cap\{Y\in\R^7:|Y|>\epsilon\}$ gives
\begin{equation}\label{eq:boruvka-boundary-estimate}
 \left|\int_{\mathcal C_\beta\cap\{|Y|>\epsilon\}}\dd\eta\right|
 =\left|\int_{\mathcal C_\beta\cap\{|Y|=\epsilon\}}\eta\right|
 \leq C_1\epsilon^2,
\end{equation}
for some $C_1>0$.
The mass formula \eqref{eq:boruvka-mass} gives
\begin{equation}\label{eq:boruvka-inner-estimate}
 \left|\int_{\mathcal C_\beta\cap\{|Y|<\epsilon\}}\dd\eta\right|
 \leq C_2\epsilon^3,
\end{equation}
for some $C_2>0$. Both $C_1$ and $C_2$ depend on $\eta$, but are independent of $\epsilon$. Combining \eqref{eq:boruvka-boundary-estimate} and
\eqref{eq:boruvka-inner-estimate} and letting $\epsilon\downarrow0$ gives $T_\beta(\dd\eta)=0$ for every
such $\eta$, so $\partial T_\beta=0$.
Because the closed form $\varphi_0$ calibrates $T_\beta$,
it is locally area minimizing by the argument in
Section~\ref{sec:area-minimizing-currents}.
  The induced metric on $\mathcal C_\beta\setminus\{0\}$ is locally conformally flat by
\eqref{eq:boruvka-pullback-metric}.

Since $\R^7\ni Y\mapsto\rho^{-1}Y\in\R^7$ preserves the cone and its
orientation for every $\rho>0$, we have
\[
 (\eta_{0,\rho})_\#T_\beta=T_\beta\qquad(\rho>0).
\]
Its tangent cone at $0$ is therefore unique and equals
$T_\beta$. The seven component functions of $\beta$ are
linearly independent, so $\beta(S^2)$ spans $\R^7$.
Hence this tangent cone is not supported on a three-plane.
Thus $0\in\R^7$ is not a branch point under
Definition~\ref{def:density-branch}.

The mass formula \eqref{eq:boruvka-mass} also gives the density
\[
 \Theta_3(\|T_\beta\|,0)=\frac{8\pi}{4\pi/3}=6.
\]
Although this density equals $6$, the current $T_\beta$ has multiplicity one.
Since no tangent cone at $0$ is supported on a three-plane,
Theorem~\ref{thm:intro-branch} does not apply at $0$.
The example therefore shows that the image of a critical point of a conformal
$\phi$-curve need not be a branch point, in the sense of
Definition~\ref{def:density-branch}, of the current of integration over its image.

\section{Further questions}\label{sec:conclusion}
The remaining questions concern the precise limits of our results.
\begin{enumerate}[leftmargin=2em]
\item \emph{Repeated components.}
For a component with $m_j\geq2$, the corresponding currents $S_{i,a}$
from Lemma~\ref{lem:matching} converge to
$m_j\llbracket A'\times\{0\}\rrbracket$, so the application of
Allard's theorem used in Lemma~\ref{lem:extract} does not yield a
single smooth graph. It remains to understand whether
the corresponding part of the current separates into distinct graphs.
\item \emph{Planar frequency two.}
Theorem~\ref{thm:local} applies when $\nu=2$ and $k\geq4$,
but the stratification used to prove the Hausdorff-dimension bound
excludes planar frequency $2$. For $k=3$, Example~\ref{ex:boruvka-quadratic}
shows that the conclusion of Theorem~\ref{thm:local} can fail when $\nu=2$.
It remains to determine what restrictions local conformal flatness
imposes on the set of branch points of planar frequency $2$.
\item \emph{Branch points outside the rapid-decay class.}
By the results of Krummel--Wickramasekera recalled in
Section~\ref{sec:setup}, the set of branch points of $T$ outside
$\bigcup_{q=2}^{\infty}\mathcal R_q(T)$ has zero $\HH^{k-2}$ measure.
That fact alone does not imply that its Hausdorff dimension is at
most $k-3$. It remains to determine whether local conformal flatness of the
induced metric on $\reg T$ implies this stronger dimension bound for that set.
\item \emph{Currents in Riemannian manifolds.}
Proposition~\ref{prop:curved} assumes \eqref{eq:power-height-convergence}
with $1<\nu<2$ in a smooth ambient metric.
To extend Theorem~\ref{thm:intro-branch} to this setting, one would
need corresponding approximations by tangent functions and
estimates for the support of the current.
\item \emph{Critical points in positive codimension.}
It remains open whether a nonconstant smooth weakly-conformal map
with connected domain of dimension at least three in positive codimension
can have critical points without the curvature hypothesis of
Theorem~\ref{thm:intro-continuation}.
Corollary~\ref{cor:curvature-growth} implies that, for every
$p\in Z_u\cap\overline{L\setminus Z_u}$,
\[
 \limsup_{\substack{x\to p\\x\notin Z_u}}
 \bigl(-r_u(x)^2\lambda(x)^2\operatorname{Scal}_\gamma(x)\bigr)
 =+\infty,
\]
and, in particular,
\[
 \limsup_{\substack{x\to p\\x\notin Z_u}}
 r_u(x)\lambda(x)|\mathrm{II}|_\gamma(x)=+\infty.
\]
The Bor\r{u}vka example shows that replacing smoothness by
$C^{2,\alpha_0}$ regularity would make both
Theorem~\ref{thm:intro-continuation} and
Corollary~\ref{cor:curvature-growth} false.
\end{enumerate}

\phantomsection


\addcontentsline{toc}{section}{References}
\begin{thebibliography}{99}
\hypersetup{urlcolor=blue}

\bibitem{Allard}
W.~K. Allard, \emph{On the first variation of a varifold},
Ann. of Math. (2) \textbf{95} (1972), 417--491,
\href{https://doi.org/10.2307/1970868}{doi:10.2307/1970868}.

\bibitem{BW}
P.~Baird and J.~C. Wood, \emph{Harmonic Morphisms Between Riemannian
Manifolds}, London Mathematical Society Monographs, New Series, vol.~29,
Oxford University Press, 2003,
\href{https://doi.org/10.1093/acprof:oso/9780198503620.001.0001}
{doi:10.1093/acprof:oso/9780198503620.001.0001}.

\bibitem{BM}
G.~Ball and J.~Madnick, \emph{The Morse index of constant curvature
$2$-spheres}, arXiv:2606.06634v1 (4 June 2026),
  \url{https://arxiv.org/abs/2606.06634v1}.

\bibitem{Besse}
A.~L. Besse, \emph{Einstein Manifolds}, Ergebnisse der Mathematik und ihrer
Grenzgebiete, vol.~10, Springer, 1987.

\bibitem{Gallier}
J.~Gallier, \emph{Notes on Spherical Harmonics and Linear Representations of
Lie Groups}, University of Pennsylvania lecture notes, 2013,
\url{https://www.cis.upenn.edu/~cis6100/sharmonics.pdf}.

\bibitem{GT}
D.~Gilbarg and N.~S. Trudinger, \emph{Elliptic Partial Differential
Equations of Second Order}, Classics in Mathematics, Springer, 2001 reprint.

\bibitem{HL}
R.~Harvey and H.~B. Lawson, Jr., \emph{Calibrated geometries},
Acta Math. \textbf{148} (1982), 47--157,
\href{https://doi.org/10.1007/BF02392726}{doi:10.1007/BF02392726}.

\bibitem{IP}
T.~Ikonen and P.~Pankka, \emph{Liouville's theorem in calibrated geometries},
arXiv:2410.02722v1 (3 October 2024),
  \url{https://arxiv.org/abs/2410.02722v1}.

\bibitem{KWI}
B.~Krummel and N.~Wickramasekera,
\emph{Analysis of singularities of area minimizing currents, Part I:
planar frequency, branch points of rapid decay, and weak locally uniform
approximation}, arXiv:2304.10653,
\url{https://arxiv.org/abs/2304.10653}.

\bibitem{KWIII}
B.~Krummel and N.~Wickramasekera,
\emph{Analysis of singularities of area-minimizing currents, Part III:
branch points of planar frequency $\neq2$, higher order asymptotics, and the
local topology}, arXiv:2607.13356v1 (15 July 2026),
  \url{https://arxiv.org/abs/2607.13356v1}.

\bibitem{PFP}
D.~P\'erez Fern\'andez and J.~G. P\'erez,
\emph{Conformal parametrizations defined by polynomials},
Differential Geom. Appl. \textbf{89} (2023), 102009,
\href{https://doi.org/10.1016/j.difgeo.2023.102009}{doi:10.1016/j.difgeo.2023.102009}.

\bibitem{SimonGMT}
L.~Simon, \emph{Introduction to Geometric Measure Theory}, Tsinghua/NTU lecture notes,
2017--2018,
\url{https://math.stanford.edu/~lms/ntu-gmt-text.pdf}.

\bibitem{Takeuchi}
H.~Takeuchi, \emph{Some conformal properties of $p$-harmonic maps and
a regularity for sphere-valued $p$-harmonic maps},
J. Math. Soc. Japan \textbf{46} (1994), no.~2, 217--234,
\href{https://doi.org/10.2969/jmsj/04620217}{doi:10.2969/jmsj/04620217}.

\bibitem{BIM}
K.~Broder, A.~Iliashenko and J.~Madnick,
\emph{Hyperbolicity and Schwarz Lemmas in Calibrated Geometry},
arXiv:2507.16313v3 (2026),
\url{https://arxiv.org/abs/2507.16313v3}.

\bibitem{CKM}
D.~R. Cheng, S.~Karigiannis and J.~Madnick,
\emph{Bubble tree convergence of conformally cross product preserving maps},
Asian J. Math. \textbf{24} (2020), no.~6, 903--984,
\href{https://doi.org/10.4310/AJM.2020.v24.n6.a1}
{doi:10.4310/AJM.2020.v24.n6.a1}.

\bibitem{CKMVar}
D.~R. Cheng, S.~Karigiannis and J.~Madnick,
\emph{A variational characterization of calibrated submanifolds},
Calc. Var. Partial Differential Equations \textbf{62} (2023), no.~6,
Paper No.~174,
\href{https://doi.org/10.1007/s00526-023-02513-7}
{doi:10.1007/s00526-023-02513-7}.

\bibitem{IK}
A.~Iliashenko and S.~Karigiannis,
\emph{A special class of $k$-harmonic maps inducing calibrated fibrations},
Math. Res. Lett. \textbf{32} (2025), no.~5, 1481--1519,
\href{https://doi.org/10.4310/mrl.251215003006}
{doi:10.4310/mrl.251215003006}.

\bibitem{IKM}
A.~Iliashenko, S.~Karigiannis and J.~Madnick,
\emph{Montel's theorem and tautness in calibrated geometry},
arXiv:2606.31393v3 (2026),
\url{https://arxiv.org/abs/2606.31393v3}.

\bibitem{Pankka}
P.~Pankka, \emph{Quasiregular curves},
Ann. Acad. Sci. Fenn. Math. \textbf{45} (2020), no.~2, 975--990,
\href{https://doi.org/10.5186/aasfm.2020.4534}
{doi:10.5186/aasfm.2020.4534}.

\bibitem{SmithOriginal}
A.~M. Smith, \emph{A theory of multiholomorphic maps},
arXiv:1112.1471v3 (2012),
\url{https://arxiv.org/abs/1112.1471v3}.

\bibitem{DLSQ}
{\color{black}C.~De Lellis and E.~N. Spadaro,
\emph{$Q$-valued functions revisited},
Mem. Amer. Math. Soc. \textbf{211} (2011), no.~991,
\href{https://doi.org/10.1090/S0065-9266-10-00607-1}{doi:10.1090/S0065-9266-10-00607-1},
\url{https://arxiv.org/abs/0803.0060v4}.}

\bibitem{LeeRM}
J.~M.~Lee, \emph{Introduction to Riemannian Manifolds}, 2nd ed.,
Graduate Texts in Mathematics 176, Springer, Cham, 2018.

\end{thebibliography}
\end{document}